\documentclass[MSc,english,francais,nonatbib]{ulthese}
\ifxetex\else \usepackage[utf8]{inputenc} \fi

\usepackage[numbers,square]{natbib}
\usepackage{amsmath}
\usepackage{amsfonts}
\usepackage{amssymb}
\usepackage{stmaryrd} 
\usepackage[mathscr]{euscript} 
\usepackage{relsize} 
\usepackage{tikz}
\usetikzlibrary{matrix,arrows}

\usepackage{hyperref}
\hypersetup{colorlinks,allcolors=ULlinkcolor}

\newenvironment{pre}[1][Preuve :]{\begin{trivlist} \item[\hskip \labelsep {\bfseries \em #1}]\em}{\end{trivlist}}
\newcommand{\qed}{\newline ${}$ \hfill \small{Q.E.D.}}
\newcommand{\Vv}[2]{\left(\begin{array}{c}#1\\#2\end{array}\right)} 
\newcommand{\Vvv}[3]{\left(\begin{array}{c}#1\\#2\\#3\end{array}\right)} 
\newcommand{\matt}[4]{\begin{pmatrix} #1 & #2 \\ #3 & #4 \\ \end{pmatrix}} 
\newcommand{\R}{\mathcal{R}} 
\newcommand{\B}{\mathcal{B}} 
\newcommand{\A}{\mathcal{A}}
\newcommand{\Q}{\mathcal{Q}} 
\newcommand{\NN}{\mathbb{N}} 
\newcommand{\ZZ}{\mathbb{Z}} 
\newcommand{\RR}{\mathbb{R}} 
\newcommand{\CC}{\mathbb{C}} 
\newcommand{\ov}[1]{\overline{#1}} 

\newcommand{\DD}{\mathcal{D}}
\newcommand{\SSS}{\mathcal{S}}
\newcommand{\Z}{\mathcal{Z}}
\newcommand{\LL}{\mathcal{L}} 
\newcommand{\MM}{\mathcal{M}}
\newcommand{\UU}{\mathcal{U}}
\newcommand{\g}{\mathfrak{g}} 
\newcommand{\Cat}{\mathcal{C}}
\DeclareMathOperator{\HHH}{H} 
\newcommand{\cohomo}[2]{\HHH^1(#1 \; ; #2)} 
\newcommand{\gcohomo}[2]{\HHH^1\big(#1 \; ; #2\big)} 
\DeclareMathOperator{\im}{im} 
\DeclareMathOperator{\Id}{Id} 
\DeclareMathOperator{\Span}{Span} 
\DeclareMathOperator{\Hom}{Hom} 
\DeclareMathOperator{\End}{End} 
\DeclareMathOperator{\Aut}{Aut} 
\DeclareMathOperator{\Ext}{Ext} 
\DeclareMathOperator{\Der}{Der} 
\DeclareMathOperator{\IDer}{IDer} 
\DeclareMathOperator{\Max}{Max} 
\DeclareMathOperator{\supp}{supp} 
\DeclareMathOperator{\Rep}{Rep} 
\DeclareMathOperator{\Bloc}{Ext-bloc} 
\DeclareMathOperator{\ev}{ev} 
\DeclareMathOperator{\nev}{nev   } 
\DeclareMathOperator{\bl}{Bl} 
\DeclareMathOperator{\lie}{Lie} 
\DeclareMathOperator{\Imagine}{Im} 
\newcommand{\secr}{\mathlarger{\mathscr{T}}} 
\newcommand{\secb}{\mathlarger{\mathscr{S}}} 
\newcommand{\F}{\mathlarger{\mathscr{F}}} 
\newcommand{\carspec}{\mathlarger{\mathscr{CS}}} 

\newcounter{con} \numberwithin{con}{section}

\newtheorem{Def}[con]{D\'efinition}
\newtheorem{Thm}[con]{Th\'eor\`eme}
\newtheorem{Prop}[con]{Proposition}
\newtheorem{Cor}[con]{Corollaire}
\newtheorem{Lem}[con]{Lemme}
\newtheorem{Rem}[con]{Remarque}
\newtheorem{Conj}[con]{Conjecture}

\newcounter{laid} \numberwithin{laid}{chapter}
\newtheorem{annexeDef}[laid]{D\'efinition}

\newtheorem{annexeCor}[laid]{Corollaire}
\newtheorem{annexeLem}[laid]{Lemme}
\newtheorem{annexeRem}[laid]{Remarque}
\newtheorem{annexeConj}[laid]{Conjecture}

\titre{Extensions des modules de dimension finie pour les algèbres de courants tordues}
\auteur{Jean Auger}
\programme{Maîtrise en mathématiques} 
\annee{2015}

\begin{document}

\frontmatter                    

\pagetitre                      

\chapter*{Résumé}                      
\phantomsection\addcontentsline{toc}{chapter}{Résumé} 

\begin{otherlanguage*}{francais}
Ce mémoire traite de la théorie des représentations d'une certaine classe d'algèbres de Lie de dimension infinie, les algèbres de courants tordues. L'objet du travail est d'obtenir une classification des blocs d'extensions d'une catégorie de modules de dimension finie pour une algèbre de courants tordue donnée. Les principales sources de cette étude sont les récentes classifications des modules simples de dimension finie pour ces algèbres et des blocs d'extensions pour les modules de dimension finie dans le cas des algèbres d'applications équivariantes. 

Ces algèbres de courants tordues comprennent entre autres les familles d'algèbres de Lie des formes tordues et des algèbres d'applications équivariantes, donc aussi les incontournables généralisations multilacets, tordues ou non, de la théorie de Kac-Moody affine.
\end{otherlanguage*}
\chapter*{Abstract}                      
\phantomsection\addcontentsline{toc}{chapter}{Abstract} 

\begin{otherlanguage*}{english}
This master's thesis is about the representation theory of a certain class of infinite dimensional Lie algebras, the twisted current algebras. The object of this work is to obtain a classification of the extension blocks of the category of finite dimensional modules for a given twisted current algebra. The principal motivations for this study are the recent classifications of simple finite dimensional modules for these algebras and of the extension blocks of the category of finite dimensional modules in the case of equivariant map algebras. 

The class of twisted current algebras includes, amongst others, the families of Lie algebras of twisted forms and equivariant map algebras, therefore the key multiloop generalisations, twisted or not, of the affine Kac-Moody setting.
\end{otherlanguage*}
\cleardoublepage

\tableofcontents                
\cleardoublepage

%
%
\dedicace{Aux étudiants de 2012}
\cleardoublepage
\epigraphe{Le mathématique est donc la présupposition fondamentale du savoir des choses.}{Martin Heidegger}
\cleardoublepage
\chapter*{Avant-propos}         
\phantomsection\addcontentsline{toc}{chapter}{Avant-propos} 

Ce mémoire de maîtrise étant écrit, une étape notable de ma vie prend fin. Merci énormément à mon directeur de recherches Michael Lau pour avoir cru en mes capacités, pour être resté toujours positif et patient malgré cet échéancier allongé, finalement avéré nécessaire. Aussi, je lui suis très reconnaissant pour toutes ces occasions qu'il a créées ou provoquées en Amérique du Nord afin de m'aider à progresser comme apprenti mathématicien. Merci aussi au professeur Lau, pour les innombrables lettres d'appui à mes dossiers de candidatures de tous genres. Tout cela, ses contacts, ses suggestions et conseils font partie intégrante de ce qui m'a propulsé dans la situation enviable dans laquelle je me trouve en 2014-2015. J'aurai évidemment beaucoup appris à ses côté, mais j'aurai surtout gagné en assurance à l'idée de poursuivre mes études supérieures.

Un merci tout spécial au professeur Erhard Neher de l'Université d'Ottawa pour m'avoir expliqué un important détail dans une preuve de l'une de ses publications que je reprends dans ce travail. Je tiens aussi à remercier le professeur Hugo Chapdelaine de m'avoir aidé à me rendre à Paris en 2011 pour participer à un congrès international portant sur mon mathématicien favori Évariste Galois. Ce fut un séjour inspirant à plusieurs égards. Je lui suis aussi bien reconnaissant, ainsi qu'au professeur Claude Lévesque, de m'avoir occupé quelques étés de temps lors de mon passage à l'Université Laval où les quelques bourses existantes pour étudiants ne sont toujours encore que réservées aux techniciens les plus efficaces des séances d'examens. Ces marques de confiance ont été pour moi bien encourageantes. Merci à Saïd El Morchid, d'être aussi formidable, humain et pour ces fois au restaurant !

Enfin, merci à la communauté de l'anneau et aux membres du C.L.A.I.R. qui m'ont encouragé à toujours croire que j'allais bien terminer ce travail un jour. Merci aux héros de 2012 pour ces extraordinaires moments ; tous ces gens qui croient eux aussi en une éducation libre. Il est des rencontres que j'y ai faites qui m'auront assurément changé : tantôt peu, tantôt démesurément. Vive la neige et la pizza ! Je salue mon existante et unique projection ainsi que Mlle D., des personnes très généreuses que je suis content de pouvoir connaître. En terminant, je salue tout spécialement les autres mousquetaires, au nombre inhabituel de trois, avec qui j'aurai vécu non-moins d'un cinquième de mon existence avant de les quitter, eux et tous les autres, pour le Mordor. Je suis d'ailleurs bien reconnaissant à P.J. de m'y avoir accueilli si naturellement.           

\mainmatter                     

\chapter*{Introduction}         
\phantomsection\addcontentsline{toc}{chapter}{Introduction} 

La théorie de Lie de dimension infinie est vaste et se développe pas à pas. Interprétée comme théorie qui doit décrire des symétries, l'étude de son volet de théorie des représentations améliore toujours la compréhension que l'on peut avoir des objets rattachés à un type donnée de symétries. Les représentations dites d'évaluation sont un outil théorique qui s'est révélé bien utile depuis son introduction dans la fin des années 80. Ces notions ont permis d'améliorer sensiblement la compréhension de plusieurs types de représentations d'algèbres de Lie. Par exemple, ces outils ont contribué à comprendre les représentations des incontournables algèbres de multilacets, tordues ou pas, qui généralisent des éléments de la théorie des algèbres de Kac-Moody affines. Récemment, les modules simples de dimension finie des algèbres de Lie qui sont des familles des « formes tordues », des « algèbres d'applications équivariantes » et des « algèbres de courants tordues » ont été classifiés via ces notions.

En théorie des représentations, il est utile de voir des représentations partageant des mêmes propriétés en tant qu'un tout ; une catégorie. Chercher à décrire les catégories de modules en tant qu'un tout est généralement un des premiers objectifs des mathématiciens qui étudient le sujet. Cela permet d'avoir une première idée générale des modules considérés et des liens entre ceux-ci. C'est dans cet esprit qu'apparaît ce travail de maîtrise qui traite d'un aspect des catégories de modules de dimension finie pour les « algèbres de courants tordues », ses blocs d'extensions.

Les algèbres de courants tordues forment une classe d'algèbre de Lie qui généralise, en quelque sorte, les formes tordue et les algèbres d'applications équivariantes. À l'intérieur de la catégorie des modules de dimension finie pour une algèbre de courants tordue donnée, tous les objets ne sont pas semi-simples. Cela signifie que les notions d'extensions, puis de blocs d'extensions sont à la fois pertinentes et d'intérêt. Une description des blocs d'extensions doit mener à une meilleure compréhension de l'agencement des modules en question dans la catégorie. Ce mémoire vise à donner une description et une classification de la décomposition en blocs d'extensions de la catégorie des modules de dimension finie pour une algèbre de courants tordue donnée. Cette étude prend ses assises sur, d'une part, la classification des objets simples de cette catégorie, mais surtout sur l'étude analogue faite dans le cas plus particulier des algèbres d'applications équivariantes par les mathématiciens E.Neher et A.Savage.

Dans le premier chapitre sont présentés les éléments essentiels d'algèbre homologique et de théorie des catégories pour comprendre la problématique de l'étude des blocs d'extensions. Le second chapitre développe d'abord des éléments de théorie des représentations et de cohomologie des algèbres de Lie qui s'avéreront nécessaires par après. Ensuite, la classification des objets simples d'une catégorie de modules de dimension finie pour une algèbre de courants tordue est rapidement abordée. Ce même chapitre se termine avec la classification des blocs d'extensions de ce type de catégorie. Le troisième chapitre présente les résultats de classification des blocs d'extensions dans deux contextes quelque peu plus spécifiques : les algèbres d'applications équivariantes et les formes tordues. Pour terminer, il y est aussi question de détailler la description d'une algèbre de Margaux et de ses modules de dimension finie. Une algèbre de Margaux est une forme tordue qui n'est pas une algèbre de multilacets tordue. Ce type d'algèbre n'est encore que bien peu compris par la communauté mathématique.

\vspace{3ex}
\textbf{\large{Notation}}

Tout au long du document, les notations suivantes seront employées.

\begin{itemize}
\item[\textbullet$\;$] $\NN = \{0,1,2,3,...\}$.
\item[\textbullet$\;$] $k$ est un corps algébriquement clos et de caractéristique 0.
\item[\textbullet$\;$] $\LL$ est une algèbre de courants tordue.
\item[\textbullet$\;$] $\LL$-$\mathbf{mod}$ est la catégorie des $\LL$-modules qui sont de dimension finie.
\item[\textbullet$\;$] $\Gamma$ est un groupe fini.
\item[\textbullet$\;$] $\g$ est une $k$-algèbre de Lie simple de dimension finie.
\item[\textbullet$\;$] $S$ est une $k$-algèbre de type fini qui est associative, commutative, unitaire et réduite. Son nilradical est nul, ce qui implique que son radical de Jacobson est nul dans ce contexte.
\item[\textbullet$\;$] $s(M) \in k \cong S / M$, où $s \in S \text{ et } M \in \Max S)$, est la classe de résidu de $s$ selon $M$.
\item[\textbullet$\;$] $\Gamma^x$ est le stabilisateur de $x \in X$ où une action de groupe $\Gamma \curvearrowright X$ est sous-entendue.
\item[\textbullet$\;$] $A^\Gamma$ est le sous-module des $\Gamma$-invariants dans le $\Gamma$-module $A$.
\item[\textbullet$\;$] $k_\lambda$ est le module de dimension 1 d'une $k$-algèbre de Lie $L$ dont l'espace vectoriel est $k$ et dont l'action est donnée par $\lambda \in (L/[L,L])^*$.
\item[\textbullet$\;$] $L'$ est l'algèbre dérivée $[L,L]$ d'une algèbre de Lie $L$.
\item[\textbullet$\;$] $V^{\,\oplus\, n} = V \oplus \cdots \oplus V$ ($n$ fois) pour $n \in \NN$ et pour tout espace vectoriel ou module $V$.
\item[\textbullet$\;$] $[V]$ est la classe d'isomorphismes du module $V$ dans la catégorie spécifiée.
\item[\textbullet$\;$] $\llbracket V \rrbracket$ est le bloc d'extensions qui contient le $\LL$-module $V$ dans la catégorie spécifiée.
\end{itemize}

À moins d'un avis contraire, tous les produits tensoriels $\,\otimes\,$ signifieront $\,\otimes_k\,$, le produit tensoriel pris sur $k$. De plus, tous les $\LL$-modules génériques considérés seront supposés de dimension finie sur $k$.

\chapter{Algèbre homologique et catégories}     

Dans ce chapitre sont présentés le cadre général et la logique dans laquelle s'inscrit ce travail de maîtrise. On y retrouve d'abord les éléments qui donnent une raison d'être à l'étude des extensions d'objets mathématiques tels les modules. Il est ensuite question de quelques outils et éléments théoriques qu'il faut saisir pour bien pouvoir analyser les extensions dans une catégorie. Ce chapitre se termine avec une présentation de quelques propriétés capitales des notions de blocs et de blocs d'extensions ainsi qu'avec un énoncé des objectifs de ce mémoire, alors placés en un contexte approprié.

\section{Analyse d'une catégorie abélienne}
\subsection{Démarche globale et généralités}

Le but de ce travail est d'améliorer la compréhension de la catégorie des modules de dimension finie pour une algèbre de courants tordue donnée. Une algèbre de courants tordue est une algèbre de Lie obtenue comme une algèbre de points fixes par l'action d'un groupe fini par $k$-automorphismes sur un produit tensoriel de la forme $\g \otimes S$, où $\g$ est une algèbre de Lie simple de dimension finie et $S$ est une $k$-algèbre de type fini, associative, commutative, unitaire et réduite. 

La définition des algèbres de courants tordues est rappelé avec un plus grand soin au début du chapitre 2. Ce qu'il est important de savoir pour aborder ce document, c'est que les catégories considérées seront toutes au moins abéliennes. 

Une caractéristique importante des catégories abéliennes est que tous les ensembles de morphismes ont des structures de groupes abéliens et que les compositions de morphismes sont des applications au moins $\ZZ$-bilinéaires. Un autre élément clé dans la définition des catégories abéliennes est que pour toute flèche entre une paire d'objets de la catégorie, il y a une notion de noyau, de co-noyau et d'image. Lorsqu'il est possible de définir une notion naturelle de quotient, les théorèmes d'isomorphisme usuels de la théorie des groupes ont habituellement des analogues valables. C'est le cas pour les catégories de modules et en particulier, pour la catégorie des modules d'une algèbre de courants tordue.

\begin{Rem} Une définition précise de ce qu'est une catégorie abélienne peut être trouvée dans plusieurs ouvrages de référence en algèbre. Par exemple, la section 7.7 de \cite{Etingof} en traite très brièvement, mais suffisamment pour le cadre de ces explications.
\end{Rem}

Dans la catégorie des modules pour une algèbre de courants tordue, les notions de produits directs, de sommes directes et de produits tensoriels d'objets existent et sont bien définies. Il s'avère que ce produit tensoriel est une opération associative à isomorphisme près qui possède aussi un « élément neutre » pour cette opération ; le module trivial de dimension 1. Ainsi, à isomorphisme près, les objets d'une catégorie de modules forment un monoïde. Une catégorie avec ce type de propriété est parfois appelée une catégorie tensorielle ou une catégorie monoïdale.

\begin{Rem} Pour aborder le reste du document, il sera bon de garder en tête le contexte d'une catégorie de modules, c'est-à-dire d'une catégorie (abélienne) monoïdale ou tensorielle.
\end{Rem}

Soit $A$ un anneau ou une algèbre, soient $M_1$ et $M_2$ deux $A$-modules et soit $f \in \Hom_A(M_1,M_2)$. Alors on définit naturellement :
\begin{align*}
\ker f &= \{m_1 \in M_1 \; | \; f(m_1) = 0 \in M_2\} \subseteq M_1\\
\im f &= \{f(m_1) \; | \; m_1 \in M_1\} \subseteq M_2
\end{align*}
Ce sont des sous-modules de $M_1$ et $M_2$ respectivement, c'est-à-dire des ensembles stables sous l'action de $A$. On a alors toujours un isomorphisme de modules
\begin{align*}
M_1 / \ker f \quad &\cong \quad \im f  \quad \subseteq \; M_2\\
x + \ker f \quad &\leftrightarrow \;\; f(x)
\end{align*}
où l'action de $A$ sur le module $M_1 / \ker f$ est donnée par $a.\big(x + \ker f\big) = (a.x) + \ker f$ pour tout $a \in A$ et $x \in M_1$. Ce dernier résultat est justement le premier théorème d'isomorphisme qui demeure fondamental dans toute la théorie des représentations. Les deux autres théorèmes d'isomorphisme sont également bien utiles, mais ils ne sont pas nécessaires pour le moment.

Soit $\LL$ une algèbre de courants tordue et soit $\LL$-$\mathbf{mod}$ la catégorie des $\LL$-modules de dimension finie. Tel que mentionné dans l'introduction, le but de ce travail est de décrire la décomposition de $\LL$-$\mathbf{mod}$ en ce qu'on appellera des blocs. Les blocs sont essentiellement des « sous-catégories pleines » dont la « somme  directe » est toute la catégorie $\LL$-$\mathbf{mod}$. Par sa définition même, cette décomposition permet de mieux comprendre la catégorie à l'étude. Il sera question des blocs à la section 1.3.1 de ce chapitre.

\begin{Def} Une sous catégorie $\mathcal{S}$ d'une catégorie $\Cat$ est dite \textbf{pleine} si pour chaque paire d'objets $A$ et $B$ de $\mathcal{S}$, on ait
$$ \Hom_\mathcal{S}(A,B) = \Hom_\Cat(A,B) $$
\end{Def}

\begin{Rem} L'expression « somme directe » du dernier paragraphe fait référence à la propriété que tout objet d'une catégorie qui a des blocs se décompose de façon unique (à isomorphismes près), comme une somme directe d'objets qui appartiennent chacun à un bloc distinct. 
\end{Rem}

\subsection{Suites exactes, foncteurs et principes de base}

Soit $L$ une algèbre de Lie. Pour la suite de cette sous-section, il sera question de la catégorie des $L$-modules. 

Des outils très pratiques ont été développés en algèbre homologique pour l'étude de nombreux concepts dont les $L$-modules. La notion de base, qui est la pierre angulaire de cette théorie, est celle d'une suite exacte de modules. Ceci mérite quelques définitions en bonnes et dues formes.

\begin{Def} Soit $\{M_i\}_{i \in \ZZ}$ une collection de $L$-modules et soient $\{f_i : M_i \rightarrow M_{i+1}\}_{i \in \ZZ}$ des homomorphismes de $L$-modules. On peut alors former une suite
$$ \cdots \; \stackrel{f_{i-2}}{\longrightarrow} \; M_{i-1} \; \stackrel{f_{i-1}}{\longrightarrow} \; M_i \; \stackrel{f_i}{\longrightarrow} \; M_{i+1} \; \stackrel{f_{i+1}}{\longrightarrow} \; M_{i+2} \; \stackrel{f_{i+2}}{\longrightarrow} \; \cdots $$
La suite ci-dessus est dite \textbf{exacte en $M_i$} si $\, \im f_{i-1} =  \ker f_i \subseteq M_i$.

La même suite est dite \textbf{exacte} si elle l'est en chacun des $M_i$ qui la compose.
\end{Def}

\begin{Def} Une \textbf{suite exacte courte} de $L$-modules est une suite exacte de la forme
$$ 0 \longrightarrow X \stackrel{f}{\longrightarrow} Y \stackrel{g}{\longrightarrow} Z \longrightarrow 0 $$
où $X$, $Y$ et $Z$ sont des $A$-modules, $f \in \Hom_L(X,Y)$ et $g \in \Hom_L(Y,Z)$.
\end{Def}

\begin{Rem} Qu'une suite $0 \longrightarrow X \stackrel{f}{\longrightarrow} Y \stackrel{g}{\longrightarrow} Z \longrightarrow 0$ soit exacte est la même chose que d'avoir $f$ injective, $g$ surjective ainsi qu'un isomorphisme $Y/X \cong Z$ donné par l'exactitude de la suite en $X$ et le premier théorème d'isomorphisme.
\end{Rem}

\begin{Rem} Soient $Y$ et $Z$ des $L$-modules et soit $g \in \Hom_L(Y,Z)$. Alors le premier théorème d'isomorphisme donne une suite exacte courte $\; 0 \rightarrow \ker g \rightarrow Y \stackrel{g}{\longrightarrow} \im g \rightarrow 0$ où la première application non-triviale est simplement l'inclusion.
\end{Rem}

Avant de passer à autre chose, il faut rappeler la notion de suite exacte scindée. Cette notion se révélera importante lorsque viendra le temps de traiter des extensions de modules.

\begin{Def} Soient $X$, $Y$ et $Z$ des $L$-modules et soit $f \in \Hom_L(X,Y)$ ainsi que $g \in \Hom_L(Y,Z)$ tels qu'on ait une suite exacte courte de $L$-modules
$$ 0 \longrightarrow X \stackrel{f}{\longrightarrow} Y \stackrel{g}{\longrightarrow} Z \longrightarrow 0 $$
Cette suite sera dite \textbf{scindée} s'il existe $s \in \Hom_L(Z,Y)$ tel que $g \circ s = \Id_Z$. Un tel homomorphisme $s$ est parfois appelé \textbf{une section}.
\end{Def}

\begin{Rem} \label{suitescindee} Il y a quelques autres définitions équivalentes d'une suite exacte courte scindée dans ce contexte. Voir par exemple les équivalences de l'exercice 37 de \cite{FDnoncom}. 

Il est même possible de préciser que si une suite exacte courte de modules est scindée, alors la section induit un isomorphisme de $L$-modules
$$ Y = \ker g \oplus \im s = f(X) \oplus s(Z) \cong X \oplus Z $$
Pour le voir, il suffit d'écrire tout $y \in Y$ comme $y = \Big(y-s\big(g(y)\big)\Big) + s\big(g(y)\big)$. 
\end{Rem}

Un autre aspect important de l'algèbre homologique est l'utilisation des foncteurs. Les foncteurs jouent le rôle des morphismes, mais entre des catégories plutôt qu'entre des objets d'une catégorie donnée. Dans toute catégorie, deux foncteurs ont des rôles plus prédominants ; ceux donnés par les morphismes vers ou dans un objet préalablement fixé. Si le contexte est celui d'une étude de modules pour $L$, alors on verra au chapitre 2 qu'en fixant un module $X$, les deux foncteurs
\begin{align*}
\Hom_L(X,-) && \Hom_L(-,X)
\end{align*}
sont de la catégorie des $L$-modules dans elle-même.  

Étant donné des foncteurs d'une catégorie dans elle-même, on peut les appliquer à des suites de morphismes et d'objets pour obtenir de nouvelles suites. Des questions de base en algèbre homologique visent à décrire l'effet de tels foncteurs (ou de foncteurs plus généraux) sur les suites exactes et sur les suites exactes courtes. Plusieurs résultats permettent parfois de déduire des informations précieuses sur les objets d'une suite exacte à laquelle on aurait appliqué les foncteurs. En particulier, cette analyse dans le cas des foncteurs $\Hom_L(X,-)$ et $\Hom_L(-,X)$ mène ultimement à plusieurs résultats bien pratiques. Pour des détails concernant ce genre de considérations, des livres de référence tels \cite{Rotman} et \cite{Weibel} sont utiles. 

Hormis à quelques reprises dans des preuves de la fin du chapitre 2, ce dont il est question dans le paragraphe précédent n'est pas d'intérêt pour le but de ce travail.

\section{Théorème de Jordan-Hölder}
\subsection{Suites de composition et énoncé du théorème}

Supposons de nouveau que $L$ soit une $k$-algèbre de Lie et que $L$-$\mathbf{mod}$ soit la catégorie des $L$-modules de dimension finie. Dans ce contexte, si $V$ est un $L$-module, soit il admet des sous-modules non triviaux ou pas. S'il n'en admet pas, c'est que c'est un module simple et on peut écrire $\{0\} \subseteq V$. Si par contre il en possède au moins un et que l'on note $V_1$ l'un de ces sous-modules, alors $\dim_k V_1 < \dim_k V$. En se posant la même question avec $V_1$, on obtient ou non un sous-module non trivial $V_2 \subseteq V_1$ avec $\dim_k V_2 < \dim_k V_1$. Puisque la dimension des $V_i$ décroît strictement, il est évident que ce processus a toujours une fin. Ainsi, à partir du module de dimension finie $V$, on obtient une suite de sous-modules emboîtés
$$ \{0\} = V_{N+1} \subseteq V_N \subseteq V_{N-1} \subseteq \cdots \subseteq V_1 \subseteq V_0 = V $$
où $N \in \NN$ est un certain entier naturel. Évidemment, on doit avoir $\dim_k V_{i+1} < \dim_k V_i$ pour chaque $i \in \{0,\,..., N\}$.

Il se trouve qu'à chaque étape de la construction du dernier paragraphe, si le module n'est pas simple, il est toujours possible de choisir un sous-module de dimension maximale. Cela vient du fait que tous les modules sont de dimension finie. Si on construit une suite de sous-modules emboîtés de $V$ avec cette attention toute particulière, il est facile de montrer que pour chaque $i \in \{0,\,...,N\}$, le module quotient $V_i/V_{i+1}$ est un module simple. Ce genre de suite fait l'objet d'une attention particulière dans la théorie et mérite donc une définition.

\begin{Def} Soit $V$ un $L$-module de dimension finie. Alors une \textbf{suite de composition pour $V$} est une suite de sous-modules emboîtés 
$$ \{0\} \subseteq V_N \subseteq V_{N-1} \subseteq \cdots \subseteq V_1 \subseteq V_0 = V $$
avec les modules quotients $V_i/V_{i+1}$ simples pour chaque $i \in \{0,\,...,N\}$. 
\end{Def}

Le théorème de Jordan-Hölder, qui est un résultat extrêmement important, explique dans plusieurs contextes, que les notions de suites de composition renferment des informations invariantes (à équivalences près) sur les objets desquels ils sont issus. Dans le cas des $L$-modules de dimension finie, ce théorème se lit comme suit.

\begin{Thm} \em \textbf{(Jordan-Hölder)} \label{jordanholder} \em Soit $V$ un $L$-module de dimension finie et soient 
\begin{align*}
&\{0\} = V_{N+1} \subseteq V_N \subseteq V_{N-1} \subseteq \cdots \subseteq V_1 \subseteq V_0 = V\\
&\{0\} = \tilde{V}_{M+1} \subseteq \tilde{V}_N \subseteq \tilde{V}_{N-1} \subseteq \cdots \subseteq \tilde{V}_1 \subseteq \tilde{V}_0 = V
\end{align*}
deux suites de composition de $V$. Alors 

\begin{itemize}
\item[\textbullet] $N = M$
\item[\textbullet] $\big\{[V_i/V_{i+1}]\big\}_{i=0}^N = \big\{[\tilde{V}_i/\tilde{V}_{i+1}]\big\}_{i=0}^M$
\item[\textbullet] chaque classe d'isomorphismes de modules simples des ensembles $\big\{[V_i/V_{i+1}]\big\}_{i=0}^N$ et \linebreak $\big\{[\tilde{V}_i/\tilde{V}_{i+1}]\big\}_{i=0}^M$ apparaît un même nombre de fois dans chacune des deux suites de composition.
\end{itemize}
\end{Thm}

Pour une preuve, voir par exemple le théorème 3.7.1 de \cite{Etingof} et les explications données sous le théorème 0.5 de \cite{FDnoncom}. À noter que des modules pour une algèbre de Lie sont des modules pour un anneau via la notion d'algèbre universelle enveloppante.

\begin{Rem} En particulier, le théorème de Jordan-Hölder permet d'associer à chaque $L$-module de dimension finie un nombre fini de classes d'isomorphismes de $L$-modules simples. Un nom sera donné à ces modules.
\end{Rem}

\begin{Def} Soit $V$ un $L$-module de dimension finie et soit une suite de composition pour ce dernier 
$$ \{0\} = V_{N+1} \subseteq V_N \subseteq V_{N-1} \subseteq \cdots \subseteq V_1 \subseteq V_0 = V $$
On appelle les modules simples $\{V_i/V_{i+1}\}_{i=0}^N$ les \textbf{facteurs de V}. Pour un indice $i$ fixé, le nombre de facteurs de $V$ isomorphes à $V_i/V_{i+1}$ est appelé \textbf{la multiplicité de $V_i/V_{i+1}$} dans $V$.
\end{Def}

\subsection{Extensions}

Dans les cas où il est applicable, le théorème de Jordan-Hölder permet de décomposer un objet donné en un ensemble fini de plus petits objets. La question de savoir comment composer de gros objets à l'aide de plus petits objets apparaît donc naturellement. Une bonne connaissance théorique de ces deux problématiques devrait permettre de faciliter l'étude de la catégorie d'intérêt. 

Une notion qui traduit bien l'idée de composer de gros objets à partir de petits objets est celle d'extension. Pour toute la suite de cette section, il sera sous-entendu que la catégorie utilisée sera celle des $L$-modules où $L$ est une $k$-algèbre de Lie.

\begin{Def} Soient $V$ et $W$ deux $L$-modules. Alors \textbf{une extension de $V$ par $W$} est la donnée d'un $L$-module $E$ et d'une suite exacte courte 
$$ 0 \rightarrow W \rightarrow E \rightarrow V \rightarrow 0 $$
\end{Def}

\begin{Rem} L'ensemble des extensions d'une paire arbitraire de $L$-modules n'est jamais vide. En effet, la somme directe des deux modules $V \oplus W$, l'inclusion $W \hookrightarrow V \oplus W$ et la projection $V \oplus W \twoheadrightarrow V$ fournit toujours une suite exacte courte.
\end{Rem}

\begin{Rem} \label{avantequivalence} Une extension n'est pas uniquement la donnée d'un $L$-module $E$. Supposons que $E$ soit une extension de $V$ par $W$ via la suite exacte courte 
$$ 0 \rightarrow W \stackrel{f}{\longrightarrow} E \stackrel{g}{\longrightarrow} V \rightarrow 0 $$ 
Alors le même module $E$ est aussi une extension de $V$ par $W$ via la suite exacte courte
$$ 0 \rightarrow W \stackrel{\frac{1}{\lambda} f}{\longrightarrow} E \stackrel{\lambda g}{\longrightarrow} V \rightarrow 0 $$
où $\lambda \in k \bs \{0,1\}$. Manifestement, ce sont là deux extensions de $V$ par $W$ qui sont, a priori, tout à fait distinctes.
\end{Rem}

Naturellement, il faut avoir une notion d'équivalence pour des extensions. La remarque précédente montre que de baser cette notion d'équivalence uniquement sur la classe d'isomorphismes du $L$-module dans la donnée d'une extension n'est pas assez précis. 

Voici comment est définie la relation d'équivalence pour des extensions d'une paire de modules donnée. 

\begin{Def} \label{GROS} Soient $V$ et $W$ deux $L$-modules et soient des extensions de $V$ par $W$ données par les suites exactes courtes
\begin{align*}
&0 \rightarrow W \stackrel{f_1}{\longrightarrow} E_1 \stackrel{g_1}{\longrightarrow} V \rightarrow 0\\
&0 \rightarrow W \stackrel{f_2}{\longrightarrow} E_2 \stackrel{g_2}{\longrightarrow} V \rightarrow 0
\end{align*}
L'extension 1 est dite équivalente à l'extension 2 s'il existe un isomorphisme de $L$-modules $\varphi : E_1 \rightarrow E_2$ tel que le diagramme suivant soit commutatif
\begin{center}\begin{tikzpicture}
y\matrix (m) [matrix of math nodes, row sep=0.5em,column sep=2.5em,text height=1.5ex,text depth=0.25ex, nodes in empty cells]
{ & & E_1 & & \\ 0 & W & & V & 0\\ & & E_2 & & \\};
\draw (m-2-1) edge[->,thick] (m-2-2)
          (m-2-2) edge[->,thick] node[above] {\footnotesize{$f_1$}} (m-1-3)
          (m-2-2) edge[->,thick] node[below] {\footnotesize{$f_2$}} (m-3-3)
          (m-1-3) edge[->,thick] node[auto] {\footnotesize{$\varphi$}} (m-3-3)
          (m-1-3) edge[->,thick] node[above] {\footnotesize{$g_1$}} (m-2-4)
          (m-3-3) edge[->,thick] node[below] {\footnotesize{$g_2$}} (m-2-4)
          (m-2-4) edge[->,thick] (m-2-5);
\end{tikzpicture} \end{center}

L'ensemble des classes d'équivalences d'extensions de $V$ par $W$ est noté $\Ext_L^1(V,W)$.
\end{Def}

\begin{Rem} Avec cette notion d'équivalence, il est facile de montrer que les deux extensions de la remarque \ref{avantequivalence} sont équivalentes. Pour le voir, il suffit de considérer $\varphi = \frac{1}{\lambda}\Id_E$.
\end{Rem}

\begin{Rem} Soit $E$ une extension de $V$ par $W$ via une suite exacte courte
$$ 0 \rightarrow W \rightarrow E \rightarrow V \rightarrow 0 $$
La remarque précédente laisse toujours sans réponse la question de savoir à quel nombre de classes d'extensions distinctes peut correspondre la classe d'isomorphismes $[E]$ du $L$-module $E$. Dans le cas où $V$ et $W$ sont tous deux des modules simples et de dimension finie, une réponse à cette question est avancée par l'exercice (d) du problème 3.9.1 de \cite{Etingof}. La question plus générale ne faisant pas partie des objectifs de ce travail, je ne développerai pas davantage sur la question.
\end{Rem}

\begin{Def} Soient $V$ et $W$ des $L$-modules. Alors une extension de $V$ par $W$ est dite \textbf{triviale} si elle est donnée par une suite exacte courte qui est scindée.
\end{Def}

\begin{Rem} \label{zeroext} Une extension triviale $0 \rightarrow W \rightarrow E \rightarrow V \rightarrow 0$ est équivalente à l'extension donnée par le module $W \oplus V$ via la suite exacte
$$ 0 \rightarrow W \stackrel{\iota}{\longrightarrow} W \oplus V \stackrel{\pi}{\longrightarrow} V \rightarrow 0 $$
où $\iota$ est l'inclusion canonique et où $\pi$ est la projection canonique.

Pour le voir, il suffit de considérer la décomposition du $L$-module $E$ comme $V \oplus W$ d'après l'isomorphisme décrit à la remarque \ref{suitescindee}. Ensuite, on vérifie que l'isomorphisme $\varphi$ entre $E$ et $V \oplus W$ qui identifie les parties correspondantes fait effectivement commuter le diagramme pour l'équivalence des extensions (voir la définition \ref{GROS}).
\end{Rem}

\begin{Rem} Si on fixe $V$ et $W$ des objets de la catégorie $L$-$\mathbf{mod}$, alors la classe d'équivalences de l'extension triviale de $V$ par $W$ est notée $0 \in \Ext_L^1(V,W)$.

Ainsi, si toutes les extensions de $V$ par $W$ correspondent à des suites exactes courtes scindées, on peut écrire $\Ext_L^1(V,W) = \{0\}$.
\end{Rem}

Pour justifier encore une fois l'importance de la notion d'extension, sont immédiatement décrits quelques constats pouvant être faits lorsque Jordan-Hölder s'applique à une catégorie et qu'on considère les facteurs d'un objet donné.

Supposons que $L$ soit une $k$-algèbre de Lie et soit $L$-$\mathbf{mod}$ la catégorie des $L$-modules de dimension finie. Fixons $V$ un objet de $L$-$\mathbf{mod}$. Supposons de plus que $V$ ait une suite de composition
$$ \{0\} = V_{N+1} \subseteq V_N \subseteq V_{N-1} \subseteq \cdots \subseteq V_1 \subseteq V_0 = V $$
Si le théorème de Jordan-Hölder s'applique à $L$-$\mathbf{mod}$, alors $V$ peut s'écrire comme une suite d'extensions successives faisant intervenir ses facteurs $\{C_i = V_i/V_{i+1}\}_{i=0}^N$. En effet, on peut écrire des suites exactes courtes
\begin{align*}
0 \rightarrow 0 \rightarrow V_N \rightarrow \; &C_N \rightarrow 0\\
0 \rightarrow V_N \rightarrow V_{N-1} \rightarrow \; &C_{N-1} \rightarrow 0\\
0 \rightarrow V_{N-1} \rightarrow V_{N-2} \rightarrow \; &C_{N-2} \rightarrow 0\\
&\;\,\vdots\\
0 \rightarrow V_{2} \rightarrow V_{1} \rightarrow \; &C_1 \rightarrow 0\\
0 \rightarrow V_1 \rightarrow V \rightarrow \; &C_0 \rightarrow 0
\end{align*}

En un sens, ceci montre qu'une connaissance des extensions entre des modules est au moins souhaitable. Vers la fin du chapitre 2, il sera question des classes d'équivalences d'extensions entre deux modules simples arbitraires pour une algèbre de courants tordue. Cela constituera un premier pas dans une compréhension plus globale de cette catégorie de modules.

\section{Blocs et blocs d'extensions}
\subsection{Définitions}
Soit $V$ un $L$-module de dimension finie, c'est-à-dire un objet de $L$-$\mathbf{mod}$. Le théorème de Jordan-Hölder et les dernières remarques de la section 1.2.2 sur les extensions expliquent que $V$ peut être assimilé à deux choses :
\begin{itemize}
\item[\textbullet] Un nombre bien précis de facteurs qui, à isomorphismes près, se répartissent en un nombre déterminé de classes d'isomorphismes qui apparaissent chacune un nombre déterminé de fois.
\item[\textbullet] Des extensions successives ne dépendant que d'une suite de composition de $V$ ; ces extensions sont déterminées par $V$.
\end{itemize}

En théorie, ces informations doivent être équivalentes à l'information contenue dans la classe d'isomorphismes $[V]$ du $L$-module $V$. Pour mieux comprendre les modules de $L$-$\mathbf{mod}$, une relation d'équivalence en lien avec les extensions est introduite.

\begin{Def} \label{defdebasebloc} Soient $V$ et $W$ des objets simples de $L$-$\mathbf{mod}$. Ces deux objets sont dit \textbf{liés} s'il existe une suite finie de $L$-modules $\{T^i\}_{i=0}^N$ avec $T^0 = V$, $T^N = W$ et
$$ \Ext_L^1(T^i,T^{i+1}) \neq \{0\} \qquad \qquad \text{ou} \qquad \qquad  \Ext_L^1(T^{i+1},T^i) \neq \{0\} $$ 
pour chaque $i \in \{0,\,..., N-1\}$.

Il est entendu que toute paire de modules simples isomorphes est liée par une suite de modules de longueur 0.

Les classes d'équivalences de cette relation sont appelées les blocs d'extensions de $L$-$\mathbf{mod}$, pour l'algèbre de Lie $L$. La notation $\Bloc(L)$ est parfois employée pour signifier l'ensemble des blocs d'extensions de l'algèbre de Lie $L$.
\end{Def}

Ceci est une relation d'équivalence sur l'ensemble des objets simples de $L$-$\mathbf{mod}$. D'après la définition même de cette relation d'équivalence, la connaissance des classes d'équivalences auxquelles elle donne lieu fournit plusieurs informations sur les extensions entre les modules simples.

\begin{Rem} En anglais, deux objets simples liés (selon la relation d'équivalence ci-dessus) sont dits « linked » et la suite des modules $T^i$ qui fait le lien entre les deux est parfois appelée une chaîne de modules. Ceci fait penser immédiatement aux $T^i$ comme à des maillons d'une chaîne ou « chain links ».
\end{Rem}

Par le théorème de Jordan-Hölder, on peut regrouper des modules de $L$-$\mathbf{mod}$ selon ce que leurs facteurs sont tous « liés » les uns aux autres. Faire de tels amalgames donne espoir qu'on puisse mieux comprendre les extensions entre des modules arbitraires de $L$-$\mathbf{mod}$. C'est de cette façon que la notion de blocs de la catégorie $L$-$\mathbf{mod}$ est introduite.

\begin{Def} Soit $b \in \Bloc(L)$. Un \textbf{bloc de la catégorie $L$-$\mathbf{mod}$}  est le sous-ensemble de modules
$$ \Cat_b \;\; = \;\; \left\{\;\; V \in L\text{-}\mathbf{mod} \quad \Bigg| \; \begin{array}{c} \text{Tous les facteurs de } V \text{ sont}\\
\text{dans le bloc d'extensions } b \end{array}\right\} $$
En particulier, les blocs de la catégorie $L$-$\mathbf{mod}$ sont en bijection avec l'ensemble $\Bloc(L)$.
\end{Def}

\begin{Prop} Chacun des blocs de la catégorie $L$-$\mathbf{mod}$ est naturellement une catégorie. Un bloc est en fait une sous-catégorie pleine de $L$-$\mathbf{mod}$.
\end{Prop}

Cette proposition est la proposition 5.2 de l'article \cite{NSext}. Le contenu de cette proposition apparaît aussi dans \cite{Etingof} au début de la section 9.5, mais dans le contexte des modules pour une algèbre de dimension finie.

Tel que mentionné précédemment, une sous-catégorie pleine signifie simplement que les espaces d'homomorphismes entre des objets de la sous-catégorie sont les mêmes peu importe que l'on voie ces objets comme issus de la sous-catégorie ou non. Concrètement, si $V$ et $W$ sont des objets de $\Cat_b$ pour un certain $b \in \Bloc(L)$, alors $\Hom_{\Cat_b}(V,W) = \Hom_L(V,W)$.

\subsection{Propriétés capitales et objectifs du mémoire}

\begin{Lem} \label{lemmedejantzen} Soient $V$ et $W$ deux objets de $L$-$\mathbf{mod}$. Pour chaque $b \in \Bloc(L)$, notons $V_b$ et $W_b$ la somme de tous les sous-modules de $V$ et $W$, respectivement, qui sont dans le bloc $\Cat_b$, respectivement. Alors 
\begin{align*}
V = \bigoplus\textstyle{_{b \, \in \, \Bloc(L)}} \; V_b && W = \bigoplus\textstyle{_{b \, \in \, \Bloc(L)}} \; W_b
\end{align*}
De plus, on a que 
$$ \Hom_L(V,W) = \bigoplus\textstyle{_{b \, \in \, \Bloc(L)}} \; \Hom_L(V_b,W_b) $$
\end{Lem}

Ce résultat est le lemme 5.5 de l'article \cite{NSext}. Ses auteurs, E.Neher et A.Savage, notent que la preuve de ce lemme est la même que celle du lemme II.7.1 de \cite{Jantzen}, qui correspond au même résultat, mais dans le cadre des représentations des groupes algébriques. Ce lemme est encore vrai dans le cadre de ce travail.

\begin{Rem} Dans la catégorie $L$-$\mathbf{mod}$, on peut montrer sans trop de difficulté que chaque objet s'écrit comme une somme directe de sous-modules indécomposables où les classes d'isomorphismes et la multiplicité des sous-modules indécomposables qui apparaissent sont bien définies. Ce résultat porte souvent le nom de « théorème de Krull-Schmidt ». On le retrouve dans plusieurs références, par exemple \cite{Etingof} (voir le théorème 3.8.1).

Le problème 9.53 (ii) de \cite{Etingof} fournit une autre perspective sur le lemme \ref{lemmedejantzen} lorsque combiné au théorème de Krull-Schmidt. Le contexte de cet exercice est cependant un peu trop précis pour pouvoir l'appliquer directement dans le cadre de ce travail.
\end{Rem}

Le lemme \ref{lemmedejantzen} explique en quelque sorte la notation $\Cat = \bigoplus_\alpha \Cat_\alpha$ qui apparait dans la proposition 5.2 de \cite{NSext}. Selon la notation adoptée jusqu'ici dans ce document, on écrira plutôt
\begin{equation*}
L\text{-}\mathbf{mod} \; = \; \bigoplus\textstyle{_{b \, \in \, \Bloc(L)}} \; \Cat_b
\end{equation*}

Cette décomposition et le lemme \ref{lemmedejantzen} justifient certainement l'utilité de la décomposition d'une catégorie en blocs pour arriver à comprendre la composition des objets par leurs facteurs simples. 

Question de préciser ce qui a été dit dans l'introduction, les objectifs de ce travail sont les suivants :

\begin{itemize}
\item[\textbullet] Étiquetter les blocs de la catégorie $\LL$-$\mathbf{mod}$ des modules de dimension finie pour une algèbre de courants tordue $\LL$.
\item[\textbullet] Décrire les ensembles $\Ext_\LL^1(V,W)$ lorsque $V$ et $W$ sont des $\LL$-modules simples de dimension finie.
\item[\textbullet] Par le fait même, décrire les blocs de la catégorie des modules de dimension finie pour les formes tordues d'algèbres de Lie simples (voir \cite{LPtfa}). 
\end{itemize}

\begin{Rem} Les blocs de la catégorie des modules de dimension finie pour les formes tordues n'ont auparavant pas été décrits.
\end{Rem}

\chapter{Algèbres de courants tordues}     

Le coeur du sujet se trouve dans ce chapitre. Tout d'abord, on y rappelle la notion d'une algèbre de courants tordue et quelques commentaires. Puis, il sera question de développer un peu de théorie spécifique à ce travail. Ensuite viendra une description précise du contexte et des objets étudiés. Nous terminerons avec une présentation de la classification des blocs d'extensions de la catégorie des représentations de dimension finie pour une algèbre de courants tordue donnée.

\section{Théorie spécifique au présent travail}
\subsection{Notion d'algèbre de courants tordue}

Rappelons qu'une algèbre de courants tordue est une algèbre de Lie obtenue comme une algèbre de points fixes par l'action d'un groupe fini par $k$-automorphismes sur un produit tensoriel de la forme $\g \otimes S$, où $\g$ est une algèbre de Lie simple de dimension finie et $S$ est une $k$-algèbre de type fini, associative, commutative, unitaire et réduite.

Une algèbre de courants tordue est avant tout une sous-algèbre de Lie de $\g \otimes S$. Son crochet de Lie est donc donné par
\begin{equation} \label{crochetL}
[x \otimes s, y \otimes r] = [x,y] \otimes sr
\end{equation}
où $x,y \in \g$ et $s,r \in S$.

Concrètement, si $\LL$ est une algèbre de courants tordue, alors on peut l'écrire de la manière suivante : 
$$ \LL = (\g \otimes S)^\Gamma = \big\{ \ell \in \g \otimes S \; | \; {}^\gamma\ell = \ell \text{ pour tout } \gamma \in \Gamma\big\} $$
où $\Gamma$ est un groupe fini qui agit (à gauche) sur $\g \otimes S$ par $k$-automorphismes (d'algèbre de Lie).

\begin{Rem} \label{Linfinie} Le plus souvent, et dans le cas de bien des exemples intéressants, les algèbres de courants tordues sont des $k$-algèbres de Lie de dimension infinie. D'un certain point de vue, l'étude de ces objets permet donc de comprendre davantage, ne serait-ce qu'un tout petit peu, une facette de la théorie des représentations d'objets de dimension infinie.
\end{Rem}

Le premier résultat important pour l'analyse d'un tel objet indique comment il est possible de décomposer l'action en utilisant la cohomologie galoisienne. Cette décomposition doit ultimement permettre d'utiliser des informations sur $S$ pour mieux comprendre $\LL$. Voici rapidement quelques définitions nécessaires à ce qui suivra.

\begin{Def} Soit $G$ un groupe. On dit que $A$ est un $G$-groupe (à gauche) si $A$ est un groupe sur lequel $G$ agit par automorphismes, c'est-à-dire où
$$ {}^g(a_1a_2) = {}^g a_1{}^g a_2 $$
pour chaque choix de $a_1,a_2 \in A$ et de $g \in G$.
\end{Def}

On peut ensuite définir la notion de 1-cocycle de $G$ dans un $G$-groupe $A$. Ces objets apparaissent souvent lorsqu'il est question d'étudier la relation entre $G$ et certains types d'objets.

\begin{Def} Un 1-cocycle de $G$ dans un $G$-groupe $A$ est une application $ \alpha : G \rightarrow A $ qui vérifie 
\begin{equation} \label{cohomoGalois}
\alpha(g_1g_2) = \alpha(g_1)\;{}^{g_1}\big(\alpha(g_2)\big)
\end{equation} 
pour tout choix de $g_1,g_2 \in G$. L'ensemble des 1-cocycles de $G$ dans $A$ est noté $Z^1(G \; ; A)$. 

Étant donné un 1-cocycle $\alpha \in Z^1(G \; ; A)$, il est usuel d'employer la notation suivante :
\begin{align*}
\alpha : \; &G \longrightarrow A \\
&g \, \longmapsto \, \alpha_g
\end{align*}
... c'est-à-dire que $\alpha(g) \in A$ est noté $\alpha_g$ en pratique.
\end{Def}

\begin{Rem} Un 1-cocycle $\alpha \in Z^1(G \; ; A)$ n'est évidemment pas un homomorphisme de groupes de $G$ dans $A$, mais le fait que (\ref{cohomoGalois}) soit à un petit symbole près d'en être la propriété requise justifie que les 1-cocycles soient parfois appelés des \textbf{homomorphismes croisés}.
\end{Rem}

Une algèbre de courants tordue correspond aux points fixes d'une algèbre $\g \otimes S$ sous l'action d'un groupe fini $\Gamma$. Dans ce contexte très particulier, on retrouve des 1-cocycles qui permettent de décomposer l'action du groupe $\Gamma$ en des parties plus faciles. La proposition qui suit est en fait la proposition 2.2 de l'article \cite{LauTCA}.

\begin{Prop} \label{1cocycle} Soit $\LL = (\g \otimes S)^\Gamma$ l'algèbre de courants tordue associée à une action de $\Gamma$ sur $\g \otimes S$. Alors, il existe une action $\Gamma \curvearrowright S$ par $k$-automorphismes (d'algèbre de Lie) et un 1-cocycle $u \in Z^1\big(\Gamma,\Aut_{S-\text{Lie}} (\g \otimes S)\big)$ tels que
\begin{align*}
{}^\gamma(x \otimes s) &= \big(u_\gamma \circ (1 \otimes \gamma)\big)(x \otimes s) \\
&= u_\gamma(x \otimes {}^\gamma s)
\end{align*}
pour tout $\gamma \in \Gamma$, $x \in \g$ et $s \in S$.
\end{Prop}

Dans ce cas-ci, le fait que $\g$ soit simple permet d'induire une action $\Gamma \curvearrowright \Max S$ à partir de l'action sur $\g \otimes S$. Puisque $S$ est réduite, son radical de Jacobson est trivial. Il devient alors possible de transférer l'action de $\Gamma$ sur $\Max S$ à une action de $\Gamma$ sur $S$. Une fois l'action $\Gamma \curvearrowright S$ établie, on peut s'intéresser à la collection des applications
$$u_\gamma : \; x \otimes s \; \longmapsto \; ^\gamma(x \otimes {}^{\gamma^{-1}} s) $$
qui mesurent une sorte de « différence entre les actions de $\gamma$ sur $\g \otimes S$ et sur $S$ ». Ces applications définissent en fait un 1-cocycle de $\Gamma$ à valeur dans le $\Gamma$-groupe $\Aut_{S-\text{Lie}} (\g \otimes S)$. 

\begin{Rem} Il n'y a aucune ambiguïté sur l'action induite de $\Gamma$ sur $S$ ; voir à ce propos la preuve du théorème 2.2 de \cite{LauTCA}. Par conséquent, le cocycle $u \in Z^1\big(\Gamma,\Aut_{S-\text{Lie}} (\g \otimes S)\big)$ ci-haut ne fait l'objet d'aucun choix.
\end{Rem}

La structure de $\Gamma$-groupe sous-entendue sur $\Aut_{S-\text{Lie}} (\g \otimes S)$ est donnée par l'action \linebreak $({}^\gamma \phi) = (1 \otimes \gamma) \circ \phi \circ (1 \otimes \gamma^{-1})$, où $\gamma \in \Gamma$ et $\phi \in \Aut_{S-\text{Lie}} (\g \otimes S)$. Des calculs montrent effectivement que chaque $u_\gamma$ est bien un automorphisme $S$-linéaire de l'algèbre de Lie $\g \otimes S$. Le 1-cocycle de la proposition s'écrit donc
\begin{align*}
u : \; &\Gamma \longrightarrow \Aut_{S-\text{Lie}} (\g \otimes S) \\
&\gamma  \longmapsto \quad u_\gamma
\end{align*}

Une première chose bien utile que l'on peut faire avec ce résultat est de montrer que $\LL$ n'est pas qu'une $k$-algèbre de Lie, mais qu'elle est aussi une $S^\Gamma$-algèbre de Lie.

\begin{Prop} \label{actrmod} Les algèbres de courants tordues sont des $S^\Gamma$-algèbres de Lie. 
\begin{pre} Le crochet de Lie de $\LL$ défini en (\ref{crochetL}) est manifestement $S$-linéaire. Il est donc suffisant de s'assurer que $\LL$ est un $S^\Gamma$-module, c'est-à-dire que $S^\Gamma.\,\LL \subseteq \LL$ . 

Fixons $\gamma \in \Gamma$ et prenons $r \in S^\Gamma$ et $\ell = \sum_i x_i \otimes s_i \in \LL$. Concrètement, il s'agit de vérifier que $^\gamma (r.\ell) =r.\ell$. Calculons 
\begin{align*}
^\gamma (r.\ell) &= {}^\gamma \left(r.\Big(\sum_i x_i \otimes s_i\Big)\right) \\
&= {}^\gamma \left(\sum_i x_i \otimes rs_i\right) \\
&= u_\gamma\left(\sum_i x_i \otimes {}^\gamma(rs_i)\right) \\
&= u_\gamma\left(\sum_i x_i \otimes {}^\gamma r ^\gamma s_i)\right) \\
\end{align*}
puisque $r \in S^\Gamma$. On peut donc poursuivre en écrivant
\begin{align*}
^\gamma (r.\ell) &= u_\gamma\left(\sum_i x_i \otimes r ^\gamma s_i\right) \\
&= u_\gamma\left(r.\Big(\sum_i x_i \otimes {}^\gamma s_i\Big)\right) \\
&= r.\left(u_\gamma\Big(\sum_i x_i \otimes {}^\gamma s_i\Big)\right) \\
&= r.^\gamma \ell \\
&= r.\ell
\end{align*}
où au passage à la troisième ligne, la $S$-linéarité de $u_\gamma$ a été utilisée. C'est donc que $r.\ell \in \LL$ et ainsi, $\LL$ est une $S^\Gamma$-algèbre de Lie.
\qed
\end{pre}
\end{Prop}

Ce sera bien utile plus tard de savoir que l'extension d'anneaux $S / S^\Gamma$ est entière.

\begin{Lem} \label{extensionentiere} L'extension d'anneaux $S/S^\Gamma$ est entière.
\begin{pre} Soit $s \in S$, alors le polynôme
$$ \prod_{\gamma \in \Gamma} (x- {}^\gamma s) $$
est un polynôme unitaire à coefficients dans $S^\Gamma$ dont une des racines est $s$. \qed
\end{pre}
\end{Lem}

L'objet de ce travail est de classifier les blocs d'extensions de la catégorie $\LL$-$\mathbf{mod}$. Certainement, il sera question de modules simples. C'est sans surprise que certains idéaux de $\LL$ joueront un rôle important.

\begin{Def} \label{maxl} L'ensemble $\Max \LL$ est l'ensemble des idéaux maximaux $\MM$ d'une algèbre de courants tordue $\LL$  qui soient tels que $\LL' \nsubseteq \MM$. 
\end{Def}

\begin{Rem} Cette définition assure que le quotient de $\LL$ par ce qu'on appelle un idéal maximal soit une algèbre de Lie qui n'est pas abélienne. Le quotient est alors une algèbre de Lie simple dans le sens usuel.

Rappelons que les idéaux d'une algèbre de Lie sont automatiquement bilatères.
\end{Rem}

En fait, il suit assez directement de cette définition que ces idéaux maximaux sont aussi des idéaux de $\LL$ en tant que $S^\Gamma$-algèbre de Lie. En effet, s'il existait $\MM \in \Max \LL$ et $r \in S^\Gamma$ tels que $r\MM \nsubseteq \MM$, alors $r\MM + \MM = \LL$ et puis
\begin{align*}
\LL' = [\LL,\LL] &= [r\MM + \MM, \LL] \\
&= [r\MM, \LL] + [\MM, \LL] \\
&= [\MM, r\LL] + [\MM, \LL] \\
&\subseteq [\MM, \LL] \\
&\subseteq \MM
\end{align*}
Manifestement, ceci contredit l'hypothèse que $\MM$ soit dans $\Max \LL$, alors une telle situation ne se produit jamais. Ainsi, tous les éléments de $\Max \LL$ sont des $S^\Gamma$-idéaux de $\LL$.

\subsection{Modules et algèbres de Lie}

Dans cette sous-section, il est question de résultats incontournables concernant les modules d'une algèbre de Lie $L$. Par ailleurs, une attention plus particulière est portée aux $L$-modules de dimension finie. L'hypothèse voulant que $k$ soit algébriquement clos est généralement superflue dans cette sous-section à moins d'une indication contraire dans les hypothèses des résultats concernés, le cas échéant. 

Dans un premier temps, il peut être amusant de classifier quelque chose puisque le but de ce travail de maîtrise est de présenter une classification.

\begin{Prop} \label{modulesdimension1} Soit $L$ une algèbre de Lie. Alors, il y a une bijection naturelle
$$ \left\{\begin{array}{c} \text{Classes d'isomorphismes de}\\
L\text{-modules de dimension 1} \end{array}\right\} \quad \stackrel{1:1}{\longleftrightarrow} \quad \Hom_k(L/L',k) = (L/L')^* $$
\begin{pre} Soit $V$ un $L$-module de dimension 1. En tant qu'espace vectoriels, $V \cong k$. Pour classifier les classes d'isomorphismes de $L$-modules de dimension 1, il suffit de classifier les actions linéaires de $L$ sur l'espace vectoriel $k$.

Une action linéaire de $L$ sur $k$ est un homomorphisme d'algèbres de Lie \linebreak $\rho : L \rightarrow \End_k(k) \cong k$. Puisque $\End_k(k)$ est une algèbre de Lie abélienne, toute fonction $f \in \Hom_k(L,k)$ qui satisfait $f|_{L'} = 0$ est un homomorphisme d'algèbres de Lie de $L$ dans $\End_k(k)$. L'ensemble des classes d'isomorphismes de $L$-modules de dimension 1 correspond donc à l'ensemble des fonctions $f \in \Hom_k(L,k)$ telles que $L' \subseteq \ker f$. En d'autres mots, toute classe d'isomorphismes correspond à un élément de $\Hom_k(L/L',k) = (L/L')^*$.

D'un autre côté, chaque élément de $(L/L')^*$ correspond à une action linéaire de $L$ sur $k$. Il ne reste donc qu'à s'assurer que deux éléments distincts de $(L/L')^*$ donnent lieu à deux actions distinctes de $L$ sur $k$. Soit donc $\lambda, \mu \in (L/L')^*$ avec $\lambda \neq \mu$. Alors, il existe $\ell + L' \in L/L'$ tel que $\lambda(\ov{\ell}) \neq \mu(\ov{\ell})$ et puis les actions de $\ell \in L$ sur $1 \in k$ sont respectivement
\begin{align*}
\ell.1 = \lambda(\ov{\ell}) \cdot 1 = \lambda(\ov{\ell}) && \ell.1 = \mu(\ov{\ell}) \cdot 1 = \mu(\ov{\ell})
\end{align*}
Elles sont donc différentes et ainsi, la proposition est vraie. \qed
\end{pre}
\end{Prop}

Dans le reste de ce document, les $L$-modules de dimension 1 seront notés $k_\lambda$, où $\lambda \in (L/L')^*$. Si $\ell \in L$, l'écriture $\lambda(\ell)$ signifiera réellement $\lambda(\ell + L')$ et ainsi, les éléments de $(L/L')^*$ seront systématiquement identifiés aux éléments de $L^*$ qui valent $0$ sur $L'$. Cette notation formellement inexacte a cependant toujours un sens bien défini là où elle sera utilisée.

Dans un second temps, voici plusieurs notions de base concernant des modules pour une algèbre de Lie $L$. Ces définitions, constructions et résultats seront bien importants pour la suite.

\begin{Def} \label{invariants} Soit $L$ une algèbre de Lie et soit $V$ un $L$-module. Alors, le \textbf{sous-module des $L$-invariants} de $V$ est l'ensemble 
$$ V^L = \{v \in V \; | \; \ell.v = 0 \text{ pour tout } \ell \in L\} $$
Cet ensemble est exactement formé des éléments de $V$ pour lesquels l'action de $L$ est triviale. C'est en fait le plus grand sous-module trivial de $V$.
\end{Def}

L'ensemble $V^L$ est bien un sous-module de $V$. En effet, $0 \in V^L$. De plus, si $a \in k$ et $v_1,v_2 \in V^L$, alors pour chaque $\ell \in L$ on peut écrire
\begin{align*}
\ell.(av_1+v_2) &= a(\ell.v_1)+\ell.v_2 \\
&= 0 + 0 \\
&= 0
\end{align*}
Ainsi, $V^L$ est un $k$-espace vectoriel et par définition, $L.V^L = \{0\} \subseteq V^L$. L'ensemble $V^L$ est donc certainement aussi un $L$-module. Il sera question de cet important sous-module dans la prochaine sous-section.

\begin{Rem} Si $V$ est un $L$-module simple non-trivial. Alors, $V^L = \{0\}$.
\end{Rem}

\begin{Def} Soient $V$ et $W$ des $L$-modules. Alors, leur \textbf{somme directe} $V \oplus W$ a une structure de $L$-module, où l'action de $\ell \in L$ donnée par
\begin{align*}
\ell \cdot (v,w) = (\ell.v,\ell.w)
\end{align*}
\end{Def}

\begin{Def} Soient $V$ et $W$ des $L$-modules. Alors leur \textbf{produit tensoriel} $V \otimes W$ a une structure de $L$-module où l'action de $\ell \in L$ est donnée par
\begin{align*}
\ell \cdot (v \otimes w) = (\ell.v) \otimes w +v \otimes (\ell.w)
\end{align*}
\end{Def}

\begin{Rem} Soit $V$ un $L$-module. Alors, l'isomorphisme naturel d'espaces vectoriels $k \otimes V \cong V$ devient un isomorhisme de $L$-modules si on munit $k$ de la structure de module trivial. En effet,
\begin{align*}
\varphi: \; &V \rightarrow V \otimes k_0\\
&v \longmapsto v \otimes 1
\end{align*}
vérifie $\; \ell.v \; \stackrel{\varphi}{\longmapsto} \; (\ell.v) \otimes 1 = (\ell.v) \otimes 1 + v \otimes 0 = (\ell.v) \otimes 1 + v \otimes (\ell.1) = \ell.(v \otimes 1) \;$ pour tout choix de $\ell \in L$ et $v \in V$.

En particulier, ce résultat atteste que la catégorie des $L$-modules $L$-$\mathbf{mod}$ est une catégorie monoïdale, c'est-à-dire qu'à isomorphisme près, l'ensemble des objets de $L$-$\mathbf{mod}$ forme un monoïde avec opération $\otimes$ et unité $k_0$. 
\end{Rem}

\begin{Lem} Soient $\lambda, \mu \in (L/L')^*$. Alors, il y a un isomorphisme naturel de $L$-modules
\begin{equation} \label{produitdimension1}
k_\lambda \otimes k_\mu \cong k_{\lambda + \mu}
\end{equation}
\begin{pre} En tant que $k$-espaces vectoriels, $k_\lambda \otimes k_\mu$ et $k_{\lambda + \mu}$ sont naturellement isomorphes via l'identification
\begin{align*}
f : \; k_\lambda &\otimes k_\mu \rightarrow k_{\lambda + \mu}\\
a &\otimes b  \longmapsto ab
\end{align*}
Cette bijection $k$-linéaire est en fait un isomorphisme de $L$-modules. En effet, soit $\ell \in L$ et soit $a \otimes b \in k_\lambda \otimes k_\mu$, alors
\begin{align*}
f\big(\ell \cdot (a \otimes b)\big) &= f\big(\lambda(\ell)a \otimes b + a \otimes \mu(\ell)b\big)\\
&= \lambda(\ell)f(a \otimes b) + \mu(\ell) f(a \otimes b)\\
&= \big(\lambda+\mu\big)(\ell)f(a \otimes b)\\
&= \ell. f(a \otimes b)
\end{align*}
\qed
\end{pre}
\end{Lem}

\begin{Rem} Ce résultat est d'autant plus intéressant au niveau théorique puisqu'il révèle qu'en regard de la bijection de la proposition \ref{modulesdimension1}, l'addition dans l'espace vectoriel $(L/L')^*$ correspond exactement au produit tensoriel du côté des classes d'équivalences de \linebreak $L$-modules de dimension 1.

En particulier, des décompositions de l'espace vectoriel $(L/L')^*$ en sommes directes correspondent à des écritures uniques des $L$-modules de dimension 1 comme des produits tensoriels de $L$-modules de dimension 1 où chacun des « facteurs » correspondent à une des facteurs directs dans la décomposition de $(L/L')^*$ donnée.
\end{Rem}

\begin{Def} Soient $V$ et $W$ des $L$-modules. Alors, l'espace vectoriel $\Hom_k(V,W)$ a une structure naturelle de $L$-module où l'action de $\ell \in L$ est donnée par
\begin{equation} \label{modulehom}
\big(\ell \cdot f\big)(v) = \ell.f(v)-f(\ell.v)
\end{equation}
où $v \in V$.
\end{Def} 

\begin{Rem} \label{invariantshom} Cette structure de module et la définition (\ref{invariants}) du sous-module des invariants entrainent directement que 
$$ \big(\Hom_k(V,W)\big)^L \cong \Hom_L(V,W) $$
Ce module est donc trivial en tant que $L$-module.
\end{Rem}

\begin{Def} Soit $V$ un $L$-module et soit $k_0$ le $L$-module trivial de dimension 1. Alors, le $L$-module $\Hom_k(V,k_0)$ est appelé \textbf{le dual de} $V$. 

Le dual d'un module $V$ est noté $V^*$.
\end{Def}

\begin{Rem} Soit $V$ un $L$-module. L'action de $\ell \in L$ sur $f \in V^*$ est  
$$ \big(\ell \cdot f\big)(v) = \ell.f(v) - f(\ell.v) = - f(\ell.v) $$
\end{Rem}

\begin{Lem} Soit $V$ un $L$-module. Alors, il y a une bijection entre les sous-modules de $V$ et les sous-modules de $V^*$.
\begin{pre} Il est possible de définir les applications ensemblistes suivantes.
\begin{align*}
\varphi : \quad &\{\text{Sous-modules de } V\} \quad \longrightarrow \quad \{\text{Sous-modules de } V^*\}\\
&\qquad\qquad\quad N \;\qquad \longmapsto \quad\;\;\; \{f \in V^* \; | \; f(n) = 0 \text{ pour tout } n \in N\}
\end{align*}
\begin{align*}
\psi : \quad &\{\text{Sous-modules de } V^*\} \quad \longrightarrow \quad \{\text{Sous-modules de } V\}\\
&\qquad\qquad G \;\qquad \longmapsto \quad\;\;\; \{v \in V \; | \; g(v) = 0 \text{ pour tout } g \in G\}
\end{align*}
En effet, on voit directement que si $N$ est un sous-module de $V$ et si $G$ est un sous-module de $V^*$, alors $\varphi(N)$ et $\psi(G)$ sont des espaces vectoriels. Puis, si $\ell \in L$ et $f \in \varphi(N)$, alors pour chaque $n \in N$, $\ell.n \in N$. Il suit assez directement que $L.\varphi(N) \subseteq \varphi(N)$ puisque
$$ \big(\ell \cdot f\big)(n) = - f(\ell.n) = 0 $$
Aussi, si $v \in \psi(G)$, alors pour chaque $g \in G$, $\ell \cdot g \in G$. Il suit assez directement, ici aussi, que $L.\psi(G) \subseteq \psi(G)$ puisque
$$ g (\ell.v) = - \big(\ell \cdot g\big)(v) = 0 $$
Ainsi, les applications $\varphi$ et $\psi$ sont bien définies. Ces applications sont en fait inverses l'une de l'autre. On a effectivement
\begin{align*}
\psi\big(\varphi(N)\big) = N && \varphi\big(\psi(G)\big) = G
\end{align*}
Par exemple, $\psi\big(\varphi(N)\big)$ correspond aux éléments de $V$ qui sont envoyés sur $0$ par tous les éléments de $\varphi(N)$. Par contre, ces derniers sont justement l'ensemble des fonctions qui valent $0$ sur $N$, d'où $N \subseteq \psi\big(\varphi(N)\big)$. Ensuite, si $v \in V \bs N$, on peut choisir un espace vectoriel $U \subseteq V$ tel que $V = N \oplus \Span_k\{v\} \oplus U$. Maintenant, il est possible de définir $\lambda \in V^*$ tel que $\lambda|_{N \oplus U} = 0$ et $\lambda(v) = 1$. Alors, $\lambda \in \varphi(N)$ par définition, puis $v \notin \psi\big(\varphi(N)\big)$ et donc $V \bs N \subseteq V \bs \psi\big(\varphi(N)\big) \; \Rightarrow \; \psi\big(\varphi(N)\big) \subseteq N$. Un raisonnement de la même nature justifie l'autre égalité.

Puisque $\varphi$ et $\psi$ sont des fonctions inverses l'une de l'autre, les deux sont des bijections. \qed
\end{pre}
\end{Lem}

\begin{Cor} \label{dualsimple} Si $V$ est un $L$-module simple, alors $V^*= \Hom_k(V,k_0)$ l'est aussi.
\begin{pre} Par la bijection donnée au lemme précédent et par l'hypothèse donnant que $V$ soit simple, il est automatique que $V^*$ n'a que deux sous-modules. Autrement dit, les seuls sous-modules de $V^*$ sont $\{0\}$ et $V^*$. Ceci donne que $V^*$ est un $L$-module simple. \qed 
\end{pre}
\end{Cor}

\begin{Lem} Soient $\lambda \in (L/L')^*$. Alors, il y a un isomorphisme naturel de $L$-modules
\begin{equation} \label{dualdimension1}
(k_\lambda)^* \cong k_{- \lambda}
\end{equation}
\begin{pre} On sait que l'espace vectoriel $(k_\lambda)^*$ est de dimension 1. En tant que $L$-module, cet espace est isomorphe à $k_\mu$ pour un certain $\mu \in (L/L')^*$. Il ne reste qu'à trouver ce $\mu$.

Soit $\ell \in L$ et soit $f \in (k_\lambda)^*$. Alors, l'action de $\ell$ sur $f : k \rightarrow k$ est décrite par
\begin{align*}
\big(\ell \cdot f\big) : \; a \mapsto - f(\ell.a) &= - f\big(\lambda(\ell)a\big)\\
&= - \lambda(\ell)f(a)
\end{align*}
On a donc que chaque $\ell \in L$ agit sur $(k_\lambda)^*$ via l'application linéaire $- \lambda \in (L/L')^*$. Ainsi, il faut avoir $(k_\lambda)^* \cong k_{- \lambda}$ en tant que $L$-modules. \qed
\end{pre}
\end{Lem}

Supposons maintenant que $V$ et $W$ sont des $L$-modules de dimension finie. Alors en tant que $L$-modules, il y a des isomorphismes naturels
\begin{align}
\label{doubledual} &\bullet \; V^{**} \cong V\\
\label{dualdistributif} &\bullet \; (V \otimes W)^* \cong V^* \otimes W^*\\
\label{homtenseur} &\bullet \; \Hom_k(V,W) \cong V^* \otimes W
\end{align}

Voici maintenant des résultats incontournables en théorie des représentations des algèbres de Lie. Rappelons que l'algèbre universelle enveloppante $\UU(L)$ d'une algèbre de Lie $L$ est la $k$-algèbre associative et unitaire définie comme étant le quotient de l'algèbre tensorielle de $L$ par les relations
$$[x,y] = x \otimes y - y \otimes x $$
pour tout choix de $x,y \in L$. Un élément arbitraire de $\UU(L)$ s'écrit comme une somme de termes composés d'un produit d'éléments de $L$.  

\begin{Thm} \label{equivalenceU} Il y a une équivalence naturelle entre la catégorie des $L$-modules et la catégorie des $\UU(L)$-modules à gauche. 

Si $V$ est un $\UU(L)$-module, alors chaque $\ell \in L$ correspond déjà à un élément de $\End_k V$ et cette croorespondance fait de $V$ un $L$-module. Si maintenant $V$ est un $L$-module, alors le $\UU(L)$-module correspondant est $V$ avec action d'un monôme $\ell_1\cdots \ell_r \;\; (\text{où } \ell_j \in L)$ sur $v \in V$ donnée par
$$ (\ell_1\cdots \ell_r).v = \ell_1.(\cdots(\ell_r.v)\cdots) $$
et où il est supposé que 1 agit comme $\Id_V$.
\end{Thm}
Des détails supplémentaires sont donnés dans \cite{Weibel} autour du théorème 7.3.3.

\begin{Thm} \em \textbf{(Théorème de Poincaré-Birkhoff-Witt)} \em Soit $L$ une algèbre de Lie. Supposons que comme $k$-espace vectoriel, $L$ soit engendrée par une base $\{l_\alpha\}_{\alpha \in I}$ où $I$ est un ensemble d'indices ordonnés. 
Alors, supposant que $\; j > i \; \Rightarrow \; \alpha_j > \alpha_i$, l'ensemble de monômes
$$ \{l_{\alpha_1}^{e_{\alpha_1}} \cdots \, l_{\alpha_r}^{e_{\alpha_r}} \; | \; r \in \NN \text{ et } e_{\alpha_j} \in \NN \text{ pour tout } \alpha_j \in I \} $$
est une base de l'algèbre universelle enveloppante $\UU(L)$. Le « produit vide » correspond ici à $1 \in k \subseteq \UU(L)$.
\end{Thm}
On le retrouve entre autre dans \cite{Weibel}, voir le théorème 7.3.7.


\begin{Cor} Si $L_1$ et $L_2$ sont deux algèbres de Lie, alors il y a un isomorphisme naturel d'algèbres associatives
$$ \UU(L_1 \oplus L_2) \cong \UU(L_1) \otimes \UU(L_2) $$
\end{Cor}

\begin{Cor} \label{sommelie} Soient $L_1$ et $L_2$ deux algèbres de Lie. Supposons que tous les modules qui apparaissent ci-dessous sont de dimension finie. Alors,
\begin{itemize}
\item[\em(i)] Soit $V$ un $L_1$-module simple et soit $W$ un $L_2$-module simple. Alors, $V \otimes W$ est un \linebreak $(L_1 \oplus L_2)$-module simple.
\item[\em(ii)] Tout $(L_1 \oplus L_2)$-module simple est de la forme $V \otimes W$, où $V$ est un $L_1$-module simple et où $W$ est un $L_2$-module simple.
\end{itemize}
\begin{pre} Le résultat suit directement du théorème 3.10.2 de \cite{Etingof} où l'on prend \linebreak $A = \UU(L_1)$ et $B = \UU(L_2)$ en utilisant l'isomorphisme du corollaire précédent ainsi que les équivalences de catégories du théorème \ref{equivalenceU} pour traduire le résultat tel quel. \qed
\end{pre}
\end{Cor}

Soit $L$ une algèbre de Lie et soit $V$ un $L$-module. Si $v \in V$, on peut considérer l'ensemble
$$ \UU(L).v = \{x.v \; | \; x \in \UU(L)\} $$
Cet ensemble est en fait un sous-module de $V$ et il est appelé \textbf{le sous-module} engendré par $v$. Il est évident que $\UU(L).v$ est un sous-espace vectoriel de $V$ étant donné que l'action de $L$ (ou de $\UU(L)$) est linéaire. Aussi, si $x.v \in \UU(L).v$, on a que $x \in \UU(L)$ et donc pour tout $\ell \in L$, $\ell.(x.v) = (\ell x).v \in \UU(L).v$ puisque $\ell x \in \UU(L)$.

\begin{Rem} En particulier, si de plus on suppose que $V$ est un $L$-module qui soit simple, alors pour tout $v \in V \bs \{0\}$, on a $$ \UU(L).v = V $$
C'est que $v \neq 0 \; \Rightarrow \; \{0\} \neq \UU(L).v \;$ qui est un sous-module du module simple $V$.
\end{Rem}

\begin{Lem} \em \textbf{(Lemme de Schur)} \label{schur} \em Soient $V$ et $W$ des $L$-modules simples. Alors, tout homomorphisme non-trivial entre $V$ et $W$ est un isomorphisme. En particulier,
$$ \Hom_L(V,W) \neq \{0\} \qquad \Rightarrow \qquad V \cong W $$
\begin{pre} Soit $f \in \Hom_L(V,W)\bs\{0\}$. Alors, puisque $\ker f \neq V$ est un sous-module du module simple $V$, $\ker f = 0$ et f est injective. Aussi, puisque $f \neq 0$, le sous-module $\im W$ du module simple $W$ est forcément $W$ au complet et donc $f$ est un isomorphisme de $L$-modules. \qed
\end{pre}
\end{Lem}

\begin{Cor} \label{dimhom} Supposons que $k$ soit un corps algébriquement clos. Soit $L$ une $k$-algèbre de Lie et soient $V$ et $W$ des $L$-modules simples de dimension finie. Alors,
\begin{equation*}
\dim_k \Hom_L(V,W) = \left\{\begin{array}{c} 1 \text{ si } V \cong W\\ 0 \text{ si } V \ncong W\end{array}\right.
\end{equation*}
\begin{pre} Si $V \ncong W$, le lemme de Schur donne aussitôt que $\Hom_L(V,W) = 0$. Si maintenant $V \cong W$, on peut choisir $f,g \in \Hom_L(V,W)\bs\{0\}$. Toujours d'après le lemme de Schur, $f$ et $g$ sont des isomorphismes. Il suit que $g^{-1} \circ f \in \End_L(V)$. Puisque $V$ est de dimension finie et que $k$ est algébriquement clos, il existe au moins une valeur propre $\lambda \in k$ de $g^{-1} \circ f$. Ensuite, on considère 
$$ g^{-1} \circ f - \lambda \Id_V \in \End_L(V) $$
Le déterminant de ce nouvel endomorphisme étant nul, on sait qu'il n'est pas inversible alors le lemme de Schur assure qu'il est, en fait, égal à 0. Si $v \in V$, on peut alors écrire
\begin{align*}
\big(g^{-1} \circ f\big)(v) &= \lambda v\\
f(v) &=\lambda g(v)
\end{align*}
Ainsi, $f = \lambda g$ et donc $\dim_k \Hom_L(V,W) = 1$. \qed
\end{pre}
\end{Cor}

\begin{Rem} Ce corollaire permet aussi de savoir que si donc $V$ et $W$ sont des $L$-module simples de dimension finie, alors on sait que $\dim_k \Hom_L(V,W) = \dim_k (V^* \otimes W)^L \leq 1$ et que l'égalité tient si et seulement si $V \cong W$.
\end{Rem}

\begin{Prop} \label{quotientreductif} Supposons que $k$ soit un corps algébriquement clos. Soit $L$ une algèbre de Lie. Supposons que $(\rho,V)$ soit une représentation de $L$ qui soit irréductible et de dimension finie. Alors, $L / \ker \rho$ est réductive et la dimension de son centre est au plus 1.
\begin{pre} Le fait que $L / \ker \rho$ soit réductive est une conséquence de la proposition 5 de la section 6.4 de \cite{Bourbakichap1}. Il ne reste qu'à justifier le fait que la dimension de son centre est d'au plus 1.  

Le reste de la preuve est faite par contradiction. Puisque $L/\ker \rho$ est réductive, écrivons
$$ L / \ker \rho = \Z \oplus \SSS $$ 
où $\Z = Z(L / \ker \rho)$ avec $\dim_k \Z \geq 2$ et où $\SSS$ est une algèbre de Lie semi-simple. 

Soient $z_1,z_2 \in \Z \bs \{0\}$ deux éléments linéairement indépendants. Chaque élément du centre $\Z$ agit sur le module simple $V$ comme la multiplication par un scalaire. En effet, en sachant que $\Id_V \in \End_L(V)$, le corollaire \ref{dimhom} démontre cette affirmation. Ainsi, peu importe $v \in V$, on peut écrire
\begin{align*}
z_1.v = \lambda v && z_2.v = \mu v
\end{align*}
où $\lambda, \mu \in k \bs \{0\}$.
Toujours peu importe $v \in V$, il est donc également vrai que
$$ \left(\frac{\mu}{\lambda}z_1\right).v = \mu v = z_2.v $$
Puisque $V$ est une représentation fidèle (de $\Z \oplus \SSS$), l'égalité précédente donne $\frac{\mu}{\lambda}z_1 = z_2$, mais cela est impossible puisque $z_1$ et $z_2$ ont été supposés linéairement indépendants. Il faut donc admettre que $\dim_k \Z \leq 1$. \qed
\end{pre}
\end{Prop}

\begin{Thm} \em \textbf{(Théorème de Weyl)} \label{weyl} \em Soit $L$ une algèbre de Lie semi-simple de dimension finie sur $k$, un corps algébriquement clos. Alors, toute représentation de dimension finie de $L$ est la somme directe de sous-représentations irréductibles.
\end{Thm}

Pour une preuve, voir par exemple le théorème 2 de la section 6.2 de \cite{Bourbakichap1}. Il vaut la peine d'élaborer encore davantage sur cette propriété qu'a une représentation de s'écrire comme la somme directe de sous-représentations irréductibles. Pour commencer, nommons-la.

\begin{Def} Un $L$-module qui, en tant que module, est isomorphe à une somme directe de module simple est appelé \textbf{un module semi-simple} ou encore \textbf{une représentation complètement réductible}.  
\end{Def}

Soit $V$ un $L$-module semi-simple et soit $\A$, un ensemble complet de représentants des classes d'isomorphismes de $L$-modules simple (sans répétition). Alors, il y a par définition un isomorphisme de $L$-modules
$$ V \cong \bigoplus_{A \in \A} \; A^{\,\oplus\, n_A(V)} $$
où $n_A(V) \in \NN$ pour chaque $A \in \A$. Ces nombres représentent simplement les multiplicités des $A \in \A$ comme facteurs de $V$. 

\begin{Rem} Si on suppose que $V$ est de dimension finie, on peut évidemment prendre la somme directe sur l'ensemble $\A_{df}$ des éléments de $\A$ qui sont des modules de dimension finie. De plus, les $n_A(V)$ valent toujours 0, sauf pour un nombre fini de $A \in \A_{df}$.
\end{Rem}

\begin{Prop} Soit $L$ une $k$-algèbre de Lie et soit $V$ un $L$-module de dimension finie. Alors,
$$ V \text{ est semi-simple} \quad \Longleftrightarrow \quad V \cong \bigoplus_{A \in \A_{df}} \big(\Hom_L(A,V) \otimes A\big) $$
où le symbole $\;\cong\;$ fait référence à un isomorphisme particulier décrit plus bas.
\begin{pre} Soit $V$ un $L$-module de dimension finie. Posons 
$$ \A_{df}^V = \{A \in \A_{df} \; | \; n_A(V) \neq 0\} $$
La remarque précédente explique en fait que $\# \A_{df}^V < \infty$. Par définition, $V$ est donc un module semi-simple si et seulement si $V \cong \bigoplus_{A \in \A_{df}^V} \; A^{\,\oplus\, n_A(V)}$. D'abord, comme la somme directe est finie, on sait que pour chaque $B \in \A_{df}^V$ on peut écrire 
\begin{align*}
\Hom_L(B,V) &\cong \Hom_L(B,\textstyle{\bigoplus_{A \in \A_{df}^V}} \; A^{\,\oplus\, n_A(V)})\\
&\cong \textstyle{\bigoplus_{A \in \A_{df}^V}} \big(\Hom_L(B,A)\big)^{\,\oplus\, n_A(V)}
\end{align*}
Ensuite, $B$ de même que tous les $A$ qui apparaissent dans l'expression précédente sont des modules simples de dimension finie. Ainsi, le corollaire \ref{dimhom} donne que 
$$ \dim_k \Hom_L(A,B) = \delta_{A,B} $$
Puis, si $A = B$, $\Hom_L(A,A)$ est un module de dimension 1 et la remarque \ref{invariantshom} nous assure que c'est un $L$-module trivial. Il est donc naturellement isomorphe à $k_0$ et on peut écrire
$$ \Hom_L(B,V) \cong k_0^{\,\oplus\, n_B(V)} $$

On peut maintenant dire que $V$ est un module semi-simple
\begin{align*}
&\Longleftrightarrow \qquad V \cong \bigoplus_{A \in \A_{df}^V} \; A^{\,\oplus\, n_A(V)}\\
&\Longleftrightarrow \qquad V \cong \bigoplus_{A \in \A_{df}^V} \; (\underbrace{A \oplus \cdots \oplus A}_{n_A(V) \text{ fois}})\\
&\Longleftrightarrow \qquad V \cong \bigoplus_{A \in \A_{df}^V} \; \big(\underbrace{(k_0 \otimes A) \oplus \cdots \oplus (k_0 \otimes A)}_{n_A(V) \text{ fois}}\big)\\
&\Longleftrightarrow \qquad V \cong \bigoplus_{A \in \A_{df}^V} \; \big(\underbrace{(k_0 \oplus \cdots \oplus k_0)}_{n_A(V) \text{ fois}} \otimes A)\big)\\
&\Longleftrightarrow \qquad V \cong \bigoplus_{A \in \A_{df}^V} \; (k_0^{\,\oplus\, n_A(V)} \otimes A)\\
&\Longleftrightarrow \qquad V \cong \bigoplus_{A \in \A_{df}^V} \big(\Hom_L(A,V) \otimes A\big)\\
&\Longleftrightarrow \qquad V \cong \bigoplus_{A \in \A_{df}} \big(\Hom_L(A,V) \otimes A\big)
\end{align*}
\qed
\end{pre}
\end{Prop}

\subsection{Cohomologie des algèbres de Lie}

Tel que décrit au chapitre 1, les notions de cohomologies sont définies à partir de complexes de chaînes donnés. Dans le cas des algèbres de Lie, ce qui est appelé  « la cohomologie » est quelque chose de bien déterminé en ce sens qu'elle est issue d'un complexe spécifique. Quoi qu'il en soit, cette notion de cohomologie permet de tirer une foule d'informations sur une algèbre de Lie donnée et sur ses modules. Dans le cas présent, il sera important de relier les notions d'extensions de modules et de cohomologie afin de tirer profit des différents résultats en cohomologie des algèbres de Lie. 

Soit $L$ une $k$-algèbre de Lie et soient $(\rho_V,V)$ et $(\rho_W,W)$ deux $L$-modules. Tel qu'expliqué dans le chapitre 1, $\Ext_L^1(V,W)$ est l'ensemble des classes d'équivalences d'extensions de $V$ par $W$. Soit $(\rho_E,E)$ une extension de $V$ par $W$ ; $E$ est un $L$-module et qu'il y a des homomorphismes de $L$-modules de façon à former une suite exacte courte
\begin{equation} \label{extensionsection}
0 \longrightarrow W \longrightarrow E \longrightarrow V \longrightarrow 0
\end{equation}
Puisque cette suite est aussi une suite exacte d'espaces vectoriels, il y a un isomorphisme d'espaces vectoriels $E = W \oplus s(V) \cong W \oplus V$ pour chaque section $s : V \rightarrow E$ de l'extension (\ref{extensionsection}). Lorsqu'une section $s$ est fixée, peut alors voir $\rho_E$ comme étant une application à valeurs dans $\End_k(W \oplus V)$. Étant donné $\ell \in L$, on peut donc écrire
\begin{equation} \label{matriceextension}
\rho_E(\ell) = \matt{\rho_W(\ell)}{d_E^s(\ell)}{0}{\rho_V(\ell)} : \Vv{w}{v} \longmapsto \Vv{\rho_W(\ell)w + d_E^s(\ell)v}{\rho_V(\ell)v} \in W \oplus V
\end{equation}
où $d_E^s(\ell) \in \Hom_k(V,W)$. 

\begin{Rem} La présence du zéro dans la matrice ci-dessus est dû au fait que $W \subseteq E$ est stable sous l'action de $L$ et au contexte qui fait que l'application $E = W \oplus V \twoheadrightarrow V$  qui définit l'extension est la projection naturelle.
\end{Rem}

\begin{Rem} L'application $d_E^s : L \rightarrow \Hom_k(V,W)$ est $k$-linéaire.
\end{Rem}

Maintenant, puisque $\rho_E : L \rightarrow \End_k(W \oplus V)$ est un homomorphisme d'algèbres de Lie, pour chaque choix de $\ell_1,\ell_2 \in L$, on doit avoir
$$ \rho_E\big([\ell_1,\ell_2]\big) = \rho_E(\ell_1)\rho_E(\ell_2) - \rho_E(\ell_2)\rho_E(\ell_1) $$
La matrice (\ref{matriceextension}) avec $\ell = [\ell_1,\ell_2]$ doit donc être égale à 
$$ \matt{\rho_W(\ell_1)}{d_E^s(\ell_1)}{0}{\rho_V(\ell_1)}\matt{\rho_W(\ell_2)}{d_E^s(\ell_2)}{0}{\rho_V(\ell_2)} - \matt{\rho_W(\ell_2)}{d_E^s(\ell_2)}{0}{\rho_V(\ell_2)}\matt{\rho_W(\ell_1)}{d_E^s(\ell_1)}{0}{\rho_V(\ell_1)} $$
Ce qui se traduit par les relations
\begin{align}
\rho_W\big([\ell_1,\ell_2]\big) &= \rho_W(\ell_1)\rho_W(\ell_2) - \rho_W(\ell_2)\rho_W(\ell_1) \\
\label{derextension1} d_E^s\big([\ell_1,\ell_2]\big) &= \rho_W(\ell_1)d_E^s(\ell_2) + d_E^s(\ell_1)\rho_V(\ell_2) - \rho_W(\ell_2)d_E^s(\ell_1) - d_E^s(\ell_2)\rho_V(\ell_1)\\
\rho_V\big([\ell_1,\ell_2]\big) &= \rho_V(\ell_1)\rho_V(\ell_2) - \rho_V(\ell_2)\rho_V(\ell_1)
\end{align}
La première et la troisième des relations apparaissant ci-haut n'apportent rien de nouveau puisque l'on supposait déjà que $V$ et $W$ étaient des $L$-modules. Par contre, la relation (\ref{derextension1}) apporte de l'information sur $E$. En gardant en tête que $d_E^s(\ell) \in \Hom_k(V,W)$ pour chaque $\ell \in L$ donné et en se souvenant de la structure naturelle de module sur $\Hom_k(V,W)$ donnée par (\ref{modulehom}), la relation (\ref{derextension1}) s'écrit directement
\begin{equation} \label{derextension2}
d_E^s\big([\ell_1,\ell_2]\big) = \big(\ell_1 \cdot d_E^s(\ell_2)\big) - \big(\ell_2 \cdot d_E^s(\ell_1)\big) 
\end{equation}

Avant d'aller plus loin, il est souhaitable d'introduire sans plus tarder la notion de dérivation pour éventuellement traduire simplement ce qui décrit $E$.

\begin{Def} Soit $U$ un $L$-module. Une \textbf{dérivation de $L$ dans} $U$ est une application linéaire $ d : L \rightarrow U$ qui satisfait la relation
$$ d\big([\ell_1,\ell_2]\big) = \ell_1.d(\ell_2) - \ell_2.d(\ell_1) $$
pour tout choix de $\ell_1,\ell_2 \in L$. L'ensemble des dérivations de $L$ dans $U$ est noté $\Der(L,U)$.
\end{Def}

\begin{Def} Une dérivation $d \in \Der(L,U)$ est dite \textbf{intérieure} s'il existe $u \in U$ tel que pour tout $\ell \in L$, on ait
$$ d(\ell) = \ell.u $$
Cette dérivation sera dorénavant notée $d^u$. 

L'ensemble des dérivations intérieures de $L$ dans $U$ est noté $\IDer(L,U)$.
\end{Def}

\begin{Rem} Les dérivations intérieures sont bien des dérivations puisque si $u \in U$ et $\ell_1, \ell_2 \in L$, alors $d^u\big([\ell_1,\ell_2]\big) = [\ell_1,\ell_2].u = \ell_1.(\ell_2.u) - \ell_2.(\ell_1.u) = \ell_1.\big(d^u(\ell_2)\big) - \ell_2.\big(d^u(\ell_1)\big)$.
\end{Rem}

\begin{Rem} Les ensembles $\Der(L,U)$ et $\IDer(L,U)$ sont des $k$-espaces vectoriels.
\end{Rem}

\begin{Rem} \label{seider} Soit $V$ un $L$-module. Alors, il existe toujours une suite exacte courte
$$ 0 \rightarrow V^L \rightarrow V \rightarrow \IDer(L,V) \rightarrow 0 $$
En effet, il y a une application surjective $app : V \twoheadrightarrow \IDer(L,V)$ qui à un élément $v \in V$ associe la dérivation intérieure $d^v \in \IDer(L,V)$. Le noyau de cette application est l'ensemble des éléments $v \in V$ qui sont tels que $d^v = 0$. Cette condition est équivalente à avoir $\ell.v = 0$ pour chaque $\ell \in L$. Ainsi, $\ker (app) = V^L$, le sous-module des $L$-invariants de $V$. 
\end{Rem}

D'après (\ref{matriceextension}), il est immédiat qu'une extension $E$ de $V$ par $W$ dans la catégorie $L$-$\mathbf{mod}$ est complètement déterminée par la connaissance de l'application linéaire $d_E^s : L \rightarrow \Hom_k(V,W)$. À la lumière des précédentes définitions, le fait que cette application $d_E^s$ satisfasse la relation (\ref{derextension2}) se traduit simplement par l'appartenance
$$ d_E^s \in \Der\big(L,\Hom_k(V,W)\big) $$
Ainsi, une extension (de $V$ par $W$) correspond à une certaine dérivation de l'algèbre de Lie $L$. Cette dérivation dépend bien entendu du choix de la section $s : V \rightarrow E$ de l'extension (\ref{extensionsection}), mais il se trouve que le choix d'une différente section, $\sigma$, de l'extension (\ref{extensionsection}) donne lieu à une relation toute particulière entre $d_E^s$ et $d_E^\sigma$ ; ils appartiennent au même co-ensemble de l'espace quotient
$$ \Der\big(L,\Hom_k(V,W)\big) / \IDer\big(L,\Hom_k(V,W)\big) $$

Pour expliquer ce lien entre $d_E^s$ et $d_E^\sigma$, commençons par écrire 
\begin{align} \label{wtfsectionext1}
E_s = W \oplus V \cong W \oplus s(V) && E_\sigma = W \oplus V \cong W \oplus \sigma(V)
\end{align}
Alors l'isomorphisme naturel suivant fait le lien entre les $L$-modules $E_s$ et $E_\sigma$ :
\begin{align*}
\Phi = \matt{\Id_W}{s-\sigma}{0}{\Id_V} : E_s = W \oplus V \rightarrow W \oplus V = E_\sigma
\end{align*}
Notons que $\Phi$ est un isomorphisme de $L$-modules d'après les définitions de la ligne (\ref{wtfsectionext1}). En particulier, la structure de module sur $E_s$ doit vérifier 
\begin{equation} \label{wtfsectionext2}
\rho_{E_s}(\ell) = \Phi^{-1} \circ \rho_{E_\sigma} \circ \Phi
\end{equation}
Grâce à cette équation (\ref{wtfsectionext2}), on peut déduire que 
\begin{align*}
\Vv{\rho_W(\ell)w + d_E^s(\ell)v}{\rho_V(\ell)v} &= \Phi^{-1} \left(\Vv{\rho_W(\ell)w + \rho_W(\ell)\big((s-\sigma)v\big) + d_E^\sigma(\ell)v}{\rho_V(\ell)v}\right) \\
&= \Vv{\rho_W(\ell)w + \rho_W(\ell)\big((s-\sigma)v\big) + d_E^\sigma(\ell)v - \Big(s-\sigma\Big)\big(\rho_V(\ell)v\big)}{\rho_V(\ell)v}
\end{align*}
pour tout $\ell \in L$, $w \in W$ et $v \in V$. Ainsi, on trouve la relation
\begin{equation} \label{wtfsectionext3}
d_E^s(\ell) = d_E^\sigma(\ell) + \big(\ell \cdot (s - \sigma)\big) \in \Hom_k(V,W)
\end{equation}
Finalement, puisque $s - \sigma : V \rightarrow E$ est réellement une application $k$-linéraire à valeurs dans le sous-espace $W \subseteq E$, on voit que le second terme de (\ref{wtfsectionext3}) implique
\begin{equation} \label{wtfsectionext4}
\ov{d_E^s} \, = \, \ov{d_E^\sigma} \in \; \Der\big(L,\Hom_k(V,W)\big) / \IDer\big(L,\Hom_k(V,W)\big)
\end{equation}

\begin{Rem} Étant donné la relation (\ref{wtfsectionext4}), le choix d'une section $s$, $\sigma$ (ou toute autre) de l'extension (\ref{extensionsection}) ne sera absolument pas important pour la suite des choses. Ceci peut-être bien compris à la lecture de la prochaine proposition.

En conséquence, l'écriture du symbole $d_E$ sous-entendra un choix préalablement arrêté d'une section de l'extension (\ref{extensionsection}) dans tout le reste du document.
\end{Rem}

Le prochain résultat permet de caractériser une classe d'équivalences de $\Ext_L^1(V,W)$ arbitraire. Ceci permet de cerner l'information réellement intéressante fournie par la donnée d'une extension $E$.

\begin{Prop} \label{dercohomo} Soit $L$ une algèbre de Lie. Deux extensions $E_1$ et $E_2$ de $V$ par $W$ sont équivalentes
$$ \qquad \quad \Longleftrightarrow \qquad \ov{d_{E_1}} = \; \ov{d_{E_2}} \;\; \in \Der\big(L,\Hom_k(V,W)\big) / \IDer\big(L,\Hom_k(V,W)\big) $$
\begin{pre} Soient $E_1$ et $E_2$ deux extensions telles que $\; \ov{d_{E_1}} = \; \ov{d_{E_2}}$. Soit $f \in \Hom_k(V,W)$ tel que $d_{E_1} = d_{E_2} + d^f$. Comme espaces vectoriels, on a encore
\begin{align*}
E_1 &= W \oplus V \\
E_2 &= W \oplus V
\end{align*}
Dans ce cas, l'application 
\begin{align*}
\varphi : \quad &\quad\; E_1 \quad \longrightarrow \quad E_2\\
&\Vv{w}{v} \mapsto \Vv{w+f(v)}{v}
\end{align*}
est une équivalence d'extensions entre $E_1$ et $E_2$. En effet, c'est une application linéaire inversible dont l'inverse est 
\begin{align*}
\varphi^{-1} : \quad &\quad\; E_2 \quad \longrightarrow \quad E_1\\
&\Vv{w}{v} \mapsto \Vv{w-f(v)}{v}
\end{align*}
Ensuite, si $\ell \in L$ et $w + v \in E_1 = W \oplus V$, alors
\begin{align*}
\ell.\varphi \Vv{w}{v} = \ell.\Vv{w + f(v)}{v} &= \rho_{E_2}(\ell)\Vv{w + f(v)}{v}\\
&= \matt{\rho_W(\ell)}{d_{E_2}(\ell)}{0}{\rho_V(\ell)}\Vv{w + f(v)}{v}\\
&= \Vv{\rho_W(\ell)w + \rho_W(\ell)f(v) + d_{E_2}(\ell)v}{\rho_V(\ell)v}
\end{align*}
est bien égal à 
\begin{align*}
\varphi\Big(\ell.\Vv{w}{v}\Big) &= \varphi\Vv{\rho_W(\ell)w + d_{E_1}(\ell)v}{\rho_V(\ell)v}\\ 
&= \Vv{\rho_W(\ell)w + d_{E_1}(\ell)v + f\big(\rho_V(\ell)v\big)}{\rho_V(\ell)v}
\end{align*}
puisque d'abord, les deuxièmes composantes de ces vecteurs sont égales et que les premières composantes sont aussi égales par l'hypothèse à l'effet que $d_{E_1} = d_{E_2} + d^f$. En effet
\begin{align*}
d^f(\ell)v + d_{E_2}(\ell)v \quad &= \quad d_{E_1}(\ell)v\\
\big(\ell \cdot f\big)(v) + d_{E_2}(\ell)v \quad &= \quad d_{E_1}(\ell)v\\
\rho_W(\ell)f(v) - f\big(\rho_V(\ell)v\big) + d_{E_2}(\ell)v \quad &= \quad d_{E_1}(\ell)v\\
\rho_W(\ell)f(v) + d_{E_2}(\ell)v \quad &= \quad d_{E_1}(\ell)v + f\big(\rho_V(\ell)v\big)\\
\rho_W(\ell)w + \rho_W(\ell)f(v) + d_{E_2}(\ell)v \quad &=\quad \rho_W(\ell)w + d_{E_1}(\ell)v + f\big(\rho_V(\ell)v\big)\\
\end{align*}
Enfin, le diagramme qui suit commute puisque l'application $\varphi$ envoie un vecteur \linebreak$(0,v) \in W \oplus V = E_1$ sur le vecteur $\big(f(v),v\big)\in W \oplus V = E_2$.
\begin{center}\begin{tikzpicture}
y\matrix (m) [matrix of math nodes, row sep=1em,column sep=2.5em,text height=1.5ex,text depth=0.25ex, nodes in empty cells]
{ & & W \oplus V \, & & \\ 0 & W & & V & 0\\ & & W \oplus V \, & & \\};
\path (m-1-3)++(-0.25ex,3.33ex) node {$\overbrace{\quad \, \quad \, \quad}^{\; E_1}$}
         (m-3-3)++(-0.25ex,-3.33ex) node {$\underbrace{\quad \quad \quad}_{\; E_2}$};
\draw (m-2-1) edge[->,thick] (m-2-2)
          (m-2-2) edge[->,thick] (m-1-3)
          (m-2-2) edge[->,thick] (m-3-3)
          (m-1-3) edge[->,thick] node[auto] {\footnotesize{$\varphi$}} (m-3-3)
          (m-1-3) edge[->,thick] (m-2-4)
          (m-3-3) edge[->,thick] (m-2-4)
          (m-2-4) edge[->,thick] (m-2-5);
\end{tikzpicture} \end{center}

Maintenant, il reste à montrer que si $E_1$ et $E_2$ sont des extensions équivalentes, alors \linebreak $\ov{d_{E_1}} = \; \ov{d_{E_2}}$ dans le quotient des dérivations par les dérivations intérieures. Par hypothèse, on a un diagramme commutatif
\begin{center}\begin{tikzpicture}
y\matrix (m) [matrix of math nodes, row sep=0.5em,column sep=2.5em,text height=1.5ex,text depth=0.25ex, nodes in empty cells]
{ & & E_1 & & \\ 0 & W & & V & 0\\ & & E_2 & & \\};
\draw (m-2-1) edge[->,thick] (m-2-2)
          (m-2-2) edge[->,thick] (m-1-3)
          (m-2-2) edge[->,thick] (m-3-3)
          (m-1-3) edge[->,thick] node[auto] {\footnotesize{$\varphi$}} (m-3-3)
          (m-1-3) edge[->,thick] (m-2-4)
          (m-3-3) edge[->,thick] (m-2-4)
          (m-2-4) edge[->,thick] (m-2-5);
\end{tikzpicture} \end{center}
où $\varphi : E_1 = W \oplus V \rightarrow W \oplus V = E_2$ est un isomorphisme de $L$-modules. Puisque le diagramme commute, si $w \in W$, alors il faut avoir $\varphi(w + 0) = w + 0$ et si $p : W \oplus V \rightarrow V$, il faut avoir $p = p \circ \varphi$. Donc, pour un $v \in V$,   on a $\varphi(w + v) = (w + 0) + (\tilde{w} + v)$ où $\tilde{w} \in W$ dépend linéairement de $v$. Tout ceci se traduit par le fait que la matrice qui représente 
$$ \varphi : E_1 = W \oplus V \rightarrow W \oplus V = E_2$$ 
est de la forme
$$ \matt{\Id_W}{f}{0}{\Id_V} $$

pour un certain $f \in \Hom_k(V,W)$. De plus, le fait que $\varphi$ soit un isomorphisme de $L$-modules permet d'écrire que peu importe $\ell \in L$, $w \in W$ et $v \in V$, on a l'égalité
\begin{align*}
\ell.\varphi(w + v) &= \varphi\big(\ell.(w + v)\big) \\
\ell.\Vv{w + f(v)}{v} &= \varphi\Big(\Vv{\rho_W(\ell)w + d_{E_1}(\ell)v}{\rho_V(l)v}\Big)\\
\Vv{\rho_W(\ell)w + \rho_W(\ell)f(v) + d_{E_2}(\ell)v}{\rho_V(\ell)v} &= \Vv{\rho_W(\ell)w + d_{E_1}(\ell)v + f\big(\rho_V(\ell)v\big)}{\rho_V(\ell)v}
\end{align*}
Cela donne directement 
\begin{align*}
\rho_W(\ell)f(v) + d_{E_2}(\ell)v &= d_{E_1}(\ell)v + f\big(\rho_V(\ell)v\big)\\
\rho_W(\ell)f(v) - f\big(\rho_V(\ell)v\big) + d_{E_2}(\ell)v &= d_{E_1}(\ell)v\\
\big(\ell \cdot f\big)(v) + d_{E_2}(\ell)v &= d_{E_1}(\ell)v\\
d^f(\ell)v + d_{E_2}(\ell)v &= d_{E_1}(\ell)v
\end{align*}
On peut donc conclure que $\; \ov{d_{E_1}} = \; \ov{d_{E_2}} \in \Der\big(L,\Hom(V,W)\big) / \IDer\big(L,\Hom(V,W)\big)$. \qed
\end{pre}
\end{Prop}

\begin{Cor} \label{zeroexttriviale} Soit $E$ une extension de $V$ par $W$ telle que $d_E \in \IDer\big(L.\Hom(V,W)\big)$, alors $E$ est une extension de $V$ par $W$ équivalente à l'extension triviale
$$ 0 \rightarrow W \rightarrow W \oplus V \rightarrow V \rightarrow 0 $$ 
En particulier, $E \cong W \oplus V$ comme $L$-module.
\end{Cor}

En fait, il y a une notion de cohomologie pour les algèbres de Lie que l'on peut lier très fortement aux quelques résultats ci-haut. L'avantage de ce faire est de donner accès à tous les théorèmes (et idées) qui ont été développés dans le cadre de la cohomologie des algèbres de Lie. Cette théorie donne beaucoup d'outils pour l'analyse des extensions dans $\LL$-$\mathbf{mod}$, la catégorie des modules de dimension finie d'une algèbre de courants tordue $\LL$. Sans introduire tout ce langage « cohomologique », des résultats de cohomologie des algèbres de Lie plus avancés seront, tout de même, mis à contribution dans la suite de ce travail. 

\begin{Def} \label{1cohomo} Soit $L$ une algèbre de Lie et soit $V$ un $L$-module. Alors, \textbf{la première cohomologie de $L$ à valeurs dans $V$} est l'ensemble
\begin{equation}
\cohomo{L}{V} = \Der(L,V) / \IDer(L,V)
\end{equation}
\end{Def}

\begin{Rem} Soit $V$ un $L$-module. Il existe une suite exacte
$$ 0 \rightarrow V^L \rightarrow V \rightarrow \Der(L,V) \rightarrow \cohomo{L}{V} \rightarrow 0 $$
Les premières applications étant les mêmes que celles décrites à la remarque \ref{seider}. La dernière application est simplement le passage au quotient.
\end{Rem}

Pour tout ce travail, il sera suffisant de ne définir que la première cohomologie. Par ailleurs, la proposition \ref{dercohomo} permet d'interpréter directement les classes d'équivalences d'extensions entre deux $L$-modules donnés en termes de la première cohomologie.

\begin{Prop} \label{cohomoext} Soient $V$ et $W$ deux $L$-modules. Il y a une bijection naturelle entre les classes d'équivalences d'extensions $\Ext_L^1(V,W)$ et la première cohomologie de $L$ à valeurs dans le $L$-module $\Hom_k(V,W)$. En somme,
\begin{equation} \label{extcohomo}
\Ext_L^1(V,W) \quad \stackrel{1:1}{\longleftrightarrow} \quad \gcohomo{L}{\Hom_k(V,W)}
\end{equation}
où $\; \stackrel{1:1}{\longleftrightarrow} \;$ fait référence à la bijection naturelle donnée par la proposition \ref{dercohomo}.
\end{Prop}

\begin{Rem} Soient $V$ et $W$ des $L$-modules. Puisque de la bijection naturelle \ref{extcohomo} et de la définition \ref{1cohomo}, on peut attribuer à $\Ext_L^1(V,W)$ la structure d'un espace vectoriel. Concrètement, si $\;Classe_\spadesuit \; \leftrightarrow \; \ov{d_\spadesuit}\;$ et $\;Classe_\clubsuit \; \leftrightarrow \; \ov{d_\clubsuit}\;$ sont des classes d'équivalences d'extensions de $V$ par $W$, alors comme modules
\begin{align}
Classe_\spadesuit &= \left(\rho_\spadesuit, W \oplus V\right) \qquad \text{où } \rho_\spadesuit(\ell) = \matt{\rho_W(\ell)}{d_\spadesuit(\ell)}{0}{\rho_V(\ell)} \quad \text{pour chaque } \ell \in L \\
Classe_\clubsuit &= \left(\rho_\clubsuit, W \oplus V\right) \qquad \text{où } \rho_\clubsuit(\ell) = \matt{\rho_W(\ell)}{d_\clubsuit(\ell)}{0}{\rho_V(\ell)} \quad \text{pour chaque } \ell \in L
\end{align}
De plus, si $a \in k$, alors $\;Classe_\spadesuit + a \cdot Classe_\clubsuit\;$ est la classe d'équivalences de l'extension donnée par $(\rho_{\spadesuit+a\clubsuit},W \oplus V)$ où 
$$ \rho_{\spadesuit+a\clubsuit}(\ell) = \matt{\rho_W(\ell)}{\big(d_\spadesuit + ad_\clubsuit\big)(\ell)}{0}{\rho_V(\ell)} \quad \text{pour chaque } \ell \in L $$
Ainsi, $\;Classe_\spadesuit + a \cdot Classe_\clubsuit \; \leftrightarrow \; \ov{d_\spadesuit+ad_\clubsuit}\;$. Cette addition est donc associative, chaque classe a un inverse additif et le corollaire \ref{zeroexttriviale} assure l'existence d'un élément neutre additif. De plus, la multiplication scalaire est compatible avec cette addition. Cette structure d'espace vectoriel sur $\Ext^1_L(V,W)$ donne un sens au traitement de tels ensembles comme des espaces vectoriels. Par exemple, il sera souvent question de sommes directes $\oplus$.
\end{Rem}

Lorsque $V$ et $W$ sont des $L$-modules de dimension finie, on peut utiliser l'isomorphisme naturel (\ref{homtenseur}) pour écrire 
\begin{equation}
\Ext_L^1(V,W) \cong \cohomo{L}{V^* \otimes W}
\end{equation}
Puis, on peut écrire
\begin{align*}
\cohomo{L}{V^* \otimes W} &\cong \gcohomo{L}{k_0 \otimes (V^* \otimes W)}\\
&\cong \gcohomo{L}{\Hom_k(k_0,V^* \otimes W)}
\end{align*}
\begin{align*}
\cohomo{L}{V^* \otimes W} &\cong \gcohomo{L}{(V^* \otimes W) \otimes k_0}\\
&\cong \gcohomo{L}{(V \otimes W^*)^* \otimes k_0}\\
&\cong \gcohomo{L}{\Hom_k(V \otimes W^*,k_0)}
\end{align*}

\begin{Cor}
Soit $L$ une algèbre de Lie et soient $V$ et $W$, des $L$-modules de dimension finie. Il y a alors des isomorphismes naturels
\begin{equation} \label{switchext}
\Ext_L^1(V \otimes W^*,k_0) \cong \Ext_L^1(V,W) \cong \Ext_L^1(k_0,V^* \otimes W)
\end{equation}
\end{Cor}

Dans un autre ordre d'idée, si $T$ est un $L$-module trivial et si $\ell_1,\ell_2 \in L$ sont des éléments arbitraires de l'algèbre de Lie, alors 
\begin{align*}
\IDer(L,T) &= \{d^t : \ell \mapsto \ell.t = 0 \; | \; t \in T\}\\ 
&= 0\\
\text{et } \quad \Der(L, T) &= \{d \in \Hom_k(L,T) \; | \; d\big([\ell_1,\ell_2]\big) = \ell_1.d(\ell_2) - \ell_2.d(\ell_1) = 0\}\\
&= \{d \in \Hom_k(L,T) \; | \; d(L') = 0\}
\end{align*}

\begin{Cor} Soit $L$ une algèbre de Lie et soit $T$ un $L$-module trivial. Il y a alors des isomorphismes
\begin{equation} \label{dertrivial}
\cohomo{L}{T} \cong \Der(L,T) \cong \Hom_k(L/L',T)
\end{equation}
En particulier, $\cohomo{L}{k_0} \cong (L/L')^*$.
\end{Cor}

\begin{Rem} Soit $V$ un $L$-module. Alors, la définition du sous-module des $L$-invariants $V^L$ (définition \ref{invariants}) en fait le plus grand sous-module de $V$ pour lequel $L$ agit comme $0$. Ainsi, pour tout $L$-module $V$, on aura toujours 
$$ \cohomo{L}{V^L} \cong \Hom_k(L/L',V^L) $$
\end{Rem} 

\begin{Cor} Soit $\lambda \in (L/L')^*$. Alors 
\begin{equation} \label{extdimension1pareils}
\Ext_L^1(k_\lambda,k_\lambda) \cong (L/L')^*
\end{equation}
\begin{pre} Par les isomorphismes (\ref{produitdimension1}), (\ref{dualdimension1}) et  (\ref{switchext}), on peut écrire
$$ \Ext_L^1(k_\lambda,k_\lambda) \cong \Ext_L^1(k_0,k_0) \cong \cohomo{L}{k_0} $$
Enfin, les isomorphismes naturels (\ref{dertrivial}) du corollaire précédent donnent le résultat. \qed
\end{pre}
\end{Cor}

Une question légitime survient : comment décrire $\Ext_L^1(k_\lambda,k_\mu)$ lorsque $\lambda \neq \mu \in (L/L')^*$ ? 

La proposition suivante donne une réponse à cette question en reprenant plusieurs notions vues jusqu'ici. Le niveau de satisfaction que procure cette réponse est cependant variable. Il dépend de façon proportionelle au niveau de confiance en la capacité du lecteur à calculer les objets considérés.

\begin{Prop} Soit $\lambda \in (L/L')^* \bs \{0\}$. Considérant $\lambda$ comme un élément de $L^*$ qui s'annulle sur $L'$, posons $K = \ker \lambda \subseteq L$. On a alors un isomorphisme naturel
\begin{equation*}
\cohomo{L}{k_\lambda} \cong \Hom_{L/K}\left(K/K',k_\lambda\right)
\end{equation*}
\begin{pre} Puisque $\lambda \neq 0$, il existe $u \in L$ tel que $\lambda(u) = 1$. Le but avoué de cette proposition est de décrire l'espace quotient des dérivations de $L$ dans $k_\lambda$ par les dérivations intérieures. Pour commencer, traitons quelque peu des dérivations intérieures.

L'ensemble $\IDer(L,k_\lambda)$ est formé des dérivations $d^a$ où $a \in k$. Soit $\ell \in L$, alors
\begin{align*}
d^a(\ell) &= \ell.a\\
&= a(\ell.1)\\
&= a\lambda(\ell)
\end{align*}
Donc, $\IDer(L,k_\lambda) = \Span_k \lambda = k \, \lambda$. En particulier, c'est un espace vectoriel de dimension 1. Ensuite, toute dérivation $d \in \Der(L,k_\lambda)$ peut s'écrire
$$ d = \big(d - d(u)\lambda\big) + d(u)\lambda $$
On peut alors remarquer directement que $d(u) \lambda \in \IDer(L,k_\lambda)$ et que $\big(d - d(u)\lambda\big)$ est une dérivation qui vaut zéro lorsqu'on l'évalue en $u \in L$. Par ailleurs, si une dérivation $d$ est telle que $d(u) = 0$, alors $\big(d - d(u)\lambda\big) = d$ et $d(u)\lambda = 0$. L'ensemble $\DD_0$ des dérivations $d$ telles que $d(u) = 0$ est un sous-espace vectoriel de $\Der(L,k_\lambda)$ et on a 
$$ \Der(L,k_\lambda) = \DD_0 \oplus \IDer(L,k_\lambda) $$
Enfin, la définition de la première cohomologie donne $\cohomo{L}{k_\lambda} \cong \DD_0$. Il suffira donc de décrire plus précisément $\DD_0$ et le tour sera joué.

D'abord, comme $\dim_k(k_\lambda) = 1$, $K \subseteq L$ est un idéal de codimension 1. Il s'en suit que chaque $\ell \in L$ peut s'écrire comme
\begin{equation} \label{codim1}
\ell = a_\ell u + \omega
\end{equation}  
où $a_\ell \in k$ et $\omega \in K$. Ceci justifie le fait que $\DD_0 \hookrightarrow \Der(K,k_\lambda)$ puisque si $d \in \DD_0$, alors $d(\ell) = a_\ell \, d(u) + d(\omega) = d(\omega)$. Puis, comme $k_\lambda$ est un $K$-module trivial, les isomorphismes (\ref{dertrivial}) donnent que 
$$ \DD_0  \quad \hookrightarrow \quad \Der(K,k_\lambda) \cong \Hom_k\big(K / K',k\big) = \big(K / K'\big)^* $$
De plus, si $d \in \DD_0$ et $\omega \in \ker \lambda$, alors on peut écrire
\begin{align*}
d(\omega) &= \big(d - d(u)\lambda\big)(\omega)\\
&= d(\omega) - \lambda(\omega)d(u)\\
&= \lambda(u)d(\omega) - \lambda(\omega)d(u)\\
&= u.d(\omega) - \omega.d(u)\\
&= d\big([u,\omega]\big)
\end{align*}
C'est donc dire que pour tout $\omega \in K$, on a $d(\omega - [u,\omega]) = 0$. Ainsi, 
$$ \DD_0 \quad \hookrightarrow \quad \Big\{f \in \big(K / K'\big)^* \; | \; f\big(\,\ov{\omega - [u,\omega]}\,\big) = 0 \text{ pour tout } \omega \in K \Big\} $$

Il faut maintenant montrer l'inclusion inverse. Soit $f \in \big(K / K'\big)^*$ tel que $f\big(\omega - [u,\omega]\big) = 0$ pour tout $\omega \in K$. Ayant alors en tête (\ref{codim1}), on peut voir $f$ comme une dérivation $\tilde{f} \in \DD_0$ par la formule $\tilde{f}(Au + \omega) = f\big(\ov{\omega}\big)$. Cette définition donne immédiatement $\tilde{f}(u) = 0$. 

Il ne reste alors qu'à s'assurer que $\tilde{f} \in \Der(L,k_\lambda)$. Si $\ell_1 = Au + \omega$ et $\ell_2 = Bu + \delta$ sont des éléments arbitraires de $L$, on a $[\ell_1,\ell_2] = A[u,\delta] - B[u,\omega] + [\omega,\delta]$ et donc, 
\begin{align*}
\tilde{f}\big([\ell_1,\ell_2]\big) &= A\tilde{f}\big([u,\delta]\big) - B\tilde{f}\big([u,\omega]\big)\\
&= Af(\ov{\delta}) - Bf(\ov{\omega})\\
&= \lambda(\ell_1)f(\ov{\delta}) - \lambda(\ell_2)f(\ov{\omega})\\
&= \lambda(\ell_1)\tilde{f}(\ell_2) - \lambda(\ell_2)\tilde{f}(\ell_1)\\
&= \ell_1.\tilde{f}(\ell_2) - \ell_2.\tilde{f}(\ell_1)
\end{align*}
Ainsi, on peut conclure que 
$$ \DD_0 = \Big\{f \in \big(K / K'\big)^* \; | \; f\big(\omega - [u,\omega]\big) = 0 \text{ pour tout } \omega \in K \Big\} $$ 
Il ne reste plus qu'à montrer qu'il est égal à $\Hom_{L/K}(K/K',k_\lambda)$. Un élément $f \in \big(K / K'\big)^*$ satisfait la propriété $f\big(\,\ov{\omega - [u,\omega]}\,\big) = 0 \text{ pour tout } \omega \in K$
\begin{align*}
&\Longleftrightarrow \qquad f\big(\,\ov{B\omega - B[u,\omega]}\,\big) = 0 \quad \text{pour tout } \omega \in K \text{ et } B \in k\\
&\Longleftrightarrow \qquad f\big(\,\ov{B\lambda(u)\omega - [Bu,\omega]}\,\big) = 0 \quad \text{pour tout } \omega \in K \text{ et } B \in k\\
&\Longleftrightarrow \qquad f\big(\,\ov{\lambda(Bu + \delta)\omega - [Bu,\omega]-[\delta,\omega]}\,\big) = 0 \quad \text{pour tout } \omega,\delta \in K \text{ et } B \in k\\
&\Longleftrightarrow \qquad f\big(\,\ov{\lambda(Bu + \delta)\omega - [Bu + \delta,\omega]}\,\big) = 0 \quad \text{pour tout } \omega,\delta \in K \text{ et } B \in k\\
&\Longleftrightarrow \qquad f\big(\,\ov{\lambda(\ell)\omega - [\ell,\omega]}\,\big) = 0 \quad \text{pour tout } \omega \in K \text{ et } \ell \in L\\
&\Longleftrightarrow \qquad \lambda(\ell)f\big(\,\ov{\omega}\,\big) - f\big(\,\ov{[\ell,\omega]}\,\big) = 0 \quad \text{pour tout } \omega \in K \text{ et } \ell \in L\\
&\Longleftrightarrow \qquad \ell.f\big(\,\ov{\omega}\,\big) - f\big(\,\ell.\ov{\omega}\,\big) = 0 \quad \text{pour tout } \omega \in K \text{ et } \ell \in L\\
&\Longleftrightarrow \qquad \lambda(\ell + \theta)f\big(\,\ov{\omega}\,\big) - f\big(\,\ov{[\ell+\theta,\omega]}\,\big) = 0 \quad \text{pour tout } \omega, \theta \in K \text{ et } \ell \in L\\
&\Longleftrightarrow \qquad (\ell + \theta).f\big(\,\ov{\omega}\,\big) - f\big(\,(\ell+\theta).\ov{\omega}\,\big) = 0 \quad \text{pour tout } \omega, \theta \in K \text{ et } \ell \in L\\
&\Longleftrightarrow \qquad f \in \Hom_{L/K}(K/K',k_\lambda)\\
\end{align*}
Enfin, on obtient $\cohomo{L}{k_\lambda} \cong \DD_0 = \Hom_{L/K}(K/K',k_\lambda)$. \qed
\end{pre}
\end{Prop}

\begin{Cor} \label{extdim1} Soient $\lambda,\mu \in (L/L')^*$. Posons $K = \ker (\mu - \lambda)$. Alors,
$$ \Ext_L^1(k_\lambda,k_\mu) \cong \Hom_{L/K}(K/K',k_{\mu - \lambda}) $$
\end{Cor}

\begin{Cor} \label{extab} Soit $L$ une algèbre de Lie abélienne et soient $\lambda, \mu \in L^* = (L/L')^*$. Alors
$$ \Ext_L^1(k_\lambda,k_\mu) \cong \left\{\begin{array}{cl} L^* & \text{si } \lambda = \mu \\ \{0\} & \text{si } \lambda \neq \mu\end{array}\right. $$
\begin{pre} Avec le corollaire précédent, on obtient directement la réponse proposée dans le cas où $\lambda = \mu$. Par contre, si $\lambda \neq \mu$, l'ensemble $\Ext_L^1(k_\lambda,k_\mu)$ vaut $\Hom_{L/K}(K/K',k_{\mu - \lambda})$ où $K = \ker (\mu - \lambda)$. 

Étant donné que $L$ est abélienne, $K/K'$ est un $L/K$-module trivial. Si on choisit un $f$ dans $\Hom_{L/K}(K/K',k_{\mu - \lambda})$, il suit directement que $f(K/K') \subseteq (k_{\mu-\lambda})^{L/K}$. Par ailleurs, $k_{\mu - \lambda}$ est un $L/K$-module simple qui est non trivial par hypothèse et donc $(k_{\mu-\lambda})^{L/K} = \{0\}$. Ainsi, $f(K/K') = \{0\}$ ou, autrement dit, $f = 0$ (la fonction nulle). \qed
\end{pre}
\end{Cor}

Après tous ces résultats sur les extensions entre modules de dimension 1, voici un résultat cohomologique relativement général, mais bien utile :

\begin{Prop} \label{extsommelie} Soient $L_1$ et $L_2$ deux algèbres de Lie et soient $V = V_1 \otimes V_2$ et \linebreak $W = W_1 \otimes W_2$ des $(L_1 \oplus L_2)$-modules simples. On peut supposer que $V_i$ et $W_i$ sont des $L_i$-modules simples d'après le corollaire \ref{sommelie}. Sous ces hypothèses, on a 
\begin{equation}
\Ext_{L_1 \oplus L_2}^1(V, W) \cong \left\{
\begin{array}{cl} \{0\} & \text{si } V_1 \ncong W_1 \text{ et } V_2 \ncong W_2\\
\Ext_{L_2}^1(V_2,W_2) & \text{si } V_1 \cong W_1 \text{ et } V_2 \ncong W_2\\
\Ext_{L_1}^1(V_1,W_1) & \text{si } V_1 \ncong W_1 \text{ et } V_2 \cong W_2\\
\Ext_{L_1}^1(V_1,W_1) \oplus \Ext_{L_2}^1(V_2,W_2) & \text{si } V_1 \cong W_1 \text{ et } V_2 \cong W_2
\end{array} \right.
\end{equation}
\end{Prop}

Pour comprendre une preuve de la proposition ci-haut, voir la proposition 2.7 b) de \cite{NSext}. Il s'agit de jouer avec les liens de dualité entre les notions générales d'homologie et de cohomologie pour algèbres de Lie afin de réinterpréter la formule de Künneth, une formule donnant des isomorphismes naturels entres certains ensembles d'homologie, mais pour les ensembles de cohomologie.

En terminant la présente sous-section, voici un résultat qui se révélera très utile dans la démarche de classification des $\Bloc$s pour les algèbres de courants tordues. La preuve de ce résultat est plutôt élémentaire et se résume à manipuler différents types de dérivations.

\begin{Prop} \label{mapremierepropyo} Supposons que $k$ soit un corps algébriquement clos. Soit $L$ une $k$-algèbre de Lie réductive et de dimension finie. 

Supposons que $T$ soit un $L$-module trivial et que $V$ soit un $L$-module simple, non-trivial et de dimension finie. Alors,
$$ \Ext_L^1(T,V) \cong \{0\} $$
\begin{pre} Par définition, $\Ext_L^1(T,V) \cong \gcohomo{L}{\Hom_k(T,V)}$. Montrer le lemme, c'est montrer que toute dérivation de $\Der\big(L,\Hom_k(T,V)\big)$ est intérieure. Soit donc $d \in \Der\big(L,\Hom_k(T,V)\big)$.

Par définition, $d \in \Hom_k\big(L,\Hom_k(T,V)\big)$ et pour toute combinaison de $x,y \in L$ et de $t \in T$, la propriété des dérivations s'écrit
\begin{align*}
\Big(d\big([x,y]\big)\Big)(t) &= \big(x \cdot d(y)\big)(t) - \big(y \cdot d(x)\big)(t)\\
&= x.\Big(\big(d(y)\big)(t)\Big) - \big(d(y)\big)(x.t) + y.\Big(\big(d(x)\big)(t)\Big) - \big(d(x)\big)(y.t)\\
&= x.\Big(\big(d(y)\big)(t)\Big) - y.\Big(\big(d(x)\big)(t)\Big)
\end{align*}
Fixons maintenant un $t \in T$. On obtient alors une dérivation de $L$ dans $V$
\begin{align} \label{deltat}
\delta_t : \; &L \; \longrightarrow \; V\\
&x \longmapsto \big(d(x)\big)(t) \nonumber
\end{align}
En effet, $\delta_t$ est $k$-linéaire et si $x,y \in L$, on a directement 
\begin{align*}
\delta_t\big([x,y]\big) &= \Big(d\big([x,y]\big)\Big)(t)\\
&= x.\Big(\big(d(y)\big)(t)\Big) - y.\Big(\big(d(x)\big)(t)\Big)\\
&= x.\delta_t(y) - y.\delta_t(x)
\end{align*}
Maintenant, les hypothèses de la courante proposition permettent de justifier que \linebreak $\Ext_L^1(k_0,V) = \{0\}$. Ainsi nous pourrons continuer en supposant que 
$$ \delta_t \in \Der(L,V) = \IDer(L,V) $$ 
pour ce même $t \in T$ préalablement fixé (et arbitraire). 

Voici d'abord comment justifier que $\Ext_L^1(k_0,V) = \{0\}$. Puisque $L$ est une algèbre de Lie réductive et de dimension fine sur un corps algébriquement clos, on peut écrire
$$ L = \Z \oplus \SSS $$
où $\Z$ est une algèbre de Lie abélienne et où $\SSS$ est une algèbre de Lie semi-simple sur un corps algébriquement clos. Ensuite, la partie (ii) du corollaire \ref{sommelie} permet d'écrire
$$ V = k_\lambda \otimes V_\SSS $$
où $\lambda \in \Z^*$ et où $V_\SSS$ est un $\SSS$-module simple de dimension finie. Notons que le $L$-module trivial $k_0$ s'écrit $k_0 \otimes k_0$ où le premier facteur peut-être vu comme un $\Z$-module et où le second peut-être vu comme un $\SSS$-module.

Il est ici possible d'utiliser la proposition \ref{extsommelie} avec $L_1 = \Z$, $L_2 = \SSS$. Celle-ci nous donne d'abord que
\begin{equation} \label{omgpourlewin}
\Ext_{L}^1(k_0, V) \cong \left\{
\begin{array}{cl} \{0\} & \text{si } 0 \neq \lambda \text{ et } k_0 \ncong V_\SSS \\
\Ext_\SSS^1(k_0,V_\SSS) & \text{si } 0 = \lambda \text{ et } k_0 \ncong V_\SSS \\
\Ext_\Z^1(k_0,k_\lambda) & \text{si } 0 \neq \lambda \text{ et } k_0 \cong V_\SSS \\
\Ext_\Z^1(k_0,k_\lambda) \oplus \Ext_\SSS^1(k_0,V_\SSS) & \text{si } 0 = \lambda \text{ et } k_0 \cong V_\SSS
\end{array} \right.
\end{equation}
Notons ensuite que par le théorème de Weyl (voir théorème \ref{weyl}), $\Ext_\SSS^1(k_0,V_\SSS) \cong \{0\}$ puisque la catégorie $\SSS$-$\mathbf{mod}$ est semi-simple. Notons également que le corollaire \ref{extab} donne que 
$$ \Ext_\Z^1(k_0,k_\lambda) \cong \left\{\begin{array}{cl} \Z^* & \text{si } 0 = \lambda \\ \{0\} & \text{si } 0 \neq \lambda\end{array}\right. $$
Combinant ces dernières informations dans l'équation (\ref{omgpourlewin}), celle-ci devient
\begin{equation} \label{omgpourlewinv2}
\Ext_{L}^1(k_0, V) \cong \left\{
\begin{array}{cl} \{0\} & \text{si } V \text{ est un } L\text{-module non-trivial} \\
\Z^* & \text{si } V \text{ est un } L\text{-module trivial}
\end{array} \right.
\end{equation}

... mais par hypothèse, le $L$-module simple $V$ est non-trivial ; ce qui donne bien 
$$ \Ext_{L}^1(k_0, V) \cong \{0\} $$ 
tel que voulu. Ainsi, pour chaque $t \in T$ fixé, la dérivation $\delta_t$ donnée à la ligne (\ref{deltat}) est intérieure, c'est-à-dire que 
$$ \delta_t \in \IDer(L,V) $$

D'après la remarque \ref{seider}, la suite d'espaces vectoriels suivante est exacte :
$$ 0 \rightarrow V^L \rightarrow V \rightarrow \IDer(L,V) \rightarrow 0 $$
Les suites exactes courtes d'espaces vectoriels sont toutes scindées. Ainsi, il existe certainement une section dans la catégorie des espaces vectoriels
$$s : \IDer(L,V) \rightarrow V$$
Autrement dit, $s$ est une application $k$-linéaire. Par définition de cette section, la propriété suivante est satisfaite :
$$ \delta_t(x) =x.s(\delta_t) \qquad \text{pour tout } x \in L $$ 
À $\delta_t \in \IDer(L,V)$, on peut donc faire correspondre un $v_t = s(\delta_t) \in V$.  Considérons maintenant l'application ensembliste
\begin{align*}
f_d : \; &T \longrightarrow V\\
&\,t \, \longmapsto \, v_t
\end{align*}
Montrons que c'est une application $k$-linéaire. Soient $t_1,t_2 \in T$ et soit $a \in k$. Alors pour tout $x \in L$, on peut écrire
\begin{align*}
\delta_{at_1+t_2}(x) &= \big(d(x)\big)(at_1+t_2)\\
&=a\big(d(x)\big)(t_1)+\big(d(x)\big)(t_2)\\
&=a\delta_{t_1}(x)+\delta_{t_2}(x)
\end{align*}
L'application $t \mapsto \delta_t$ est donc $k$-linéaire. Puisque la section $s : \IDer(L,V) \rightarrow V$ est aussi $k$-linéaire, la composition de ces deux applications, c'est-à-dire $f_d$, est $k$-linéaire. Autrement dit, $f_d \in \Hom_k(T,V)$.

Enfin, pour tout $t \in T$, on peut écrire
\begin{align*}
\big(d(x)\big)(t) &= \delta_t(x)\\
&= x.v_t\\
&= x.f_d(t)\\
&= x.f_d(t) - f_d(x.t)\\
&= \big(x \cdot f_d\big)(t)
\end{align*}
Ainsi, $d(x) = x \cdot f_d$ pour tout $x \in L$. La dérivation arbitraire $d \in \Der\big(L,\Hom_k(T,V)\big)$ est donc intérieure. En conclusion, on obtient
$$ \Ext_L^1(T,V) \cong \gcohomo{L}{\Hom_k(T,V)} \cong \{0\} $$
\vspace{-2em}
\qed
\end{pre}
\end{Prop}


%
\subsection{Modules d'évaluation pour une algèbre de courants tordue}

Dans cette sous-section, $k$ sera un corps algébriquement clos.

La notion de module d'évaluation est apparue dans les années 1980. Elle a notamment été utilisée pour classifier des représentations d'algèbres de lacets. Depuis, cette même analyse a pu être adaptée à des contextes plus généraux. Récemment, les modules de dimension finie des algèbres d'applications équivariantes et des formes tordues d'algèbres de Lie simples ont été étudiés en détail grâce à cette approche. Dans cette courte sous-section, les modules d'évaluation sont introduits dans le contexte des algèbres de courants tordues. Ce type de module joue un rôle de premier plan dans le présent effort de classification.

Soit $M \in \Max S$. Puisque $S$ est une $k$-algèbre de type fini, $S = k[x_1,...\,,x_t] / I$ où $I$ est un idéal de $k[x_1,...\,,x_t]$. Par la correspondance des idéaux d'un quotient, on sait que $M = \MM/I$ où $\MM \in \Max k[x_1,...\,,x_t]$. Le Nullstellensatz d'Hilbert (voir \cite{GWsym}, théorème A.1.3) donne que $\MM = \langle \, x_1-a_1,...\,,x_t-a_t \, \rangle$ pour certains scalaires $a_i \in k$. Ainsi, on obtient
\begin{equation} \label{résidu} 
S / M \cong \frac{k[x_1,...\,,x_t]/I}{\MM/I} \cong \frac{k[x_1,...\,,x_t]}{\MM} = \frac{k[x_1,...\,,x_t]}{\langle \, x_1-a_1,...\,,x_t-a_t \, \rangle} \cong k
\end{equation}

En fait, il n'existe qu'un tel isomorphisme. En effet, si $f : k \rightarrow S/M$ est un isomorphisme de $k$-algèbres, alors pour chaque $a \in k$,
$$ f(a) = a \cdot f(1) = a \cdot (1 + M) = a+M $$

\begin{Def} Soit $s \in S$ et soit $M \in \Max S$. On note $s(M)$ l'élément de $k$ qui correspond à $s + M \in S/M$ par l'isomorphisme (\ref{résidu}). On appelle $s(M)$ le \textbf{résidu} de $s$ selon $M$. 
\end{Def}

Le résultat qui suit caractérise simplement cette notion.

\begin{Lem} Soit $s \in S$. Alors, $s(M)$ est l'unique élément de $k$ tel que 
$$ s - s(M) \in M $$
\begin{pre} Soit $s \in S$. Par définition de $s(M)$, on a $s + M \, \leftrightarrow \, s(M)$ selon l'isomorphisme $S/M \cong k$. Ainsi, puisque $1+M \, \leftrightarrow \, 1 \in k$, on peut écrire $\big(s - s(M) + M\big) \, \leftrightarrow \, 0 \in k$ et par conséquent,
$$ 0 \, \leftrightarrow \, 0 + M = s-s(M) + M $$
Ce qui montre que $s-s(M) \in M$.

Supposons maintenant qu'il existe $a \in k$ avec $a \neq s(M)$, mais $s-a \in M$. Alors, on a aussi $s-s(M) \in M$ et donc
$$ (s - a) - \big(s - s(M)\big) = s(M) - a \in M \bs \{0\} $$
Mais ceci est impossible puisque $s(M)-a \in k^\times \subseteq S^\times$ ce qui entrainerait $M = S$. \qed
\end{pre}
\end{Lem}

\begin{Cor} L'application $\varepsilon_M : s \mapsto s(M)$ est un homomorphisme surjectif de \linebreak $k$-algèbres entre $S$ et $k$.
\begin{pre} Il est facile de voir que
$$ \varepsilon_M : s \mapsto s + M \mapsto s(M) $$
n'est que la composition de l'isomorphisme (\ref{résidu}) et de la projection naturelle de $S$ dans $S/M$ qui sont tous deux des homomorphismes de $k$-algèbres surjectifs. \qed
\end{pre}
\end{Cor}

\begin{Rem} L'application $\varepsilon_M : S \rightarrow k$ est parfois appelée une évaluation de $S$ en l'idéal maximal $M$. Une telle appellation prend tout son sens en géométrie algébrique. Néanmoins, il peut-être suffisant de considérer un exemple tout simple pour comprendre la terminologie. La \linebreak $\CC$-algèbre, $\CC[z]$ est unitaire, associative, commutative et de type fini, puis $\CC$ est un corps algébriquement clos. Dans ce cas, nous savons que 
$$ \Max \CC[z] = \{M_a = (x-a)\CC[z] \; | \; a \in \CC\} $$
Voir par exemple \cite{GWsym}, théorème A.1.3 avec $X = \CC$.

L'évaluation d'un élément de $\CC[z]$ en un idéal maximal $M_a$ est donc
$$ \varepsilon_{M_a} : p(z) \longmapsto p(a) \in \CC $$
Il s'avère que dans le quotient $\CC[z]/M_a$, nous avons $\ov{z} = \ov{a}$. L'évaluation d'un polynôme \linebreak $p(z) \in \CC[z]$ en l'idéal maximal $M_a$ est donc tout simplement l'évaluation de $p$ en $a \in \CC$. 
\end{Rem}

Si maintenant un groupe $\Gamma$ agit sur $S$ par $k$-automorphismes (de $k$-algèbre), on peut trouver une relation utile entre $s(M)$ et ${}^\gamma s(M)$. D'abord, notons que l'action de $\Gamma$ sur $S$ induit une action $\Gamma \curvearrowright \Max S$. En effet, l'image d'une idéal d'un algèbre par un automorphisme est encore un idéal. Supposons donc que $M \in \Max S$ et montrons que ${}^\gamma M \in \Max S$ pour un $\gamma \in \Gamma$ fixé. Soit $s \in S \bs {}^\gamma M$, alors puisque $M$ est maximal, il existe un élément $r \in S$ tel que $(r + M) \cdot({}^{\gamma^{-1}} s + M) = 1 + M$. C'est donc que $r \cdot {}^{\gamma^{-1}}s - 1 \in M$ alors
$$ {}^\gamma (r \cdot {}^{\gamma^{-1}}s - 1) = {}^\gamma r \cdot s - 1 \in  {}^\gamma M $$
C'est donc que $({}^\gamma r + {}^\gamma M) \cdot ( s + {}^\gamma M) = (1 + {}^\gamma M)$. Puisque $s \in S$ est arbitraire, $S/{}^\gamma M$ est un corps et ${}^\gamma M$ est maximal.

Soient $\gamma \in \Gamma$, $s \in S$ et $M \in \Max S$. Notons $m = s-s(M) \in M$. Alors,
\begin{align*} 
{}^\gamma m &=  {}^\gamma \big(s-s(M)\big) \\
&= {}^\gamma s - {}^\gamma \big(s(M)\big) \\
&= {}^\gamma s - s(M)
\end{align*}
On en conclut la relation $s(M) = {}^\gamma s( {}^\gamma M)$. Cette relation peut aussi s'écrire
\begin{equation}
{}^\gamma s(M) = s({}^{\gamma^{-1}} M) 
\end{equation}
En particulier, si on choisit $\gamma \in \Gamma^M$, l'équation précédente s'écrit simplement ${}^\gamma s(M) = s(M)$.

Maintenant, comment peut-on combiner évaluations et algèbres de Lie ? Soit $L$ une $k$-algèbre de Lie et soit l'algèbre de Lie $L \otimes S$. En tant qu'espaces vectoriels, il est relativement facile d'observer que
$$ L \cong L \otimes k \cong L \otimes S/M $$
Ainsi, étant donné un élément de $L \otimes S$, on peut penser à évaluer les éléments de $S$ en un idéal maximal $M \in \Max S$ pour obtenir un élément de $L$. L'application qui réalise ce dessin est
\begin{align*}
\ev_M = \Id \otimes \, \varepsilon_M : \; L &\otimes S \rightarrow L \otimes S/M \cong L \\
\ell &\otimes s \quad \longmapsto \quad s(M)\ell
\end{align*}

\begin{Rem} Soit $M \in \Max S$. Alors, l'application $\ev_M : L \otimes S \rightarrow L$ est un homomorphisme d'algèbres de Lie. Rappelons que le crochet de Lie de $L \otimes S$ est donné par 
$$ [x \otimes s,y \otimes r] = [x,y] \otimes sr $$
où $x, y \in L$ et $s,r \in S$.
\end{Rem}

Évaluer une algèbre de Lie $L \otimes S$ en un idéal maximal de $S$ est donc un véritable moyen de produire des $(L \otimes S)$-modules à partir de modules pour l'image de l'application d'évaluation via le \linebreak « pullback » de l'action des $\big(\im(\ev_M)\big)$-modules par $\ev_M$.

\begin{Def} Soit $L$ une $k$-algèbre de Lie, alors un $(L \otimes S)$-module $(\rho,E)$ est dit \textbf{d'évaluation} s'il existe $M_1,...\,,M_r \in \Max S$ tel que l'action $\rho : L \otimes S \rightarrow \End_k(E)$ se factorise à travers une action de l'algèbre de Lie
$$ \tilde{\rho} : \quad \im\left(\bigoplus_{i=1}^r \ev_{M_i}\right) \quad \longrightarrow \quad \End_k(E) $$

Posons $\ev_\mathbf{M} = \bigoplus_{i=1}^r \ev_{M_i}$. Rappelons que « l'action $\rho$ se factorise à travers l'action $\tilde{\rho}$ » si le diagramme suivant commute
\begin{center}\begin{tikzpicture}
y\matrix (m) [matrix of math nodes, row sep=2.5em,column sep=2.5em,text height=1.5ex,text depth=0.25ex, nodes in empty cells]
{ L \otimes S & \End_k(E) \\ \im(\ev_\mathbf{M}) &  \\};
\draw (m-1-1) edge[->,thick] node[auto] {\footnotesize{$\rho$}} (m-1-2)
	  (m-1-1) edge[->>,thick] node[left] {\footnotesize{$\ev_\mathbf{M}$}}(m-2-1)
	  (m-2-1)++(5.25ex,1.85ex) edge[->,thick,dashed] node[below] {\footnotesize{$\tilde{\rho}$}} (m-1-2);
\end{tikzpicture} \end{center}
\end{Def}

Cette notion de module d'évaluation est bien définie dans le cas d'une algèbre de courants tordue $\LL = (\g \otimes S)^\Gamma$. Supposons que $M \in \Max S$, alors on peut montrer que
$$ \ev_M(\LL) \cong \g^M $$
où $\g^M \subseteq \g$ est une sous-algèbre de Lie définie comme suit :
\begin{equation} \label{definitiongm}
\g^M = \{x \in \g \; | \; \gamma(M)\big(x\big) = x \quad \text{pour tout } \gamma \in \Gamma^M\} 
\end{equation}
Dans cette définition, $\gamma(M) \in \Aut_k(\g)$ est défini par $ \gamma(M)\big(x\big) = \ev_M\big({}^\gamma(x \otimes 1)\big)$ tel qu'expliqué dans \cite{LauTCA} (voir l'équation 2.6).

\begin{Rem} Le fait que $\ev_M(\LL) \cong \g^M$ est montré dans la démonstration de la proposition 2.9 de \cite{LauTCA}.
\end{Rem}

\begin{Rem} Dans le cas des algèbres de courants tordues, un module d'évaluation pour $\LL$ est avant tout un module pour une algèbre de Lie de dimension finie sur $k$. C'est là un constat significatif étant donné que la théorie des représentations en dimension finie est bien mieux comprise que celle pour les algèbres de dimension infinie. La remarque \ref{Linfinie} prend donc un sens différent dans ce contexte.
\end{Rem}

Avant de terminer cette sous-section, il est fondamental de savoir que les $\g^M$ sont des algèbres de Lie réductives. 

\begin{Prop} \label{gmreductive} Soit $M \in \Max S$. Alors, $\g^M \subseteq \g$ est une algèbre de Lie réductive.
\end{Prop}

Ce fait est expliqué dans la preuve du théorème 3.2 de \cite{LauTCA}. Le résultat est dû à l'isomorphisme $\g^M \cong (\g \otimes S/M)^{\Gamma^M}$ mentionné plus tôt dans l'article et à la proposition 14 de la section 1.4 de \cite{Bourbakichap78}.

\section{Objets simples de $\LL$-$\mathbf{mod}$}
Les deux prochaines sous-sections résument les résultats principaux de l'article \cite{LauTCA} dont le but est précisément d'étudier ce type d'objet simple.
\subsection{Description des objets simples de $\LL$-$\mathbf{mod}$}

La validité du théorème de Jordan-Hölder dans la catégorie des modules  de dimension finie pour une algèbre de courants tordue $\LL$ donne espoir de bien décrire cette catégorie si ses objets simples et ses blocs d'extensions peuvent être bien décrits. La classification des objets simples de cette catégorie est le premier pas qui s'inscrit dans cette démarche. C'est d'ailleurs quelque chose de nécessaire pour une étude pertinente des blocs d'extensions de cette catégorie et de leur classification. À la base, deux modules simples qui sont isomorphes doivent être dans un même bloc d'extensions et la relation d'équivalence qui définit les classes d'équivalences donne que l'appartenance d'un module à une classe est déterminée à partir de ses facteurs simples. 

Dans cette première sous-section, il est montré comment décrire, peut-être grossièrement, un objet simple de la catégorie $\LL$-$\mathbf{mod}$ où $\LL$ est une algèbre de courants tordue. Rappelons que $\LL = (\g \otimes S)^\Gamma$ avec toutes les hypothèses de la définition donnée au tout début de la section 2.1.1.  Rappelons aussi que l'action du groupe $\Gamma$ sur $\g \otimes S$ induit une action de $\Gamma$ sur $\Max S$ de laquelle est déduite une action sur $S$ (voir la proposition \ref{1cocycle} et les explications qui la suivent).

Pour la première partie de cette sous-section, $(\rho,V)$ sera un $\LL$-module simple de dimension finie. D'abord, si $(\rho,V)$ est un $\LL$-module simple de dimension finie, alors la proposition \ref{quotientreductif} donne que $\LL/\ker \rho$ est une algèbre de Lie réductive avec un petit centre, i.e. : 
\begin{equation} \label{quotientreductifirrep}
\LL / \ker \rho \cong \Z \oplus \SSS = \Z \oplus \Big(\bigoplus_{i=1}^r \SSS_i\Big)
\end{equation}
où $\Z \cong Z(\LL / \ker \rho)$ est une algèbre de Lie abélienne de dimension inférieure ou égale à 1 et où $\SSS \cong [\LL / \ker \rho,\LL / \ker \rho]$ est une algèbre de Lie semi-simple de dimension finie. 

La première chose qui peut être dite du $\LL$-module $V$ est qu'il s'écrit de manière unique (à isomorphismes près) comme un produit tensoriel d'un $\Z$-module simple et d'un $\SSS$-module simple. Ceci est dû au corollaire \ref{sommelie}. De plus, l'algèbre de Lie $\Z$ est abélienne alors ses modules simples sont des modules de dimension 1. On peut donc écrire
\begin{equation} \label{premierechose}
(\rho,V) \cong (\rho_W \otimes \rho_E, W \otimes E)
\end{equation}
où $W$ est un $\Z$-module de dimension 1 et où $E$ est un $\SSS$-module de dimension finie. 

De même, on peut décomposer le $\SSS$-module qu'est $E$ en un produit tensoriel de $r$ modules simples ; un pour chaque $\SSS_i$. Pour chaque $i \in \{1,...\,,r\}$, notons le $\SSS_i$-module en question $(\rho_i,E_i)$. On peut donc écrire
\begin{equation} \label{modevprodtens}
(\rho_E,E) \cong (\rho_1 \otimes \cdots \otimes \rho_r, E_1 \otimes \cdots \otimes E_r)
\end{equation}

Ensuite, chacune des algèbres de Lie simples $\SSS_i$ qui apparaissent dans (\ref{quotientreductifirrep}) peuvent être vues comme le quotient de $\LL$ par un idéal maximal $\MM_i \in \Max \LL$ (voir la définition \ref{maxl}). Une étude approfondie de ces idéaux a permis de comprendre qu'il est possible de leur associer des idéaux maximaux $I_i \in \Max S^\Gamma$ d'une certaine façon bien précise (voir \cite{LauTCA}, lemme 2.7). Ce qu'il faut retenir, c'est que si l'on choisit des idéaux maximaux $M_i$ de $S$ qui sont au dessus des idéaux $I_i \in \Max S^\Gamma$ correspondants, alors selon la proposition 2.9 de \cite{LauTCA}, il existe un homomorphisme surjectif de l'algèbre de Lie $\g^{M_i}$ dans $\SSS_i$ pour chacun des indices $i$ ;
\begin{equation} \label{surjectionobscure}
\g^{M_i} \quad \twoheadrightarrow \quad \LL/\MM_i \; = \; \SSS_i
\end{equation}
où $\g^{M_i}$ est définie par (\ref{definitiongm}).

\begin{Rem} Puisque l'extension d'anneaux $S/S^\Gamma$ est entière (lemme \ref{extensionentiere}), choisir un idéal de $S$ au dessus d'un idéal maximal de $S^\Gamma$ est toujours possible. De plus, le fait que l'idéal de $S^\Gamma$ soit maximal implique que les idéaux de $S$ au dessus de lui sont maximaux. Ceci est montré dans \cite{Bourbakialgcomm} à la proposition 1 de la section 2.1.

De plus, si $J \in \Max S^\Gamma$, alors $\Gamma$ agit transitivement sur l'ensemble des idéaux maximaux de $S$ au dessus de $J$ et en particulier, cet ensemble est fini. Ceci est montré dans \cite{Bourbakialgcomm} au théorème 2 de la section 2.2.
\end{Rem}

\begin{Rem} En particulier, puisque pour chaque $i \in \{1,...\,,r\}$, l'homomorphisme \linebreak (\ref{surjectionobscure}) est surjectif, le « pullback » de chaque $\rho_i$ par la surjection correspondante fait des $\SSS_i$-modules simples $E_i$ des représentations irréductibles de $\g^{M_i}$.

Autrement dit, chaque $E_i$ est lui-même un $\LL$-module simple d'évaluation.
\end{Rem}

Chacun des $E_i$ est donc un $\LL$-module simple d'évaluation, mais qu'en est-il pour le module $E = \otimes_{i=1}^r E_i$ de la ligne (\ref{modevprodtens}) ? Le résultat s'avère également être vrai ; $E$ sera aussi un $\LL$-module simple d'évaluation. 

Supposons qu'un sous-ensemble d'indices $A \subseteq \{1,...\,,r\}$ est tel que pour tout $a,b \in A$, les idéaux maximaux $M_a \cap S^\Gamma = M_b \cap S^\Gamma = I_A \in \Max S^\Gamma$. Alors la proposition 2.9 de \cite{LauTCA} donne que pour tout choix de $M_a \in \Max S$ au dessus de $I_A$, il existe un homomorphisme surjectif
\begin{equation} \label{surjectiontresobscure}
\g^{M_a} \quad \twoheadrightarrow \quad \bigoplus_{c \in A} \LL/\MM_c \; = \; \bigoplus_{c \in A} \SSS_c
\end{equation}

\begin{Rem} Que le module $E = \otimes_{i=1}^r E_i$ soit un $\LL$-module simple d'évaluation signifie que $E$ est un objet simple de $\LL_{\ev}$-$\mathbf{mod}$.
\end{Rem}

En conséquence, on peut décomposer un $\LL$-module simple de dimension finie $(\rho,V)$ comme
$$ V \cong W \otimes E $$
où $W$ est un $\Z$-module de dimension 1 et $E$ est un $\left(\bigoplus_{i=1}^r \g^{M_i}\right)$-module simple de dimension finie avec un ensemble d'idéaux $\{M_i\}_{i=1}^r \subseteq \Max S$ qui proviennent de $r$ $\Gamma$-orbites distinctes.

Le $\LL$-module $V$ peut donc être interprété comme une représentation irréductible de l'algèbre de Lie $\Z \oplus \big(\bigoplus_{i=1}^r \g^{M_i}\big)$. Ce résultat est précisément l'énoncé du théorème 3.1 de \cite{LauTCA}.
%

Le théorème 3.2 dans \cite{LauTCA} montre ensuite que tous les produits tensoriels d'un $\LL$-module de dimension 1 avec un tel objet simple de $\LL_{\ev}$-$\mathbf{mod}$ sont des $\LL$-modules simples de dimension finie.

\begin{Rem} \label{equivalencemodulesimple} En conséquence, les notions suivantes apparaissent comme équivalentes :
\begin{itemize}
\item[\em(1)] Un objet simple de $\LL$-mod.
\item[\em(2)] Un élément de $(\LL/\LL')^*$ et un objet simple de $\LL_{\ev}$-$\mathbf{mod}$ qui fait intervenir un nombre fini d'idéaux de $\Max S$ issus de $\Gamma$-orbites toutes deux-à-deux distinctes.
\end{itemize}
\end{Rem}

Ceci suggère qu'afin d'étudier les modules simples de dimension finie de $\LL$, une démarche naturelle consisterait en l'étude des modules simples d'évaluation, puis celle de leur relation avec les modules de dimension 1.

\subsection{Classification des objets simples de $\LL$-$\mathbf{mod}$}

Dans cette sous-section, la classification des modules simples de dimension finie pour les algèbres de courants tordues est rapidement exposée. Les résultats de cette sous-section viennent eux-aussi de l'article \cite{LauTCA}. Le lecteur pourra y retrouver tous les détails incluant les démonstrations complètes qui sont volontairement omises dans ce document. La notion de module d'évaluation est particulièrement importante dans la classification de cette sous-section, tout comme dans la suite de ce travail.

Reprenons à partir des constats établis à la remarque \ref{equivalencemodulesimple}. Commençons par supposer que $(\rho_E,E)$ soit un objet simple de $\LL_{\ev}$-$\mathbf{mod}$ qui fait intervenir un nombre fini d'idéaux de $\Max S$ issus de $\Gamma$-orbites toutes distinctes entre elles. Tel qu'expliqué dans la sous-section précédente, on peut supposer qu'il existe un ensemble d'idéaux $\{M_i\}_{i=1}^r \subseteq \Max S$ tels que $M_i \cap S^\Gamma \neq M_j \cap S^\Gamma$ si $i \neq j$ et que 
$$ (\rho_E,E) \cong (\rho_1 \otimes \cdots \otimes \rho_r, E_1 \otimes \cdots \otimes E_r) $$
où chaque $E_i$ est un $\g^{M_i}$-module simple.

À la lumière de la remarque \ref{equivalencemodulesimple}, il est naturel d'identifier la représentation $E$ à la collection finie des modules $E_i$. Pour chaque $i$, $E_i$ est un $\g^{M_i}$-module. Pour détailler quoi que ce soit d'autre, il sera important de comprendre ce que représentent ces idéaux $M_i \in \Max S$ au regard des modules simples $E_i$. En outre, si l'information contenue dans la donnée des $E_i$ s'avère suffisante pour décrire $E$, il faudra savoir sous quelles conditions de telles informations donnent lieu à une même classe d'isomorphismes de $\LL$-modules d'évaluation. 

Avant de présenter le résultat qui répond essentiellement à ces questions, il sera utile de décrire efficacement l'information contenue dans la donnée d'un $\LL$-module simple d'évaluation $E$. Pour $M \in \Max S$, notons $\Rep(\g^M)$ l'ensemble des classes d'isomorphismes de $\g^M$-modules simples de dimension finie. C'est ainsi que nous considérerons la fibration
$$ \R \; = \bigsqcup_{M \in \Max S} \Rep(\g^M) \twoheadrightarrow \Max S $$

\begin{Rem} L'introduction d'une telle fibration peut être vue, dans ce contexte, comme l'introduction d'un simple outil pratique. Cette notion permet remarquablement bien d'exposer beaucoup d'information de façon très concise. D'ailleurs, les résultats qui suivent en témoignent.
\end{Rem}

Dans le cadre de cette démarche, il est naturel de considérer des sections de cette fibration $\R$. Posons donc 
$$ \secr = \{ \text{ sections de } \R \text{ } \} $$

Soit un $f \in \secr$. Alors, pour chaque $M \in \Max S$, $f$ lui fait correspondre une classe d'isomorphismes de $\g^M$-module simple de dimension finie. Aussi, on définit le support de $f$ par
$$ \supp f  = \{M \in \Max S \; | \; f(M) \neq [k_0] \in \Rep(\g^M)\}$$

\begin{Rem} \label{supraremarque} Tel que compris dans la sous-section précédente, l'information contenue dans la donnée d'un objet simple de $\LL_{\ev}$-$\mathbf{mod}$ comprend essentiellement :
\begin{itemize}
\item[\em(1)] Un nombre fini de modules simples, chacun étant un module simple pour une algèbre de la forme $\g^M$ (pour un $M \in \Max S$).
\item[\em(2)] Si l'on prend un de ces modules simples et que c'est un $\g^M$-module simple ($M \in \Max S$), alors il doit être le même $\LL$-module (à isomorphisme près) si on « le ré-interprète » comme un $\g^{{}^\gamma M}$-module ; et ce quel que soit $\gamma \in \Gamma$.
\end{itemize}
Il apparaît donc naturel de faire un rapprochement entre les objets simples de $\LL_{\ev}$-$\mathbf{mod}$ et des sections à support fini si ces dernières ont une certaine propriété d'invariance par rapport au groupe $\Gamma$ qui refléterait $(2)$.
\end{Rem}

Soit $M \in \Max S$ et soit $\gamma \in \Gamma$. Concrètement, le sens de « ré-interpréter » un $\g^M$-module simple $(\rho,V)$ comme un $\g^{{}^\gamma M}$-module est de prendre le « pullback » de l'action $\rho$ via l'isomorphisme d'algèbres de Lie 
\begin{equation} \label{isopourpullback}
\gamma^{-1}(M) : \; \g^{\,{}^\gamma M} \longrightarrow \; \g^M
\end{equation}
Le $\g^M$-module simple $(\rho,V)$ devient donc le $\g^{{}^\gamma M}$-module simple $\big(\rho \,\circ\, \gamma^{-1}(M),V\big)$. Voir \cite{LauTCA}, lemme 3.8 pour les détails concernant (\ref{isopourpullback}). 

\begin{Rem} Dans la même référence, le lemme 3.8 et les équations (3.5), (3.6) et (3.7) montre que pour « ré-interpréter » un $\g^{{}^{\gamma^{-1}}M}$-module comme un $\g^M$-module, il faut plutôt utiliser l'isomorphisme d'algèbres de Lie $\big(\gamma(M)\big)^{-1} : \g^M \rightarrow \g^{{}^{\gamma^{-1}}M}$ pour le « pullback » de l'action des modules.
\end{Rem}

Soit $f \in \secr$ une section de $\R$. Quelque soit $M \in \Max S$, on sait que $f(M)$ est une classe d'isomorphismes de $\g^{M}$-module simple alors c'est possible d'écrire
$$ f(M) = \left[(\rho_{f,M},V)\right] \in \Rep(\g^M) $$
Si $\gamma \in \Gamma$, on définit une nouvelle section ${}^\gamma f \in \secr$ par 
\begin{equation} \label{actionfgamma}
{}^\gamma f : \; M \longmapsto \left[\Big(\rho_{f,{}^{\gamma^{-1}}M} \circ \big(\gamma(M)\big)^{-1},V\Big)\right]
\end{equation}

\begin{Rem} Pour éviter d'avoir recours à cette écriture plutôt lourde, l'équation (\ref{actionfgamma}) sera dorénavant abrégée par l'écriture
$$ \big({}^\gamma f\big)(M) = \left[f({}^{\gamma^{-1}}M) \circ \big(\gamma(M)\big)^{-1}\right] $$ 
\end{Rem}

\begin{Rem} \label{actionsecr} La définition de ${}^\gamma f$, donnée par (\ref{actionfgamma}), est bien définie et est telle que l'application
\begin{align*}
&\Gamma \times \secr \longrightarrow \secr\\
&\;(\gamma,f) \mapsto {}^\gamma f
\end{align*}
est une action du groupe $\Gamma$ sur l'ensemble des sections $\secr$ de la fibration $\R$. C'est là l'objet du lemme 3.10 de \cite{LauTCA}.
\end{Rem}

Il n'est pas trop difficile de se convaincre que pour qu'une section $f \in \secr$ rende compte du point (2) de la remarque \ref{supraremarque}, il faut et il suffit d'avoir $f \in \secr^\Gamma$, c'est-à-dire que $f$ soit invariante sous l'action $\Gamma \curvearrowright \secr$ présentée à la remarque précédente.

Pour donner suite à l'ensemble de la remarque \ref{supraremarque}, voici les sections de $\R$ qui ont le plus d'intéret jusqu'ici. Ces sections méritent une définition.

\begin{Def} L'ensemble des sections de $\R$ dont le support est fini et qui sont $\Gamma$-invariantes est noté $\F^\Gamma$  
\end{Def}

La proposition qui suit donne enfin la classification des objets simples de $\LL_{\ev}$-$\mathbf{mod}$. Celle-ci ne devrait désormais pas être trop surprenante. Il s'agit de la proposition 3.13 de \cite{LauTCA}.

\begin{Prop}\label{classmodulesev} L'ensemble des classes d'isomorphismes des objets simples de $\LL_{\ev}$-$\mathbf{mod}$ est en bijection naturelle avec l'ensemble $\F^\Gamma$ des sections $\Gamma$-invariantes à support fini de la fibration $\R$.
\begin{equation*} 
\F^\Gamma \quad \stackrel{1:1}{\longleftrightarrow} \quad
\left\{\begin{array}{c} \text{Classes d'isomorphismes de}\\
\LL\text{-modules simples d'évaluation}\\ \text{de dimension finie} \end{array}\right\}
\end{equation*}

Explicitement, supposons que l'on commence avec une classe d'isomorphismes $[E]$ d'un objet simple de $\LL_{\ev}$-$\mathbf{mod}$ $E$. Pour représenter la classe $[E]$, on peut en prendre un représentant $E_1 \otimes \cdots \otimes E_r$ où pour chaque indice $i \in \{1,...\,,r\}$, $(\rho_i,E_i)$ est un $\g^{M_i}$-module simple pour un certain idéal $M_i \in \Max S$. Sans perdre de généralité, il peut être supposé que $\Gamma.M_i \neq \Gamma.M_j$ lorsque $i \neq j$.

La bijection associe à $[E]$ la section $f_{[E]} = f_E \in \F^\Gamma$ définie comme
\begin{align*}
f_E : \; &\Max S \longrightarrow \bigsqcup_{M \in \Max S} \Rep(\g^M) \quad = \quad \R\\
 &\quad\;\; M \longmapsto \left\{
\begin{array}{cl} \left[ \Big(\rho_i \circ \big(\gamma(M)\big)^{-1}, E_i \Big) \right] & \text{si } M = {}^\gamma M_i\\
\text{$[\,k_0\,]$} & \text{sinon}\\
\end{array} \right.
\end{align*}

Si par contre, on commence avec une section $\Gamma$-invariante à support fini $f \in \F^\Gamma$, nous pouvons écrire son support comme une union de $\Gamma$-orbites disjointes dans $\Max S$ :
$$ \supp f = \bigsqcup_{i=1}^r \Gamma.M_i $$

Enfin, la bijection associe à $f$ la classe d'isomorphismes de $\LL$-module simple
$$ [E_f] = \left[E_{f(M_1)} \otimes \cdots \otimes E_{f(M_r)}\right] $$
où $E_{f(M_i)}$ est un $g^{M_i}$-module simple issu de la classe d'isomorphismes $f(M_i) \in \Rep(g^{M_i})$ pour chacun des $i$.
\end{Prop}

Maintenant que les $\LL$-modules simples d'évaluation de dimension finie sont bien décrits par les sections $\F^\Gamma$, il ne reste qu'à s'assurer d'avoir une description des modules de la forme $W \otimes E$ où $W$ est un module de dimension 1 et $E$ en est un simple d'évaluation et de dimension finie. Naturellement, on doit vouloir penser à associer au $\LL$-module simple de dimension finie $W \otimes E$ le couple $(\LL/\LL')^* \times \F^\Gamma$ correspondant à $W$ et à $E$, mais il faut prendre soin d'éliminer les cas où la représentation donnée par deux couples est la même. Ainsi, certains couples sont laissés de côté et certains sont conservés. C'est dans cet esprit, que la classification complète des objets simples de $\LL$-$\mathbf{mod}$ est donnée dans \cite{LauTCA} par le théorème 3.14. Voici ce qui est montré :

\begin{Thm} \label{classmodules} L'ensemble des classes d'isomorphismes de $\LL$-modules simples de dimension finie est en bijection avec l'ensemble des couples $(\lambda,f) \in \LL^* \times \F^\Gamma$ satisfaisant
\begin{align*}
&\bullet \quad \lambda|_{[\LL,\LL]} = 0\\
&\bullet \quad \LL/\ker \rho_f \text{ est semi-simple}\\
&\bullet \quad \ker (\lambda+\rho_f) = \ker \lambda \cap \ker \rho_f
\end{align*}
où $\rho_f : \LL \rightarrow \End(E_f)$ est l'application décrivant l'action de $\LL$ sur $E_f$, un représentant de la classe d'isomorphismes de modules simples d'évaluation associée à $f$.
\end{Thm}

\begin{Rem} Cette proposition révèle en fait que dans chaque classe d'isomorphismes de $\LL$-module simple, il existe un unique représentant $(\lambda \otimes \rho_f, k_\lambda \otimes E_f)$ qui satisfait les trois conditions énumérées ci-haut.
\end{Rem}

Il est tout de même important de préciser le sens exact de la terminologie, possiblement ambigüe, utilisée dans la troisième condition de ce dernier théorème. A priori, le sens de $\lambda + \rho_f$ n'est pas défini. Pour fixer une notation, supposons que le couple $(\lambda,f) \in \LL^* \times \F^\Gamma$ correspond à la représentation irréductible $(\lambda \otimes \rho_f, k_\lambda \otimes E_f)$. Pour $a \otimes e \in k_\lambda \otimes E_f$, l'action d'un élément $\ell \in \LL$ est donné par
\begin{align*}
\ell.(a \otimes e) &= (\ell.a) \otimes e + a \otimes (\ell.e)\\
&= \big(\lambda(\ell)a\big) \otimes e + a \otimes \rho_f(\ell)e\\
&= a \otimes \lambda(\ell)e + a \otimes \rho_f(\ell)e\\
&= a \otimes \big(\lambda(\ell)\Id+\rho_f(\ell)\big)e
\end{align*}

Comme espace vectoriel, $k_\lambda \otimes E_f \cong E_f$ avec $a \otimes e = 1 \otimes ae \; \leftrightarrow \; ae$. Sous cette identification, $\lambda + \rho_f$ prend le sens de l'application
\begin{align*}
(\lambda\Id + \rho_f) : \; &\LL \longrightarrow \End_k(E_f) \\
&\; \ell \longmapsto \lambda(\ell)\Id + \rho_f(\ell)
\end{align*}

Plus simplement ces quelques lignes établissent que $(\lambda \otimes \rho_f, k_\lambda \otimes E_f) \cong (\lambda \Id + \rho_f,E_f)$ en tant que représentations de $\LL$.
\begin{Rem} Il est maintenant apparent que $\ker(\lambda \Id + \rho_f) = \ker(\lambda \otimes \rho_f)$. La troisième condition de la proposition \ref{classmodules} se lit donc  $ \ker(\lambda \otimes \rho_f) = \ker \lambda \cap \ker \rho_f $.
\end{Rem}

Maintenant, soit $(\lambda \otimes \rho_f,k_\lambda \otimes E_f) \, \leftrightarrow \, (\lambda, f) \in \LL^* \times \F^\Gamma$ satisfaisant aux trois conditions de la proposition \ref{classmodules}. Alors, il y a un isomorphisme naturel
\begin{align}
\label{restesChinois} \alpha : \quad &\LL / \ker(\lambda \otimes \rho_f) \; \longrightarrow \; \LL / \ker \lambda \oplus \LL/\ker \rho_f \\
&\ell + \ker(\lambda \otimes \rho_f) \; \mapsto \; \big(\ell + \ker \lambda,\ell + \ker \rho_f\big) \nonumber
\end{align}

En effet, il suffit de montrer que $\ker \lambda + \ker \rho_f = \LL$ et alors la démonstration classique du théorème des restes chinois tient et garantit le résultat. Si $\lambda = 0$ ou si $\rho_f = 0$, il n'y a rien à montrer. Supposons donc que $\lambda \neq 0$ et que $\rho_f \neq 0$.

Dans ce cas, on sait que $\dim_k \LL / \ker \lambda = 1$ et donc, 
$$ d = \dim_k \left(\frac{\ker \lambda + \ker \rho_f}{\ker \lambda}\right) \in \{0,1\} $$
Si $d = 0$, alors le second théorème d'isomorphisme permet d'écrire
$$ \dim_k \left(\frac{\ker \rho_f}{\ker \lambda \cap \ker \rho_f}\right) = 0 $$
Ceci entraine directement que $\ker \rho_f \subset \ker \lambda$. Il y a donc une surjection
$$ \LL / \ker \rho_f \twoheadrightarrow \LL / \ker \lambda $$
Puisque $0 \neq \LL / \ker \rho_f$ est une algèbre de Lie semi-simple, tous ses quotients sont semi-simples. C'est donc que $0 \neq \LL / \ker \lambda$ doit être à la fois une algèbre de Lie abélienne et semi-simple, mais cela est impossible.

Il suit que $d = 1$. Maintenant, puisque 
$$ \frac{\ker \lambda + \ker \rho_f}{\ker \lambda} \; \subseteq \; \frac{\LL}{\ker \lambda} $$ 
et que ce sont deux espaces de dimension 1, ils n'ont d'autre choix que d'être égaux. Ainsi, $\ker \lambda + \ker \rho_f = \LL$ et l'isomorphisme (\ref{restesChinois}) tient.

\begin{Rem} L'isomorphisme (\ref{restesChinois}) permet de donner une décomposition explicite de $\LL / \ker (\lambda \otimes \rho_f)$ en tant qu'algèbre de Lie réductive dans le cas général. Comme de raison, le centre de cette algèbre a une dimension d'au plus 1.

\end{Rem}

\section{Classification des blocs d'extensions de $\LL$-$\mathbf{mod}$}
\subsection{Résultats utiles}

C'est ici que commence la démarche pour la classification des blocs d'extensions de la catégorie $\LL$-$\mathbf{mod}$. L'essentiel des résultats sont repris du contexte des algèbres d'applications équivariantes, des travaux de E.Neher et A.Savage de l'Université d'Ottawa. Ils ont présenté et détaillé des travaux très complets dans leur article \cite{NSext}. 

Pour classifier les blocs d'extensions analogues dans le cas des algèbres d'applications équivariantes, ces professeurs se basaient sur leur précédente classification, avec P.Senesi, des classes d'isomorphismes de modules simples de dimension finie présentée dans l'article \cite{NSSreps}. Dans le cas des algèbres de courants tordues, une classification des modules simples de dimension finie a été montrée et détaillée par M.Lau dans son article \cite{LauTCA}. Cet article de M.Lau est, en quelque sorte, l'analogue de \cite{NSSreps} pour E.Neher et A.Savage dans leurs travaux de classification de blocs d'extensions.

La démarche proposée pour classifier les blocs d'extensions est de s'attarder d'abord aux $\LL$-modules simples d'évaluations, puis d'utiliser la classification des modules simples pour classifier les blocs d'extensions de tous les $\LL$-modules simples pour une algèbre de courants tordue $\LL$.

Cette première sous-section contient des résultats plus techniques ou qui, pris à part, n'apportent rien à la classification des blocs d'extensions. Néanmoins, ces résultats sont importants.  

\begin{Prop} \label{suitespectrale} Soit $\rho : \LL \rightarrow \End_k(V)$ une représentation de dimension finie de $\LL$ et soit $K \subseteq \ker \rho$ un idéal tel que $\LL/K$ soit une algèbre de Lie réductive de dimension finie. 

Écrivons $\LL/K = \Z \oplus \SSS$ où $\Z$ est une algèbre de Lie abélienne et $\SSS$ est une algèbre de Lie semi-simple. Finalement, notons $\tilde{\rho}$ l'action correspondante de $\LL/K$ sur $V$. Alors,

\begin{itemize}
\item[\em(1)] Si $\LL/K$ est semi-simple ou si $V$ est complètement réductible et qu'il existe $z \in \Z$ tel que $\tilde{\rho}(z)$ est inversible, alors
$$ \cohomo{\LL}{V} \cong \Hom_{\LL/K}(K/K',V) $$
\item[\em(2)] Si $\Z.V = 0$, alors il y a un isomorphisme induit $\cohomo{\LL/K}{V} \cong \Hom_k\big(\Z,V^\SSS\big)$ et une suite exacte
$$ 0 \rightarrow \cohomo{\LL/K}{V} \stackrel{inf}{\longrightarrow} \cohomo{\LL}{V} \stackrel{res}{\longrightarrow} \Hom_{\LL/K}(K/K',V) \longrightarrow \HHH^2(\LL/K;V) $$
\end{itemize}
\end{Prop}

Pour une preuve, voir ce qui est dit dans \cite{NSext}. La preuve donnée fait essentiellement appel à la suite spectrale de Hochschild-Serre en lui déduisant de information voulue de diverses façons. Voir au besoin le théorème 6 de \cite{HSlie}.

\begin{Rem} La suite exacte qui apparaît dans la partie (2) de la proposition précédente se lit aussi comme
$$ 0 \rightarrow \cohomo{\LL/K}{V} \stackrel{inf}{\longrightarrow} \cohomo{\LL}{V} \stackrel{res}{\longrightarrow} \big(\Hom_k(K/K',V)\big)^{\LL/K} \longrightarrow \HHH^2(\LL/K;V) $$
D'après l'isomorphisme (\ref{dertrivial}), on peut écrire que $\Hom_k(K/K',V) \cong \cohomo{K}{V}$ et ainsi, cette suite exacte peut aussi s'écrire
$$ 0 \rightarrow \cohomo{\LL/K}{V} \stackrel{inf}{\longrightarrow} \cohomo{\LL}{V} \stackrel{res}{\longrightarrow} \big(\cohomo{K}{V}\big)^{\LL/K} \longrightarrow \HHH^2(\LL/K;V) $$
Cette dernière forme est la plus révélatrice de la nature de la suite spectrale. Elle fait un peu penser à une suite exacte qui aurait été obtenue à partir de la suite exacte courte de $\LL$-modules $0 \rightarrow K \rightarrow \LL \rightarrow \LL/K \rightarrow 0$ par un foncteur contravariant exact à gauche bien que cela ne soit pas le cas.
\end{Rem}

\begin{Rem} Dans le contexte des modules pour les groupes, il existe des suites exactes similaires. Par exemple, le théorème 9.84 de \cite{Rotman} donne une suite exacte parfaitement analogue à celle dans la partie \em{$(2)$} \em de la proposition \ref{suitescindee}. Pour le voir, on peut se servir de l'isomorphisme (\ref{dertrivial}) avec $L = K$ et de la remarque \ref{invariantshom} avec $L = \LL/K$.
\end{Rem}

La proposition précédente est capitale dans toute cette étude des Ext-blocs. Elle joue un rôle de premier plan notamment dans les travaux de E.Neher et A.Savage sur les algèbres d'applications équivariantes. Ceci atteste, une fois parmi tant d'autres, que les notions de suites spectrales sont vraiment très performantes. 

\begin{Lem} \label{entiertypefini} Soit $A$ un anneau associatif, commutatif et unitaire et soit $B$ une $A$-algèbre. Soit $b \in B$. Alors, les propriétés suivantes sont équivalentes :
\begin{itemize}
\item[\em(i)] $b$ est entier sur $A$.
\item[\em(ii)] $A[b]$ est un $A$-module de type fini.
\end{itemize}
\end{Lem}

C'est un résultat bien connu. On en retrouve une preuve par exemple dans \cite{Lang}, voir l'équivalence démontrée à la page 334 (au tout début du chapitre VII).

\begin{Lem} \label{entiermax} Soit $A \subseteq B$ une extension entière d'anneaux. Si $P$ est un idéal maximal de $B$, alors $P \cap A$ est un idéal maximal de $A$.
\begin{pre} Soit $P \in \Max B$. Puisque $P$ est premier et que l'inclusion $A \subseteq B$ est un homomorphisme d'anneaux, on sait que $P \cap A$ est un idéal premier de $A$. Considérons l'homomorphisme d'anneaux 
\begin{align*}
\iota : \; &A / P \cap A \longrightarrow B / P\\
&a + P \cap A \mapsto a + P
\end{align*}
Montrons qu'il est injectif. Si $\iota (a + P \cap A) = a + P = 0 + P$, alors $a \in P$ et donc 
$$ a + P \cap A = 0 + P \cap A $$
Ainsi, $\iota$ est injectif. Montrons maintenant que $A / P \cap A \cong \iota(A / P \cap A)$ est un corps. Soit $a + P \in \iota(A / P \cap A) \bs \{0\}$. Comme $P \in \Max B$, il existe un unique $b + P \in B$ tel que 
$$ (a + P)(b + P) = ab + P = 1 + P $$
Il suffit donc de montrer que $b \in \iota(A / P \cap A)$. Comme l'extension d'anneaux $A \subseteq B$ est entière, il existe un $n \in \NN \bs \{0\}$ et des $r_i \in A$ tels que
$$ b^n + \sum_{i=1}^n r_i b^i = 0 $$
On peut donc écrire 
\begin{align*}
1 + P = (a^n + P)(b^n + P) &= - \sum_{i=1}^{n-1} (r_i + P)(a^n + P)(b^i + P)\\
&= - \sum_{i=1}^{n-1} (r_i + P)(a^{n-i} + P)\\
&= (a + P)\Big(- \sum_{i=1}^{n-1} (r_ia^{n - i - 1} + P)\Big)
\end{align*}
Par l'unicité de l'inverse dans un corps, on obtient
$$ b + P = - \sum_{i=1}^{n-1} (r_ia^{n - i - 1} + P) $$
Enfin, puisque $r_i, a \in \iota(A / P \cap A)$, $b + P$ aussi. C'est donc que $\iota(A / P \cap A) \cong A / P \cap A$ est un corps et que $A \cap P \in \Max A$. \qed
\end{pre}
\end{Lem}

\begin{Rem} En fait, l'énoncé du lemme précédent peut être vu comme un cas particulier d'un résultat un peu plus général qui dit que si $A \subseteq B$ est une extension entière d'anneaux et $Q$ est un idéal premier de $B$, alors $Q$ est maximal si et seulement si $Q \cap A$ est un idéal maximal de $A$.
\end{Rem}

\begin{Lem} Soit $A$, un anneau associatif, commutatif, unitaire et Noethérien. Alors, toute $A$-algèbre de type fini est Noethérienne.
\begin{pre} Soit $B = A[b_1,...\,,b_n]$ une $A$ algèbre de type fini. On peut alors écrire $$ B = A[x_1,...\,,x_n]/I $$
où $I$ est un idéal de $A[x_1,...\,,x_n]$.

Par le théorème de la base d'Hilbert, l'anneau $A[x_1,...\,,x_n]$ des polynômes en $n$ variables à coefficients dans $A$ est un anneau Noethérien. Ses quotients, dont fait partie $B$, le sont donc aussi. \qed
\end{pre}
\end{Lem}

\begin{Lem} \label{toutestnoetherien} $S^\Gamma$ est une $k$-algèbre de type fini et donc, $S^\Gamma$ est Noethérienne. De plus, $S$ et $\LL$ sont des $S^\Gamma$-modules de type fini et Noethériens.
\begin{pre} Le lemme \ref{extensionentiere} donne que $S^\Gamma \subseteq S$ est une extension entière d'anneaux.  

Montrons d'abord que $S^\Gamma$ est une $k$-algèbre de type fini. Par hypothèse, $S$ est une $k$-algèbre de type fini ; soient donc $s_1,...\,,s_t \in S$ tels que $S = k[s_1,...\,, s_t]$ (c'est-à-dire un ensemble de générateurs de $S$ comme $k$-algèbre). Soient $m_1(x),...\,,m_t(x) \in S^\Gamma [x]$ leurs polynômes minimaux respectifs.

Posons $C$ comme étant la $k$-algèbre engendrée par tous les coefficients des $m_i(x)$. Notons que $C$ est un anneau Noethérien par le lemme précédent. Par ailleurs, on a directement les propriétés suivantes :
\begin{align*}
&\bullet \quad C\text{ est une } k\text{-algèbre de type fini}\\
&\bullet \quad k \subseteq C \subseteq S^\Gamma \subseteq S\\
&\bullet \quad S / C \text{ est une extension entière d'anneaux}\\
&\bullet \quad S = C[s_1,...\,,s_t]
\end{align*}
Par induction, on peut montrer que $S$ est un $C$-module de type fini. Pour commencer, $C[s_1]$ est un $C$-module de type fini par le lemme \ref{entiertypefini}. Puis, supposons que $C[s_1,...,s_{n-1}]$ (avec $1 \leq n \leq t$) soit un $C$-module de type fini. Puisque $s_n$ est entier sur $C$, il est entier sur $C[s_1,...\,,s_{n-1}]$ et donc
$$ C[s_1,...\,, s_n] = \big(C[s_1,...\,,s_{n-1}]\big)[s_n] $$
est un $\big(C[s_1,...\,,s_{n-1}]\big)$-module de type fini. Puisque $\big(C[s_1,...\,,s_{n-1}]\big)$ est un $C$-module de type fini (c'est l'hypothèse d'induction), on conclut que $C[s_1,...\,, s_n]$ est un $C$-module de type fini. Si on prend $n = t$, ceci montre que $S$ est un $C$-module de type fini.

Puisque $C$ est un anneau Noethérien et que $S$ est un $C$-module de type fini, $S$ est un $C$-module Noethérien. En particulier, le sous-module $S^\Gamma \subseteq S$ est un $C$-module de type fini.

Soient $g_1,...\,,g_n \in S^\Gamma$ tels que $S^\Gamma = C.g_1 + \cdots + C.g_n$. Alors puisque $C$ est une $k$-algèbre de type fini, $S^\Gamma$ aussi est une $k$-algèbre de type fini. En particulier, le lemme précédent donne que $S^\Gamma$ est un anneau Noethérien.

Il reste à montrer que $S$ et $\LL$ sont des $S^\Gamma$-modules de type fini. Le fait qu'ils soient des $S^\Gamma$-modules Noethérien suivra automatiquement puisque $S^\Gamma$ est un anneau Noethérien.

Les $s_i$ ont été choisis de telle sorte que $S = k[s_1,...\,,s_t]$ en tant que $k$-algèbre. En conséquence, on a certainement $S = S^\Gamma[s_1,...\,,s_t]$. La même preuve par induction que celle utilisée un peu plus haut montre que $S$ est un $S^\Gamma$-module de type fini.

Le fait que $\LL = (\g \otimes S)^\Gamma$ soit un $S^\Gamma$-module est établi par la proposition \ref{actrmod}. Il reste à montrer que c'est un module de type fini. Fixons $\{e_j\}_{j=1}^n$ une $k$-base de la $k$-algèbre de Lie $\g$ de telle sorte qu'on puisse écrire l'égalité de $S$-modules
$$ \g \otimes S = \bigoplus_{j=1}^n (ke_j \otimes S) $$
Ainsi, $\g \otimes S \cong S^{\,\oplus\, n}$ comme $S$-modules. Puisque l'action de $S^\Gamma$ sur $\g \otimes S$ est également donnée par la multiplication sur le $S$, on a en fait que $\g \otimes S \cong S^{\,\oplus\, n}$ en tant que $S^\Gamma$-modules. Comme $S$ est un $S^\Gamma$-module de type fini, $\g \otimes S$ l'est également. Enfin, $S^\Gamma$ est un anneau Noethérien alors $\g \otimes S$ est un $S^\Gamma$-module Noethérien. Ainsi, $\LL$ qui est un sous-module de $\g \otimes S$ est un module de type fini. En somme, c'est que $\LL$ est un $S^\Gamma$-module de type fini et Noethérien. \qed
\end{pre}
\end{Lem}

\begin{Lem} \label{lemmekr} Soit $\{M_i\}_{i=1}^r \subseteq \Max S$ un ensemble d'idéaux maximaux de $S$ issus de $r$ $\Gamma$-orbites distinctes. Posons 
\begin{align*}
I = \{s \in S \, | \, s(^\gamma M_i) = 0 \text{ pour tout } \gamma \in \Gamma \text{ et pour chaque } i\} \\
\end{align*}
Alors, $I^\Gamma = I \cap S^\Gamma$ est un idéal de $S^\Gamma$ et de plus, en tant que $k$-algèbres, il y a un isomorphisme naturel 
$$ S^\Gamma / I^\Gamma \cong k^r $$
\begin{pre} Par le lemme \ref{entiermax}, les $M_i \cap S^\Gamma$ sont tous des idéaux maximaux de $S^\Gamma$. Par notre hypothèse, ces idéaux sont deux-à-deux distincts et donc, le théorème des restes chinois donne un homomorphisme surjectif d'anneaux unitaires
\begin{align*} 
f : \; S^\Gamma &/ I^\Gamma \longrightarrow \bigoplus_{i=1}^r S^\Gamma / M_i \cap S^\Gamma \quad \cong \quad k^r\\
&\ov{s} \longmapsto \big(s(M_1 \cap S^\Gamma), ... \,, s(M_r \cap S^\Gamma)\big) 
\end{align*}
Cette application est bien définie puisque si $x \in I^\Gamma \subseteq S$, alors $x(M_i) = 0 \in S/M_i \cong k$ pour tout $i$ et donc 
$$ x(M_i \cap S^\Gamma) = 0 \in S^\Gamma / M_i \cap S^\Gamma \; \hookrightarrow \; S / M_i \cong k $$

Pour s'assurer du résultat, il suffit de montrer que $f$ est injective. Soit $\ov{x} \in \ker f$. Pour chaque $i$, on a alors $x(M_i \cap S^\Gamma) = 0$ ; donc $x \in M_i$. Soit $\gamma \in \Gamma$. Alors, ${}^\gamma x = x$ et donc $x \in {}^\gamma M_i$. En somme, on a
$$ x({}^\gamma M_i) = 0 \quad \text{pour tout }\gamma \in \Gamma \text{ et pour tout } i $$ 
C'est donc que $x \in I^\Gamma$ et que $\ov{x} = 0 + I^\Gamma$. Il suit que $f$ est injective et est un isomorphisme. \qed
\end{pre}
\end{Lem}

\begin{Prop} \label{structureKab} Soit $\{M_i\}_{i=1}^r \subseteq \Max S$ un ensemble d'idéaux maximaux de $S$ issus de $r$ $\Gamma$-orbites distinctes. Considérons l'application d'évaluation 
$$ \ev_\mathbf{M} = \oplus_{i=1}^r \ev_{M_i} : \; \LL \twoheadrightarrow \bigoplus_{i=1}^r \g^{M_i} $$
Posons 
\begin{align*}
K &= \ker (\ev_\mathbf{M}) \\
I &= \{s \in S \, | \, s(^\gamma M_i) = 0 \text{ pour tout } \gamma \in \Gamma \text{ et pour chaque } i\} \\
&= \bigcap_{\gamma \in \Gamma, \, i \in \{1,...,r\}}^{\,} {}^\gamma M_i \\
I^\Gamma &= I \cap S^\Gamma \\
&= \bigcap_{i=1}^r (M_i \cap S^\Gamma) \\
N &= \{n \in K \, | \, I^\Gamma.\,n \subseteq K'\} 
\end{align*}
Alors, $K$, $K'$ et $N$ sont des $S^\Gamma$-idéaux de $\LL$ tels que $K' \unlhd N \unlhd K$. L'action (adjointe) de $\LL$ induit des actions de $\big(\bigoplus_{i=1}^r \g^{M_i}\big) \cong \LL / K$ sur les quotients $N/K'$, $K/K'$ et $K/N$. 

Sous ces hypothèses, il est vrai que comme des $\LL / K$-modules : 

\begin{itemize}
\item[\em(1)] $K/N$ est trivial.
\item[\em(2)] $N/K' = \bigoplus_{i=1}^r T_i$ o\`u les $\LL/K$-modules $T_i$ peuvent être vus comme des $g^{M_i}$-modules de dimension finie avec la propri\'et\'e que $i \neq j \; \Rightarrow \; g^{M_j}.T_i = \{0\}$.

En particulier, $ \dim_k N/K' < \infty$.
\end{itemize}
\begin{pre} On sait déjà que $K$ et $K'$ sont des $k$-idéaux de $\LL$. Montrons d'abord que $N$ est aussi un idéal de $\LL$.

Pour commencer, $0 \in N$ puisque $I^\Gamma.0 = \{0\} \subseteq K'$. Ensuite, si $a, b \in N$, 
$$ I^\Gamma.(a+b) \subseteq I^\Gamma.a + I^\Gamma.b \subseteq K' + K' \subseteq K' $$
Aussi, si $a \in N$, alors $-a$ l'est aussi car la multiplication par $-1 \in k$ commute avec l'action de $I^\Gamma \subseteq S^\Gamma$ sur $N \subseteq \LL$. Finalement, par la proposition \ref{actrmod}, $\LL$ est une $S^\Gamma$-algèbre de Lie et le crochet de Lie est $S^\Gamma$-linéaire. On peut donc écrire
$$ I^\Gamma.[\LL,N] = [\LL,I^\Gamma.N] \subseteq [\LL,K'] \subseteq K' $$
Ceci montre que $\LL.N = [\LL,N] \subseteq N$. Ainsi, $N$ est un $k$-idéal de $\LL$.

Montrons ensuite que $N$, $K$ et $K'$ sont des $S^\Gamma$-modules. Commençons par $K$. Soit $a \in S^\Gamma$ et soit $\sum_i x_i \otimes s_i \in K$, alors pour chacun des $j \in \{1,...\,,r\}$, on a que
\begin{align*}
\ev_{M_j}\left(a.\sum_i x_i \otimes s_i\right) &= \ev_{M_j}\left(\sum_i x_i \otimes as_i\right)\\
&= \sum_i a(M_j)s_i(M_j)x_i\\
&= a(M_j)\sum_i s_i(M_j)x_i\\
&= a(M_j) \ev_{M_j}\left(a.\sum_i x_i \otimes s_i\right)\\
&= a(M_j) \cdot 0\\
&= 0
\end{align*}
Ainsi, $K$ est un $S^\Gamma$-module. Pour ce qui est de $K'$, la $S^\Gamma$-linéarité du crochet de Lie permet d'écrire 
$$ S^\Gamma.K' = S^\Gamma.[K,K] = [S^\Gamma.K,K] \subseteq [K,K] = K' $$
Il ne reste que le cas de $N$. Soit $a \in S^\Gamma$ et soit $n \in N$. Alors puisque $I^\Gamma$ est un idéal de $S^\Gamma$, on peut écrire
\begin{align*}
I^\Gamma.(a.n) &= (I^\Gamma a).n\\
&\subseteq I^\Gamma.n\\
&\subseteq K'
\end{align*}
Maintenant, on a évidemment les inclusion $K' \subseteq K$ et $N \subseteq K$ alors il faut justifier $K' \subseteq N$. Cette dernière inclusion est vérifiée puisque 
$$ I^\Gamma.K' \subseteq S^\Gamma.K' \subseteq K' $$
Il est donc établi que $K$, $N$ et $K'$ sont des $S^\Gamma$-idéaux de $\LL$. Les quotients $K/N$, $N/K'$ et $K/K'$ sont donc des $\LL$-modules qui sont aussi invariants sous l'action de $S^\Gamma$. Il est facile de voir que $K \subseteq \LL$ agit comme 0 sur ces trois modules quotients. En effet, si $\ov{k} \in K/N$, on a
$$ K.\ov{k} = \ov{[K,k]} \in K/N $$
Puisque $K' \subseteq N$, on a $K.\ov{k} = 0 \in K/N$. La même démarche fonctionne directement pour $K/K'$. Quant à $N/K'$, prenons $\ov{n} \in N/K'$ et alors, on peut écrire
$$ K.\ov{n} = \ov{[K,n]} \in N/K' $$
Puisque $n \in N \subseteq K$, on a $[K,n] \subseteq K'$ et donc que $K.\ov{n} = 0 \in N/K'$. Ainsi, ces trois quotients sont des modules pour l'algèbre de Lie $\LL/K \cong \big(\bigoplus_{i=1}^r \g^{M_i}\big)$.

Pour terminer cette preuve, il faut vérifier les points (1) et (2) par rapport au structures de modules. Montrons d'abord (1).

Soit $\ell + K \in \LL/K$ et soit $y + N \in K/N$. Alors
\begin{equation} \label{ksurntrivial}
(\ell + K).(y + N) = [\ell, y] + N \in K/N
\end{equation}
Il suffit donc de montrer que $[\ell,y] \in N$. Par un calcul précédent dans cette démonstration, on voit que si $s \in S^\Gamma$ et $\ell \in \LL$, alors pour chaque $j \in \{1,...\,,r\}$, on a
\begin{equation} \label{igammatuelsurk}
\ev_{M_j}(s.\ell) = s(M_j)\ev_{M_j}(\ell)
\end{equation}
Si on prend $s \in I^\Gamma$ dans l'équation précédente, on obtient $\ev_{M_j}(s.\ell)$ quel que soit $j$. Ainsi, $I^\Gamma.\LL \subseteq K$. On peut donc écrire 
$$ I^\Gamma.[\LL,N] = [I^\Gamma.\LL,N] \subseteq K' \subseteq N $$
Ceci prouve, d'après l'équation (\ref{ksurntrivial}), que $\LL/K$ agit comme 0 sur le module $K/N$. Il ne reste plus qu'à montrer le point (2) de la proposition.

Par définition de $N$, $I^\Gamma$ agit comme 0 sur le module $N/K'$. Aussi, l'équation (\ref{igammatuelsurk}) montre que $I^\Gamma$ agit comme 0 sur l'algèbre de Lie $\LL/K$. En conséquence, $\LL/K$ et $N/K'$ ont des structures de $(S^\Gamma/I^\Gamma)$-module, c'est-à-dire de $k^r$-module d'après le lemme \ref{lemmekr}. De plus, on a que $N/K$ est un module pour la $k^r$-algèbre de Lie qu'est $\LL/K \cong \bigoplus_{i=1}^r \g^{M_i}$.

En vertu du lemme 3.1 de \cite{NSext}, $N/K'$ peut s'écrire comme
\begin{equation} \label{decompositionnsurkprime}
N/K' = \bigoplus_{i=1}^r T_i
\end{equation}
où $\g^{M_j}.T_i = 0$ si $j \neq i$.

Enfin, puisque $N$ est un sous-$S^\Gamma$-module du $S^\Gamma$-module Noethérien $\LL$ (voir le lemme \ref{toutestnoetherien}), il est de type fini en tant que $S^\Gamma$-module. Ceci implique que $N/K'$ est un $k^r$-module de type fini et donc, que les $T_i$ de (\ref{decompositionnsurkprime}) sont tous des $k$-espaces vectoriels de dimension finie. On en conclut que $\dim_k N/K' < \infty$. \qed
\end{pre}
\end{Prop}

\begin{Rem} Si $r = 1$ et que $M = M_1$, alors c'est que $K = \ker (\ev_M)$ et $\LL / K \cong g^M$. Le résultat est alors simplement que comme $\g^M$-modules, $K/N$ est trivial et $N/K'$ est de dimension finie.

Notons aussi que dans ce cas, on a aussi $K = (\g \otimes I)^\Gamma$ et que les trois modules quotients ne dépendent vraiment que de $\Gamma.M \subseteq \Max S$.
\end{Rem}

\subsection{Classification}

C'est dans cette sous-section qu'est enfin présentée la classification des blocs d'extensions de $\LL$-$\mathbf{mod}$ lorsque $\LL$ satisfait une condition technique. Cette sous-section débute aussitôt avec quelques résultats qui apparaîtront soit rassurants ou soit révélateurs aux yeux du lecteur attentif.
  
\begin{Prop} \label{suppdisjoints} Soient $E$ et $F$ des $\LL$-modules simples d'\'evaluation de dimension finie dont les supports sont disjoints. Alors,
$$ \Ext_\LL^1(E,F) \cong \{0\} $$
\begin{pre} Soit $\{M_i\}_{i=1}^r \subseteq \Max S$ un ensemble complet de représentants des $\Gamma$-orbites de $\supp E \, \sqcup \; \supp F$. Sans perdre de généralité, choisissons les indices de telle sorte que \linebreak $M_1, ... , M_s \in \supp E$ et $M_{s+1}, ..., M_r \in \supp F$ et écrivons
\begin{align*}
\g_{E} &= \bigoplus_{i=1}^s g^{M_i} \\
\g_{F} &= \bigoplus_{i=s+1}^r g^{M_i} \\
\g_\text{Total} &= \g_E \oplus \g_F = \bigoplus_{i=1}^r g^{M_i}
\end{align*}
Posons $\ev_\mathbf{M} = \bigoplus_{i=1}^r \ev_{M_i} : \LL \twoheadrightarrow \g_\text{Total}$ et appelons $K$ son noyau. 

Posons aussi $U = \Hom_k(E,F) \cong E^* \otimes F$ de telle sorte que $\Ext_\LL^1(E,F) \cong \cohomo{\LL}{U}$ d'après la proposition \ref{cohomoext}. 

Notons que $U$ est un $\LL$-module simple. Le module $E^*$ est un $\LL$-module simple (voir le corollaire \ref{dualsimple}) qui est d'évaluation. En fait, avant d'être un $\LL$-module, $E$ est un $\g_E$-module et donc, $E^*$ est aussi un $\g_E$-module. Aussi, on a sans difficulté l'égalité $\supp E^* = \supp E$. Enfin, puisque les supports de $E^*$ et $F$ sont disjoints, que $E^*$ est un $\g_E$-module simple et que $F$ est un $\g_F$-module simple, le corollaire \ref{sommelie} donne que $E^* \otimes F$ est un module simple pour l'algèbre de Lie $\g_E \oplus \g_F = \g_\text{Total}$.

Les hypothèses et le choix de notation entrainent directement que $E$, $F$, puis $U$ peuvent tous être vus comme des $\g_\text{Total}$-modules pour lesquels $\g_F.E = \{0\}$, $\g_E.F = \{0\}$ et aussi
\begin{equation} \label{propriosu}
g_E.U \neq \{0\} \qquad \text{et} \qquad g_F.U \neq \{0\}
\end{equation}

Le but de cette démonstration est donc de montrer que $\cohomo{\LL}{U} \cong \{0\}$. Dans un premier temps, justifions que pour arriver au résultat, il est suffisant d'avoir 
\begin{equation} \label{suffisant}
\Hom_{\g_\text{Total}}(K/K',U) \cong \{0\}
\end{equation} 
Cette condition sera vérifiée dans un second temps. Notons que $\LL/K \cong \g_\text{Total}$ est une algèbre de Lie réductive parce que c'est une somme directe d'algèbres de Lie réductives (voir la proposition \ref{gmreductive}). Écrivons donc 
$$ \g_\text{Total} = \Z \oplus \SSS $$
où $\Z$ est une algèbre de Lie abélienne et $\SSS$ est une algèbre de Lie semi-simple. Si $\Z.U = \{0\}$, alors la partie (2) de la proposition \ref{suitespectrale} est valide. Dans la présente situation, le premier terme de la suite exacte $\Hom_k(\Z,U^\SSS)$ vaut $\{0\}$. Effectivement, $U$ est un $\g_\text{Total}$-module simple non-trivial et on suppose que $\Z.U = \{0\}$, ainsi, $U^\SSS = U^{\g_\text{Total}} = \{0\}$. Ceci donne directement $\Hom_k(\Z,U^\SSS) = \Hom_k\big(\Z,\{0\}\big) \cong \{0\}$. De cette information, la suite exacte donnée dans la proposition \ref{suitespectrale} se lit en fait 
$$ 0 \longrightarrow \cohomo{\LL}{U} \stackrel{res}{\longrightarrow} \Hom_{\g_\text{Total}}(K/K',U) $$
Si par contre $\Z.U \neq \{0\}$, comme le module $U$ est simple (donc complètement réductible), il existe $z \in \Z$ tel que $\tilde{\rho}(z) \in \End_k(U) \bs \{0\}$. Puisque $\Z = Z(\LL/K)$, $\tilde{\rho}(z) \in \End_{\g_\text{Total}}(U)$. Puis, $U$ est simple alors le lemme de Schur (lemme \ref{schur}) donne que $\tilde{\rho}(z)$ est un isomorphisme ; en particulier, c'est une application linéaire inversible. Les conditions sont donc remplies pour pouvoir utiliser la partie (1) de la proposition \ref{suitespectrale} qui indique que 
$$ \cohomo{\LL}{U} \cong \Hom_{\g_\text{Total}}(K/K',U) $$
Dans les deux cas, c'est manifeste que la condition (\ref{suffisant}) implique la conclusion voulue, à savoir $\cohomo{\LL}{U} \cong \{0\}$.

Maintenant, il s'agit de vérifier (\ref{suffisant}). Soit $f \in \Hom_{\g_\text{Total}}(K/K',U)$ alors si $N/K' \subseteq \ker f$, on a un diagramme commutatif de $\g_\text{Total}$-modules 
\begin{center}\begin{tikzpicture}
y\matrix (m) [matrix of math nodes, row sep=2.5em,column sep=2.5em,text height=1.5ex,text depth=0.25ex, nodes in empty cells]
{ K/K' & f(K/K') \subseteq U \\ K/N &  \\};
\draw (m-1-1) edge[->>,thick] node[auto] {\footnotesize{$f$}} (m-1-2)
	  (m-1-1) edge[->>,thick] (m-2-1)
	  (m-2-1) edge[->>,thick] node[below] {\footnotesize{$\tilde{f}$}} (m-1-2);
\end{tikzpicture} \end{center}
D'après la proposition \ref{structureKab}, $K/N$ est un $\g_\text{Total}$-module trivial et donc, 
$$ f(K/K') \subseteq U^{\g_\text{Total}} $$
Puisque $U$ est un module simple non-trivial, $U^{\g_\text{Total}} = \{0\}$ et ainsi, $\tilde{f} = 0 \; \Rightarrow f = 0$. 

Supposons enfin que $N/K' \nsubseteq \ker f$. Alors la proposition \ref{structureKab} donne une décomposition de $N/K'$ en une somme directe de sous-modules $T_i$ qui satisfont à la propriété 
$$ i \neq j \; \Rightarrow \; g^{M_j}.T_i = \{0\} $$ 
Puisque $f|_{N/K'} \neq 0$, il faut admettre que $f|_{T_{i_0}} \neq 0$ pour un certain indice $i_0$. Ensuite, $U$ est simple donc $f(T_{i_0}) = U$ et ainsi, ce module simple qu'est $U$ apparaît comme un quotient du module $T_{i_0}$. Il suit directement que $\g_E.U = 0$ ou $\g_F.U = 0$, mais dans les deux cas, la conclusion obtenue contredit directement les propriétés préalablement établies et présentées à la ligne (\ref{propriosu}). 

Par ce qui a été expliqué ci-haut, la preuve est terminée et l'on peut conclure tour à tour que $\Hom_{\g_\text{Total}}(K/K',U) \cong \{0\}$ et puis que $\Ext_\LL^1(E,F) \cong \{0\}$. \qed
\end{pre}
\end{Prop}

\begin{Thm} \label{exteval} Soient $E$ et $F$ deux $\LL$-modules simples d'évaluation de dimension finie. Fixons $\{M_i\}_{i=1}^r \subseteq \Max S$ un ensemble complet de représentants des $\Gamma$-orbites de \linebreak $\supp E \cup \supp F$. Supposons donc que $E = \bigotimes_{i=1}^r E_i$ et  $F = \bigotimes_{i=1}^r F_i$ où les $E_i$ et $F_i$ sont des $\g^{M_i}$-modules simples. 

Soient $f_E$ et $f_F$ les sections $\Gamma$-invariantes à support fini de $\Max S$ dans $\mathcal{R}$ correspondant aux classes d'isomorphismes de $E$ et $F$ respectivement. Alors,

\begin{itemize}
\item[(1)] Si $f_E$ et $f_F$ diff\`erent sur plus d'une $\Gamma$-orbite de $\Max S$, alors 
$$ \Ext_\LL^1(E,F) \cong \{0\} $$
\item[(2)] Si $f_E$ et $f_F$ ne diff\`erent que sur la $\Gamma$-orbite de $M_j \in \Max S$, alors 
$$ \Ext_\LL^1(E,F) \cong \Ext_\LL^1(E_j,F_j) $$
\item[(3)] Si $f_E = f_F$, alors 
$$ \big((\LL/\LL')^*\big)^{\oplus\,r-1} \oplus \Ext_\LL^1(E,F) \cong \bigoplus_{i=1}^r \Ext_\LL^1(E_i,F_i) $$
\end{itemize}
\begin{pre} Fixons d'abord $\g_\text{Total} = \bigoplus_{i=1}^r g^{M_i}$. 

Pour commencer, supposons que $f_E$ et $f_F$ diffèrent sur au moins une $\Gamma$-orbite et fixons $j$ tel que $f_E(M_j) \neq f_F(M_j)$. Alors par (\ref{switchext}), on peut écrire
$$ \Ext_\LL^1(E,F) \cong \Ext_\LL^1\Big(\bigotimes_{i \neq j} \big(E_i \otimes F_i^*\big),E_j^* \otimes F_j\Big) $$
Maintenant, puisqu'un produit tensoriel (fini) de représentations complètement réductibles résulte en une représentation complètement réductible, chacun des deux $\g_\text{Total}$-modules qui apparaissent dans la précédente expression sont semi-simples et ainsi, ils s'écrivent comme des sommes finies de modules simples. 
\begin{align*}
\bigotimes_{i \neq j} \big(E_i \otimes F_i^*) \cong \bigoplus_{a=1}^m X_a && E_j^* \otimes F_j \cong \bigoplus_{b=1}^n Y_b
\end{align*}
où les $X_a$ sont des $\big(\bigoplus_{i \neq j} \g^{M_i}\big)$-modules simples et où les $Y_b$ sont des $\g^{M_j}$-modules simples. Comme le « bifoncteur » qu'est en fait $\Ext_\LL^1(-,-)$ commute avec les sommes directes finie dans ses deux composantes (voir \cite{Rotman}, propositions 7.21 et 7.22 ), on peut écrire
\begin{equation} \label{extdiffere}
\Ext_\LL^1(E,F) \cong \bigoplus_{a=1}^m\bigoplus_{b=1}^n\Ext_\LL^1(X_a,Y_b)
\end{equation}
Ensuite, l'hypothèse de $f_E(M_j) \neq f_F(M_j)$ donne que $\dim_k(E_j^* \otimes F_j)^L = 0$ (voir le corollaire \ref{dimhom}). Ainsi, aucun des modules $Y_b$ n'est trivial.

Supposons qu'il existe au moins un deuxième indice $k \neq j$ tel que $f(M_k) \neq g(M_k)$. Sachant que les $X_a$ sont des produits tensoriels de $\g^{M_i}$-modules simples (pour $i \neq j$), cette hypothèse combinée au corollaire \ref{dimhom} assure que $X_a$ n'est également jamais trivial (quel que soit $a$). Enfin, le fait qu'aucun $X_a$ ou $Y_b$ ne soient trivial donne immédiatement que 
$$ \supp X_a \cap \supp Y_b = \emptyset $$ en tant que $\g_\text{Total}$-modules. Ainsi, la proposition \ref{suppdisjoints} est applicable et on conclut que $\Ext_\LL^1(X_a,Y_b) \cong \{0\}$ quels que soient $a$ et $b$ et l'isomorphisme (\ref{extdiffere}) s'écrit 
$$ \Ext_\LL^1(E,F) \cong \{0\} $$

Autrement, si $f_E$ et $f_F$ ne diffèrent que sur la $\Gamma$-orbite de $M_j$, c'est que tous les $E_i \cong F_i$ si $i \neq j$. Le corollaire \ref{dimhom} donne encore que $k_0$ n'apparaît qu'une et une seule fois dans chaque $E_i \otimes F_i^*$ et ainsi, il existe un et un seul $X_a$ qui soit $k_0$, le module trivial. À l'aide de cette information et du raisonnement développé à la fin du paragraphe précédent, on voit que l'isomorphisme (\ref{extdiffere}) s'écrit
\begin{align*}
\Ext_\LL^1(E,F) &\cong \bigoplus_{b=1}^n\Ext_\LL^1(k_0,Y_b)\\
&\cong \Ext_\LL^1\Big(k_0,\bigoplus_{b=1}^n Y_b\Big)\\
&\cong \Ext_\LL^1(k_0,E_j^* \otimes F_j)\\
&\cong \Ext_\LL^1(E_j,F_j)
\end{align*}

Il ne reste désormais qu'à traiter du cas où $f_E = f_F$. Dans ce cas, le corollaire \ref{dimhom} donne que pour chaque $i \in \{1,...,r\}$, on ait des décompositions en somme de modules simples de la forme
$$ E_i \otimes F_i^* \cong k_0 \oplus \Big(\bigoplus_{s_i=1}^{n_i} C_{s_i}\Big) $$
où les $C_i^s$ sont tous des $g^M_i$-modules simples non-triviaux. Pour alléger la notation de ce qui suit, écrivons $N_i = \bigoplus_{s_i=1}^{n_i} C_{s_i}$. De cette façon, on peut écrire
\begin{align*}
\Ext_\LL^1(E,F) &\cong \Ext_\LL^1(E \otimes F^*,k_0)\\
&\cong \Ext_\LL^1\Big(\bigotimes_{i=1}^r (E_i \otimes F_i^*),k_0\Big)\\
&\cong \Ext_\LL^1\Big(\bigotimes_{i=1}^r \big(k_0 \oplus N_i\big),k_0\Big)
\end{align*}
Ensuite, on peut développer $\bigotimes_{i=1}^r \big(k_0 \oplus N_i\big)$ comme suit :
\begin{align*}
\textstyle{\bigotimes_{i=1}^r} \big(k_0 \oplus N_i\big) &\cong k_0 \oplus N_1 \oplus \cdots \oplus N_r\\
&\quad\; \oplus (N_1 \otimes N_2) \oplus (N_1 \otimes N_3) \oplus \cdots \oplus (N_1 \otimes N_r) \oplus \cdots \oplus (N_{r-1} \otimes N_r)\\
&\quad\; \oplus (N_1 \otimes N_2 \otimes N_3) \oplus (N_1 \otimes N_2 \otimes N_4) \oplus \cdots \oplus (N_{r-2} \otimes N_{r-1} \otimes N_r)\\
&\quad\; \oplus \cdots \oplus (N_1 \otimes \cdots \otimes N_r)
\end{align*}
C'est donc que $\Ext_\LL^1(E,F)$ est isomorphe à la somme directe de $\Ext_\LL^1(k_0,k_0)$, des $\Ext_\LL^1(N_i,k_0)$ et de termes $\Ext_\LL^1(N_{i_1} \otimes \cdots \otimes N_{i_c},k_0)$ où $c \in \NN \bs \{0,1\}$ et où $\#\{i_k\}_{k=1}^c = c$. Pour chaque terme de la dernière forme, on peut écrire
$$ \Ext_\LL^1(N_{i_1} \otimes \cdots \otimes N_{i_c},k_0) \cong \Ext_\LL^1(N_{i_1} \otimes \cdots \otimes N_{i_{c-1}},N_{i_c}^*) $$
où chacun des $N_{i_k}$ est un $\g^{M_{i_k}}$-module semi-simple avec chacune des composantes simples qui est non-trivial. C'est donc que comme des $\LL$-modules, $N_{i_1} \otimes \cdots \otimes N_{i_{c-1}}$ et $N_{i_c}^*$ sont des sommes de $\LL$-modules simples d'évaluation dont les supports des sont disjoints. Ainsi, la proposition \ref{suppdisjoints} donne 
$$ \Ext_\LL^1(N_{i_1} \otimes \cdots N_{i_c},k_0) \cong \Ext_\LL^1(N_{i_1} \otimes \cdots \otimes N_{i_{c-1}},N_{i_c}^*) \cong \{0\} $$
Puis, en combinant toutes les informations obtenues jusqu'ici, on a un isomorphisme
\begin{equation} \label{eurk}
\Ext_\LL^1(E,F) \cong \Ext_\LL^1(k_0,k_0) \oplus \Big(\bigoplus_{i=1}^r \Ext_\LL^1(N_i,k_0)\Big)
\end{equation}
D'un autre côté, si on calcule un peu, on trouve
\begin{align*}
\bigoplus_{i=1}^r \Ext_\LL^1(E_i,F_i) &\cong \bigoplus_{i=1}^r \Ext(E_i \otimes F_i^*,k_0)\\
&\cong \bigoplus_{i=1}^r \Ext_\LL^1(k_0 \oplus N_i,k_0)\\
&\cong \bigoplus_{i=1}^r \big(\Ext_\LL^1(k_0,k_0) \oplus \Ext_\LL^1(N_i,k_0)\big)\\
&\cong \Big(\bigoplus_{i=1}^r \Ext_\LL^1(k_0,k_0)\Big) \oplus \Big(\bigoplus_{i=1}^r \Ext_\LL^1(N_i,k_0)\Big)\\
&\cong \Big(\bigoplus_{i=1}^{r-1} \Ext_\LL^1(k_0,k_0)\Big) \oplus \Ext_\LL^1(k_0,k_0) \oplus \Big(\bigoplus_{i=1}^r \Ext_\LL^1(N_i,k_0)\Big)\\
&\cong \Big(\bigoplus_{i=1}^{r-1} \Ext_\LL^1(k_0,k_0)\Big) \oplus \Ext_\LL^1(E,F)
\end{align*}
où à la dernière ligne, l'isomorphisme (\ref{eurk}) a été utilisé. Pour terminer, le corollaire \ref{extdim1} donne que $\Ext_\LL^1(k_0,k_0) \cong (\LL / \LL')^*$ alors l'isomorphisme qui précède se réécrit
$$ \big((\LL / \LL')^*\big)^{\oplus\,r-1} \oplus \Ext_\LL^1(E,F) \cong \bigoplus_{i=1}^r \Ext_\LL^1(E_i,F_i) $$
Le cas $f_E = f_F$ est donc lui aussi traité. \qed
\end{pre}
\end{Thm}

Le précédent théorème joue un rôle clé pour bien décrire la décomposition en Ext-blocs de la catégorie des $\LL$-modules simples de dimension finie. Il réduit essentiellement la problématique de comparer les Ext-blocs de deux modules simples d'évaluation à comparer les Ext-blocs des modules sur chacune de leurs $\Gamma$-orbites communes associées.

Pour bien présenter cette description, il faut d'abord introduire quelques notions. Encore, $\LL$ est ici une algèbre de courants tordue et $\LL_{\ev}$-$\mathbf{mod}$ désignera sa catégorie de modules d'évaluation de dimension finie pour tout ce qui suit.

\begin{Def} Pour un idéal maximal $M \in \Max S$, notons 
$$ \Bloc(\LL \, | \, \Gamma.M) $$ 
l'ensemble des blocs d'extensions de la catégorie des $\LL$-modules d'évaluation de dimension finie dont tous les facteurs ont un support égal à l'orbite $\Gamma.M \subseteq \Max S$.
\end{Def}

\begin{Rem} La notation $\Bloc(\LL \, | \, \Gamma.M)$ ci-dessus ne doit pas être confondue avec $\Bloc(\g^M)$ qui, selon la définition \ref{defdebasebloc}, n'est que l'ensemble des blocs d'extensions de \linebreak $\g^M$-$\mathbf{mod}$. Ceci est une distinction très importante !
\end{Rem}

De manière analogue au cadre de la classification des classes d'isomorphismes de $\LL$-modules simples d'évaluation, nous considérerons ici la fibration
$$ \B \; = \bigsqcup_{M \in \Max S} \Bloc(\LL \, | \, \Gamma.M) \twoheadrightarrow \Max S $$

Il est naturel de considérer des sections de cette nouvelle fibration $\B$. Posons donc
$$ \secb = \{ \text{ sections de } \B \text{ }\} $$

Notons bien que chaque bloc d'extensions de la catégorie $\LL$-$\mathbf{mod}$ contient au moins une classe d'isomorphismes de $\g^M$-modules simple et qu'inversément, toute classe d'isomorphismes de $\g^M$-module simple appartient à un bloc d'extensions de la catégorie $\LL$-$\mathbf{mod}$. On peut donc voir une section de $\B$ comme une section de $\R$ suivie des applications bien définies
\begin{align*}
\bl_M : \; \Rep(&\g^M) \rightarrow \Bloc(\LL \, | \, \Gamma.M)\\
&[V] \; \longmapsto \; \llbracket V \rrbracket
\end{align*}
où $V$, un $\g^M$-module simple, est ré-interprété comme un $\LL$-module simple d'évaluation. 

Plus précisément, si $\chi$ est une section de $\B$, alors pour chaque $M \in \Max S$, il existe au moins un $\g^M$-module simple $V_M$ dans le bloc d'extensions $\chi(M) \in \Bloc(\LL \, | \, \Gamma.M)$ alors on peut écrire $\chi(M) = \llbracket V_M \rrbracket$. Enfin, on peut construire une section de $\R$, à savoir $f_\chi : M \mapsto [V_M]$, et écrire 
\begin{equation} \label{fchi}
\chi(M) = \bl_M\big(f_\chi(M)\big) = \big(\bl_M \circ f_\chi\big)(M) \qquad \text{pour chaque } M \in \Max S 
\end{equation}

\begin{Rem} D'après la définition même des applications $\bl_M$, les choix des classes d'isomorphismes $[V_M] \in \chi(M)$ pour définir $f_\chi$ ne sont pas importants dans l'écriture des relations (\ref{fchi}).
\end{Rem}

S'il est exact que les objets simples de $\LL_{\ev}$-$\mathbf{mod}$ sont classifiés par $\F^\Gamma$, les sections $\Gamma$-invariantes à support fini de $\R$, on peut envisager de considérer les sections de $\B$ qui ont des propriétés analogues. 

Soit $\chi : \Max S \rightarrow \B$, une section de $\B$. Son support est défini comme
$$ \supp \chi = \big\{M \in \Max S \; | \; \chi(M) \neq \llbracket k_0 \rrbracket \in \Bloc(\LL \, | \, \Gamma.M)\big\}$$

Maintenant, soit $\gamma \in \Gamma$ et soit $\chi \in \secb$. Choisissons $f_\chi \in \secr$ de façon à respecter les relations de la ligne (\ref{fchi}). On définit ${}^\gamma \chi \in \secb$ par la formule
\begin{equation} \label{actioncarspec}
\big({}^\gamma \chi\big)(M) = \big(\bl_M \circ \, {}^\gamma f_\chi\big)(M)
\end{equation}
où l'action de $\gamma \in \Gamma$ sur $f_\chi \in \secr$ est celle définie à la ligne (\ref{actionfgamma}). Plus explicitement, si $M \in \Max S$, on peut choisir un $\g^M$-module simple $(\rho_{f_\chi,M},V)$ dans la classe d'isomorphismes $f_\chi(M) \in \Rep(\g^M)$ de sorte que $\chi(M) = \left\llbracket(\rho_{f_\chi,M},V)\right\rrbracket$. Ensuite, (\ref{actioncarspec}) indique que
\begin{equation} \label{actioncarspecexplicite}
{}^\gamma \chi : \; M \longmapsto \left\llbracket\Big(\rho_{f,{}^{\gamma^{-1}}M} \circ \big(\gamma(M)\big)^{-1},V\Big)\right\rrbracket
\end{equation}

\begin{Rem} Pour éviter d'avoir recours à l'écriture plutôt lourde de (\ref{actioncarspecexplicite}), cette écriture sera abrégée par
$$ \big({}^\gamma \chi\big)(M) = \left\llbracket \,\chi({}^{\gamma^{-1}}M) \circ \big(\gamma(M)\big)^{-1}\right\rrbracket $$ 
\end{Rem}

\begin{Rem} La définition ${}^\gamma \chi$ à la ligne (\ref{actioncarspec}) et le fait que $\gamma \in \Gamma$ agisse sur $f \in \secr$ par $f \mapsto {}^\gamma f$, révèlent ensemble que l'application
\begin{align*}
&\Gamma \times \secb \longrightarrow \secb\\
&\;(\gamma,\chi) \mapsto {}^\gamma \chi
\end{align*}
est une action du groupe $\Gamma$ sur l'ensemble des sections $\secb$ de la fibration $\B$.
\end{Rem}

Suit maintenant la définition anlalogue, pour la fibration $\B$, des sections $\F^\Gamma$ de la fibration $\R$. Comme $\F^\Gamma$ classifie les objets simples de $\LL_{\ev}$-$\mathbf{mod}$ (voir le théorème \ref{classmodulesev}), on peut s'attendre à ce que leurs analogues pour $\B$ soient intéressant d'un point de vue théorique.

\begin{Def} Soit $\LL = (\g \otimes S)^\Gamma$ une algèbre de courants tordue. Alors un \textbf{caractère spectral} de $\LL$ est une section $\Gamma$-invariante de la fibration $\B$ dont le support est fini.

L'ensemble des caractères spectraux est noté $\,\carspec^\Gamma$.
\end{Def}

\begin{Rem} Le nom donné à ces sections a une raison d'être historique. L'idée d'introduire des fonctions de ce genre pour l'étude des blocs d'extensions est éprouvée. Tout récemment par exemple, une notion de caractère spectral joue un rôle de premier plan dans \cite{NSext}, un article de 2011 sur lequel ce travail est construit. Une autre notion de caractère spectral apparaît aussi dans \cite{CMspec}, un article de 2004. Les auteurs de ce dernier article ont été inspirés par une notion de « caractère elliptique » donnée dans un article de P.Etingof et A.Moura, paru en 2003, pour l'étude des algèbres quantiques affines. L'idée originale des caratères spectraux n'est donc pas nouvelle, mais chose sûre, elle est toujours pertinente au $\text{XXI}^{\text{e}}$ siècle ! 
\end{Rem}

Sachant que les classes d'isomorphismes d'objets simples de $\LL_{\ev}$-$\mathbf{mod}$ correspondent aux sections $\F^\Gamma \subseteq \secr$ de la fibration $\R$ (voir le théorème \ref{classmodulesev}), on peut associer à chacune de ces classes d'isomorphismes une section de la fibration $\B$ via les applications $\bl_M$ comme dans (\ref{actioncarspec}). Plus concrètement, ceci revient à faire l'association
\begin{align*} 
&\left\{\begin{array}{c} \text{Classes d'isomorphismes de}\\
\text{d'objets simples de $\LL_{\ev}$-$\mathbf{mod}$} \end{array}\right\} \; \leftrightarrow \; \F^\Gamma \qquad \qquad \longrightarrow \qquad \quad \; \secb\\
&\qquad \qquad \qquad \quad [E] \qquad \qquad \longleftrightarrow \qquad \;\; f_E \qquad \qquad \longmapsto \qquad \quad \bl \circ f_E 
\end{align*}
où $\big(\bl \circ f_E\big) : M \mapsto \bl_M\big(f_E(M)\big) \in \Bloc(\LL \, | \, \Gamma.M)$ pour tout $M \in \Max S$.

\begin{Rem} Le fait que $f_E \in \F^\Gamma$ donne immédiatement que $\bl \circ f_E \in \carspec^\Gamma$ pour toute telle classe d'isomorphismes $[E]$. En effet, $\bl \circ f_E \in \secb$ a un support fini parce que c'est le cas pour $f_E \in \F^\Gamma$ et le fait que $\bl \circ f_E$ soit une section $\Gamma$-invariante découle de la définition de l'action de $\Gamma$ à la ligne (\ref{actioncarspec}) et que $f_E$ le soit par hypothèse.
\end{Rem}

\begin{Def} Soit $E$ un $\LL$-module simple d'évaluation de dimension finie. Alors \textbf{le caractère spectral associé à E} est défini comme 
$$ \chi_E \, = \, \bl \circ f_E \in \carspec^\Gamma $$
De façon équivalente, si $f \in \F^\Gamma$, je noterai $\chi_f$ le caractère spectral $\, \bl \circ f \in \carspec^\Gamma$.

Explicitement, si $E = \bigotimes_{i=1}^r E_i$ où pour chaque indice $i \in \{1,...\,,r\}$, $(\rho_i,E_i)$ est un \linebreak $\g^{M_i}$-module simple et où $\{M_i\}_{i=1}^r \subseteq \Max S$ est un ensemble d'idéaux maximaux issus de $r$ $\Gamma$-orbites distinctes, alors
\begin{align*}
\chi_E : \; &\Max S \longrightarrow \bigsqcup_{M \in \Max S} \Bloc(\LL \, | \, \Gamma.M) \quad = \quad \B\\
 &\quad M \longmapsto \left\{
\begin{array}{cl} \left\llbracket \Big(\rho_i \circ \big(\gamma(M)\big)^{-1}, E_i \Big) \right\rrbracket & \text{si }M = {}^\gamma M_i\\
\text{$\llbracket \, k_0 \, \rrbracket$} & \text{sinon}\\
\end{array} \right.
\end{align*}
\end{Def}

Les résultats qui suivent complètent graduellement la classification des blocs d'extensions pour les algèbres de courants tordues $\LL$ qui vérifient une certaine propiété technique. Les caractères spectraux jouent un rôle de premier plan dans le ces efforts de classification. Le prochain résultat peut en témoigner.

\begin{Prop} \label{extblocev} Soient $E$ et $F$ deux $\LL$-modules simples d'évaluation de dimension finie. 

Alors $\llbracket E \rrbracket = \llbracket F \rrbracket$ dans $\LL_{\ev}$-$\mathbf{mod}$
$$ \qquad \Longleftrightarrow \quad \chi_E = \chi_F \in \carspec^\Gamma $$
\begin{pre} Supposons d'abord que $E$ et $F$ soient dans un même Ext-bloc. Si $E \cong F$, la classification des classes d'isomorphismes de $\LL$-modules simples d'évaluation donne que \linebreak $f_E = f_F \in \F^\Gamma$. Par la définition même des caractères spectraux, cela implique que 
$$ \chi_E = \chi_F \in \carspec^\Gamma $$
Nous pouvons donc supposer que $E \ncong F$ sont dans le même Ext-bloc et qu'il existe une suite finie de $\LL$-modules simples d'évaluation de dimension finie $\{T^j\}_{j=1}^n$ avec $T^1 = E$, $T^n = F$, puis, pour chaque $m \in \{1,...\,,n-1\}$, on a $T^m \ncong T^{m+1}$ et
$$ \Ext^1_\LL(T^m,T^{m+1}) \neq \{0\} \qquad \text{ou} \qquad \Ext^1_\LL(T^{m+1},T^m) \neq \{0\} $$
Pour conclure, il suffit donc de montrer que pour chacun de ces indices $m$, on ait $\chi_{T^m} = \chi_{T^{m+1}}$. Fixons donc $m$ et, sans perte de généralité, supposons que $\Ext^1_\LL(T^m,T^{m+1})$. Le théorème \ref{exteval} permet alors de déduire que $f_{T^m}$ et $f_{T^{m+1}}$ diffèrent exactement sur une $\Gamma$-orbite de $\Max S$ (seule l'éventualité (2) est envisageable). Disons qu'il s'agit de l'orbite de $M \in \Max S$ et notons $T^m_M$ et $T^{m+1}_M$ les $\g^M$-modules simples non-isomorphes qui apparaissent dans les écritures de $T^m$ et $T^{m+1}$, respectivement. 

La proposition \ref{exteval} (2) donne alors que 
$$ \Ext^1_\LL(T^m,T^{m+1}) = \Ext^1_\LL(T^m_M,T^{m+1}_M) \neq \{0\} $$
C'est donc que $\llbracket T^m_M \rrbracket = \llbracket T^{m+1}_M \rrbracket \in \Bloc(\LL \, | \, \Gamma.M)$. Ainsi, $\chi_{T^m}(M) = \chi_{T^{m+1}}(M)$, d'où on peut déduire
$$ \chi_{T^m} = \chi_{T^{m+1}} \in \carspec^\Gamma $$

Supposons maintenant que $\chi_E = \chi_F \in \carspec^\Gamma$. Alors, il faut montrer que $\llbracket E \rrbracket = \llbracket F \rrbracket$ dans $\LL_{\ev}$-$\mathbf{mod}$. Fixons $\{M_i\}_{i=1}^r \subseteq \Max S$ un ensemble complet de représentants des $\Gamma$-orbites de $\supp E \cup \supp F$. Supposons donc que $E = \bigotimes_{i=1}^r E_i$ et  $F = \bigotimes_{i=1}^r F_i$ où les $E_i$ et $F_i$ sont des $\g^{M_i}$-modules simples pour chacun des $i$.

Ce qui suit est une induction sur le nombre de $\Gamma$-orbites de $\Max S$ telles que les fonctions de $\F^\Gamma$ associées diffèrent. Soient $V$ et $W$ des $\LL$-modules d'évaluation tels que ce nombre vaut $0$, c'est donc que $f_V = f_W \in \F^\Gamma$, donc $V \cong W$ et ainsi, $\llbracket V \rrbracket = \llbracket W \rrbracket$ dans $\LL_{\ev}$-$\mathbf{mod}$.

Fixons maintenant le nombre de $\Gamma$-orbites de $\Max S$ où $f_E \neq f_F$ à $n > 1$ et supposons que pour toute paire de $\LL$-modules simples d'évaluation de dimension finie $V$ et $W$ qui vérifie que $f_V \neq f_W$ sur $n-1$ $\Gamma$-orbites ou moins, on ait $\llbracket V \rrbracket = \llbracket W \rrbracket$ dans $\LL_{\ev}$-$\mathbf{mod}$. Il faut montrer que ce résultat tienne pour la paire $E$ et $F$ également (c'est-à-dire qu'il tienne pour $n$ $\Gamma$-orbites où les sections diffèrent).

Fixons $M_k \in \Max S$ tel que $f_E(M_k) \neq f_F(M_k)$. Alors, comme $\chi_E(M_k) = \chi_F(M_k)$ par hypothèse, on a que $\llbracket E_k \rrbracket = \llbracket F_k \rrbracket \in \Bloc(\g^{M_k})$. Il existe donc une suite finie de $\g^{M_k}$-modules simples d'évaluation de dimension finie $\{T^j\}_{j=1}^n$ avec $T^1 = E_k$, $T^N = F_k$, puis, pour chaque $m \in \{1,...\,,n-1\}$, on a $T^m \ncong T^{m+1}$ et
$$ \Ext^1_\LL(T^m,T^{m+1}) \neq \{0\} \qquad \text{ou} \qquad \Ext^1_\LL(T^{m+1},T^m) \neq \{0\} $$
Posons 
$$ \tilde{E} = \bigotimes_{\substack{i=1 \\ i \neq k}}^r E_i \qquad \text{et} \qquad \tilde{F} = \bigotimes_{\substack{i=1 \\ i \neq k}}^r F_i $$
et puis $\tilde{T}^m = T^m \otimes \tilde{E}$ pour chacun des $m$. Par construction, $f_{\tilde{T}^m} \neq f_{\tilde{T}^{m+1}}$ que sur une $\Gamma$-orbite de $\Max S$. Par le théorème \ref{exteval}, on a donc
$$ \Ext^1_\LL(\tilde{T}^m,\tilde{T}^{m+1}) = \Ext^1_\LL(T^m,T^{m+1}) \qquad \text{et} \qquad \Ext^1_\LL(\tilde{T}^{m+1},\tilde{T}^m) = \Ext^1_\LL(T^{m+1},T^m) $$
Pour chacun des $m$, l'un de ces deux ensemble n'est pas $\{0\}$ par hypothèse et donc, \linebreak $\llbracket \tilde{T}^1 \rrbracket = \llbracket \tilde{T}^N \rrbracket$ dans $\LL_{\ev}$-$\mathbf{mod}$. Par ailleurs, $\tilde{T}^1 = E$ et $\tilde{T}^N = T^N \otimes \tilde{E}$.

Ensuite, c'est que $\tilde{T}^N = F_k \otimes \tilde{E}$ et donc $f_{\tilde{T}^N} \neq f_F$ sur précisément $n-1$ des $\Gamma$-orbites de $\Max S$. De l'hypothèse d'induction, il suit alors que $\llbracket \tilde{T}^N \rrbracket =  \llbracket F \rrbracket$ dans $\LL_{\ev}$-$\mathbf{mod}$. C'est ainsi que $\llbracket E = \tilde{T}^1 \rrbracket = \llbracket \tilde{T}^N\rrbracket = \llbracket F \rrbracket$ dans $\LL_{\ev}$-$\mathbf{mod}$.
\qed
\end{pre}
\end{Prop}

On est maintenant en mesure de donner la décomposition de la catégorie $\LL_{\ev}$-$\mathbf{mod}$ en blocs. Il faut cependant préciser une définition avant de ce faire.

\begin{Def} Soit $\chi \in \carspec^\Gamma$. Alors, on définit $\Cat^\chi$ comme étant la sous-catégorie (pleine) de $\LL_{\ev}$-$\mathbf{mod}$ des modules d'évaluation dont chaque facteur irréductible a $\chi$ pour caractère spectral.
\end{Def}

\begin{Cor} \label{blocev} La décomposition de $\LL_{\ev}$-$\mathbf{mod}$ en blocs est
$$ \text{$\LL_{\ev}$-$\mathbf{mod}$} = \bigoplus_{\chi \in \carspec^\Gamma} \Cat^\chi $$
En particulier, les blocs et blocs d'extensions de $\LL_{\ev}$-$\mathbf{mod}$ sont en bijection avec les caractères spectraux.
\end{Cor}

Maintenant, il faudrait pouvoir décrire les blocs de la catégorie $\LL$-$\mathbf{mod}$. Il faut se souvenir que d'après la classification des objets simples de $\LL$-$\mathbf{mod}$ et la proposition \ref{classmodules}, tout tel objet est isomorphe à un unique $k_\lambda \otimes E_f$ correspondant à un couple $(\lambda,f) \in \LL^* \times \F^\Gamma$ satisfaisant aux propriétés 
\begin{align*}
&\bullet \quad \lambda|_{[\LL,\LL]} = 0.\\
&\bullet \quad \LL/\ker \rho_f \text{ est semi-simple}.\\
&\bullet \quad \ker(\lambda+\rho_f) = \ker \lambda \cap \ker \rho_f.
\end{align*}
Appelons momentanément de tels couples, des couples « admissibles ». Il suit que l'ensemble des Ext-blocs de $\LL$-$\mathbf{mod}$ est égal à l'ensemble 
$$ \big\{\llbracket k_\lambda \otimes E_f \rrbracket \; | \; (\lambda,f) \text{ est un couple « admissibles »}\big\} $$ 
Ceci n'est, par contre, pas précis du tout. Il reste à déterminer selon quels critères les couples « admissibles » donnent lieu à un même bloc d'extensions dans $\LL$-$\mathbf{mod}$, c'est-à-dire qu'il reste à déterminer quelles classes d'isomorphismes de $\LL$-modules simples sont dans le même bloc d'extensions de $\LL$-$\mathbf{mod}$.

La classification des blocs d'extensions de $\LL_{\ev}$-$\mathbf{mod}$ est donnée par la proposition \ref{extblocev} via les caractères spectraux $\carspec^\Gamma$. La définition même du caractère spectral associé à un objet simple de $\LL_{\ev}$-$\mathbf{mod}$ laisse penser que certains éléments de l'ensemble $\LL^* \times \carspec^\Gamma$ peuvent prétendre à classifier les blocs d'extensions de $\LL$-$\mathbf{mod}$. L'idée est plutôt naturelle et consisterait à identifier les blocs d'extensions des couples $(\lambda,f)$ « admissibles » à l'élément $(\lambda,\chi_f) \in \LL^* \times \carspec^\Gamma$. Idéalement, il s'agirait donc de réaliser une identification bijective qui ressemblerait à 
$$ \llbracket k_\lambda \otimes E_f \rrbracket \qquad \longleftrightarrow \qquad (\lambda,\chi_f) $$

Une telle identification va pouvoir fonctionner à quelques modifications près. Quoi qu'il en soit, pour vraiment procéder à la classifications des blocs d'extensions de $\LL$-$\mathbf{mod}$ par cette méthode, il faut exiger que les algèbres de courants tordues $\LL$ vérifient une condition technique. 

\begin{Def} \label{extlocal} Une algèbre de courants tordue $\LL$ est dite \textbf{locale vis-à-vis des extensions} si la propriété suivante est toujours respectée :
$$ \left.\begin{array}{c} E \text{ un objet simple de $\LL_{\ev}$-$\mathbf{mod}$} \\ \lambda \in \LL^* \text{ tel que $k_\lambda \notin$ $\LL_{\ev}$-$\mathbf{mod}$}\end{array}\right] \quad \Longrightarrow \quad \Ext^1_\LL(k_\lambda,E) = \{0\} $$
\end{Def}

\begin{Rem} Les auteurs de \cite{NSext}, E.Neher et A.Savage notent dans leur article sur les algèbres d'applications équivariantes que dans le cas de tous les exemples de base, les algèbres vérifient cette propriété. Ils reconnaissent même ne pas connaître d'algèbre d'applications équivariantes qui ne soit pas locale vis-à-vis des extensions. Voir la remarque 5.14 de \cite{NSext} pour plus de détails.
\end{Rem}

\begin{Rem} \label{parfaiteextlocal} Il est également bon de prendre note qu'une algèbre de courants tordue pour laquelle tous les modules simples sont d'évaluation satisfait automatiquement la propriété de localité vis-à-vis des extensions. Par exemple, c'est le cas des formes tordues ; il en sera question plus en détail dans le chapitre 3.

C'est le cas par exemple si $\LL$ est parfaite (i.e. $[\LL,\LL] = \LL$) puisqu'alors, $\LL/\LL' \cong \{0\}$ et la seule représentation de $\LL$ de dimension 1 est $k_0$ ; un module d'évaluation.
\end{Rem}

\begin{Rem} Il serait très intéressant si l'on pourrait se débarasser du recours à la propriété de localité vis-à-vis des extensions. À l'image des auteurs de \cite{NSext}, je ne connais pas d'algèbre de courants tordue qui ne vérifie pas cette propriété. Il faut dire qu'il est plutôt difficile d'expliciter des algèbres de courants tordues avec tous leurs paramètres et ma connaissance d'exemples demeure limitée. 

Par ailleurs, aucune raison que je connais ne tend à favoriser l'hypothèse que toutes les algèbres de courants tordues serait locale vis-à-vis des extensions. De la même façon, aucune raison ne pousse à croire que la localité vis-à-vis des extension est une hypothèse superflue ou non pour la classification des blocs d'extensions.
\end{Rem}

Poursuivons néanmoins la démarche de classification. Posons 
$$ (\LL/\LL')_{\ev}^* = \{\lambda \in (\LL/\LL')^* \; | \; k_\lambda \text{ est un $\LL$-module d'évaluation}\} $$
On a d'abord évidemment que $0 \in (\LL/\LL')_{\ev}^*$, puis les isomorphismes \ref{produitdimension1} et \ref{dualdimension1} montrent que cet ensemble est un sous-espace vectoriel du $k$-espace $(\LL/\LL')^*$.

Pour la suite, fixons un choix d'un espace vectoriel complémentaire au sous-espace-vectoriel $(\LL/\LL')_{\ev}^*$ de $(\LL/\LL')^*$. Notons ce complémentaire $(\LL/\LL')_{\nev}^*$ de sorte qu'on ait la décomposition en somme directe 
\begin{equation} \label{decompositionnev}
(\LL/\LL')^* = (\LL/\LL')_{\ev}^* \oplus (\LL/\LL')_{\nev}^*
\end{equation}

D'après la classification des objets simples de $\LL$-$\mathbf{mod}$ (la proposition \ref{classmodules}), les $\LL$-modules simples de dimension finie ont une écriture unique comme $k_\lambda \otimes E_f$ où $\lambda \in \LL^*$ et $f \in \F^\Gamma$ satisfont certaines conditions. Avec une décomposition fixée comme (\ref{decompositionnev}), on peut alors décomposer $\lambda$ de façon unique comme 
$$ \lambda = \lambda_{\ev} + \lambda_{\nev} $$ 
De façon équivalente, la décomposition (\ref{decompositionnev}) choisie mène à une façon unique de réécrire le module simple comme 
\begin{align*}
k_\lambda \otimes E_f &\cong k_{\lambda_{\ev}} \otimes k_{\lambda_{\nev}} \otimes E_f\\
&\cong k_{\lambda_{\nev}} \otimes (k_{\lambda_{\ev}} \otimes E_f)
\end{align*}
où $k_{\lambda_{\ev}} \otimes E_f$ est un $\LL$-module (simple) d'évaluation de dimension finie. En conséquence, $k_{\lambda_{\ev}} \otimes E_f$ a un caractère spectral.

\begin{Def} \label{partieeval} Supposons qu'une décomposition (\ref{decompositionnev}) a été fixée. La \textbf{partie d'évaluation} d'une classe d'isomorphismes de $\LL$-module simple $[k_\lambda \otimes E_f]$ identifiable à $(\lambda,f) \in \LL^* \times \F^\Gamma$ vérifiant les conditions de la proposition \ref{classmodules} est la classe d'isomorphismes du $\LL$-module simple d'évaluation $k_{\lambda_{\ev}} \otimes E_f$. Autrement dit, la partie d'évaluation de $(\lambda,f)$ est la classe d'isomorphismes $[k_{\lambda_{\ev}} \otimes E_f]$. 

Le caractère spectral de la partie d'évaluation de la classe d'isomorphismes de $\LL$-module simple classifiée par $(\lambda,f) \in \LL^* \times \F^\Gamma$ sera noté $\chi_f^\lambda \in \carspec^\Gamma$. 
\end{Def}

Voici enfin le résultat qui classifie les blocs d'extensions des algèbres de courants tordues qui sont locales vis-à-vis des extensions : 

\begin{Thm} Soit $\LL$ une algébre de courants tordue qui soit locale vis-à-vis des extensions et soient $k_\lambda \otimes E_f$ et $k_\mu \otimes E_g$ des $\LL$-modules simples correspondant aux couples $(\lambda,f)$ et $(\mu,g) \in \LL^* \times \F^\Gamma$ vérifiant les trois conditions du théorème \ref{classmodules}.

Fixons une décomposition 
$$ (\LL/\LL')^* = (\LL/\LL')_{\ev}^* \oplus (\LL/\LL')_{\nev}^* $$
Alors, $\llbracket k_\lambda \otimes E_f \rrbracket = \llbracket k_\mu \otimes E_g \rrbracket$ dans $\LL$-$\mathbf{mod}$
$$ \qquad \Longleftrightarrow \quad (\,\lambda_{\nev},\,\chi_f^\lambda\,) = (\,\mu_{\nev},\,\chi_g^\mu\,) \;\in\; (\LL/\LL')_{\nev}^* \times \; \carspec^\Gamma $$
\begin{pre} Supposons d'abord que $(\,\lambda_{\nev},\,\chi_f^\lambda\,) = (\,\mu_{\nev},\,\chi_g^\mu\,) \in (\LL/\LL')_{\nev}^* \times \; \carspec^\Gamma$. 

Selon la proposition \ref{extblocev}, l'hypothèse donne que $\llbracket k_{\lambda_{\ev}} \otimes E_f \rrbracket = \llbracket k_{\mu_{\ev}} \otimes E_g \rrbracket$ dans $\LL_{\ev}$-$\mathbf{mod}$. Il existe donc une suite finie de $\LL$-modules d'évaluation de dimension finie $\{T^j\}_{j=1}^n$ avec \linebreak $T^1 = k_{\lambda_{\ev}} \otimes E_f$, $T^n = k_{\mu_{\ev}} \otimes E_g$, puis, pour chaque $m \in \{1,...\,,n-1\}$, on a $T^m \ncong T^{m+1}$ et
$$ \Ext^1_\LL(T^m,T^{m+1}) \neq \{0\} \qquad \text{ou} \qquad \Ext^1_\LL(T^{m+1},T^m) \neq \{0\} $$
Puisque $\lambda_{\nev} = \mu_{\nev}$, on a les égalités
\begin{align*} 
\Ext^1_\LL(T^m,T^{m+1}) &= \Ext^1_\LL(T^m \otimes k_0,T^{m+1}) = \Ext^1_\LL(T^m \otimes k_{\lambda_{\nev}},T^{m+1} \otimes k_{\mu_{\nev}})\\
\Ext^1_\LL(T^{m+1},T^m) &= \Ext^1_\LL(T^{m+1} \otimes k_0,T^m) = \Ext^1_\LL(T^{m+1} \otimes k_{\lambda_{\nev}},T^m \otimes k_{\mu_{\nev}})
\end{align*}
Ainsi, les $\LL$-modules de dimension finie $\{T^j \otimes k_{\lambda_{\nev}}\}_{j=1}^n$ montrent que $\llbracket k_\lambda \otimes E_f \rrbracket = \llbracket k_\mu \otimes E_g \rrbracket$ dans $\LL$-$\mathbf{mod}$.

Supposons ensuite que $\llbracket k_\lambda \otimes E_f \rrbracket = \llbracket k_\mu \otimes E_g \rrbracket$ dans $\LL$-$\mathbf{mod}$. Alors il existe une suite $\LL$-modules de dimension finie $\{\lambda_j \otimes T^j\}_{j=1}^n$ avec 
\begin{align*}
&\bullet \quad \lambda_j \in (\LL/\LL')_{\nev}^* \text{ pour tout } m \in \{1,...\,,n-1\}.\\
&\bullet \quad T^j \text{ est un module d'évaluation pour tout } m \in \{1,...\,,n-1\}.\\
&\bullet \quad T^j \ncong T^{j+1} \text{ pour tout } j.\\
&\bullet \quad \lambda_1 \otimes T^1 \cong k_{\lambda} \otimes E_f \text{ et } \lambda_1 = \lambda_{\nev}.\\
&\bullet \quad \lambda_n \otimes T^n \cong k_{\mu} \otimes E_g \text{ et } \lambda_n = \mu_{\nev}.\\
&\bullet \quad \text{Pour chaque } m \in \{1,...\,,n-1\} \text{, soit } \Ext^1_\LL(k_{\lambda_m} \otimes T^m,k_{\lambda_{m+1}} \otimes T^{m+1}) \neq \{0\}\\
& \quad \,\,\,\, \text{ ou soit } \Ext^1_\LL(k_{\lambda_{m+1}} \otimes T^{m+1}, k_{\lambda_m} \otimes T^m) \neq \{0\}.
\end{align*}
Par ailleurs, on a les isomorphismes
\begin{align*}
\Ext^1_\LL(k_{\lambda_m} \otimes T^m,k_{\lambda_{m+1}} \otimes T^{m+1}) &\cong \Ext^1_\LL(k_{\lambda_m - \lambda_{m+1}}, (T^m)^* \otimes T^{m+1})\\
\Ext^1_\LL(k_{\lambda_{m+1}} \otimes T^{m+1},k_{\lambda_m} \otimes T^m) &\cong \Ext^1_\LL(k_{\lambda_{m+1} - \lambda_m}, (T^{m+1})^* \otimes T^m)
\end{align*}
Comme $\LL$ est locale vis-à-vis des extensions (voir la définition \ref{extlocal}), alors les dernières équations et les hypothèses de la petite liste ci-haut entrainent que $\lambda_m = \lambda_{m+1}$ pour tout $m \in \{1,...\,,n-1\}$ ; c'est donc que $\lambda_{\nev} = \mu_{\nev}$. Alors, les deux égalités précédentes entrainent directement que $\llbracket T_1 \rrbracket = \llbracket T^n \rrbracket$ dans $\LL_{\ev}$-$\mathbf{mod}$. Puisque $T^1 = k_{\lambda_{\ev}} \otimes E_f$ et $T^n = k_{\lambda_{\ev}} \otimes E_g$, cela est équivalent à dire que $ \llbracket k_{\lambda_{\ev}} \otimes E_f \rrbracket = \llbracket \mu_{\nev} \otimes E_g \rrbracket$ dans $\LL_{\ev}$-$\mathbf{mod}$.

En somme, l'information conclue par les explications du paragraphe précédent et la proposition \ref{extblocev} est équivalente à la conclusion que 
$$ (\,\lambda_{\nev},\,\chi_f^\lambda\,) = (\,\mu_{\nev},\,\chi_g^\mu\,) \;\in\; (\LL/\LL')_{\nev}^* \times \; \carspec^\Gamma $$
\vspace{-2em}
\qed
\end{pre}
\end{Thm}

Il est maintenant possible de donner la décomposition de la catégorie $\LL$-$\mathbf{mod}$ en blocs.

\begin{Def} Soit $(\alpha,\chi) \in (\LL/\LL')_{\nev}^* \times \; \carspec^\Gamma$. Alors on définit $\Cat^{(\alpha,\chi)}$ comme étant la sous-catégorie (pleine) de $\LL$-$\mathbf{mod}$ des modules dont chaque facteur irréductible a une partie d'évaluation dont le caractère spectral est $\chi$ et dont la partie de dimension 1, qui n'est pas d'évaluation, est $\alpha$.
\end{Def}

\begin{Cor} \label{bloc} Si $\LL$ est une algèbre de courants tordues locale vis-à-vis des extensions, la décomposition de $\LL$-$\mathbf{mod}$ en blocs est
$$ \text{$\LL$-$\mathbf{mod}$} = \bigoplus_{(\alpha,\chi) \,\in\, (\LL/\LL')_{\nev}^*\times \, \mathsmaller{\carspec}^\Gamma} \Cat^{(\alpha,\chi)} $$
En particulier, les blocs et les blocs d'extensions de $\LL$-$\mathbf{mod}$ sont en bijection avec les éléments de $(\LL/\LL')_{\nev}^*\times \, \mathsmaller{\carspec}^\Gamma$.
\end{Cor}

\begin{Rem} Ceci complète donc la classification des blocs d'extensions des algèbres de courants tordues locales vis-à-vis des extensions. Bien entendu, cette classification n'est pas idéale dans ce qu'elle n'est pas exprimée en les même termes que la classification des objets simples de $\LL$-$\mathbf{mod}$ (voir la proposition \ref{classmodules}). Aussi, le fait qu'il faille fixer un choix de décomposition de $(\LL/\LL')^*$ comme \ref{decompositionnev} est, certes, embêtant.
\end{Rem}

\subsection{Sommaire, précisions et implications}

Cette sous-section a pour but de relever et préciser les résultats clés de la classification des blocs d'extensions d'une algèbre de courants tordue, dans le but de les insérer dans un contexte plus pratique. Pour ce faire, il faut d'abord préciser, autant que possible, la valeur de $\Ext_\LL^1(V,W)$ où $\LL$ est une algèbre de courants tordue et où $V$ et $W$ sont des objets simples de $\LL$-$\mathbf{mod}$ arbitraires. 

\begin{Rem} Pour avoir accès à la classification des blocs d'extensions de la précédente sous-section, on devra supposer que  $\LL$ soit locale vis-à-vis des extensions (voir la définition \ref{extlocal}). Par ailleurs, l'important théorème \ref{exteval} et quelques résultats qui suivent sont indépendants de cette hypothèse.
\end{Rem}

Soit donc $\LL$ une algèbre de courants tordue locale vis-à-vis des extensions. Commençons par préciser ce que vaut $\Ext_\LL^1(V,W)$ où $V$ et $W$ sont des objets simples de $\LL$-$\mathbf{mod}$ non-isomorphes. Par la classification des objets simples de $\LL$-$\mathbf{mod}$ (la proposition \ref{classmodules}), on peut, sans perte de généralité, supposer que
\begin{align*}
V = k_\lambda \otimes E_f && W = k_\mu \otimes E_g
\end{align*}
où $(\lambda,f) \neq (\mu,g) \in \LL^* \times \F^\Gamma$ vérifient tous deux les trois conditions de la classification. Fixons aussi une décomposition 
$$ (\LL/\LL')^* = (\LL/\LL')_{\ev}^* \oplus (\LL/\LL')_{\nev}^* $$
Selon cette décomposition, fixons les décompositions correspondantes
\begin{align*}
V \; \cong \; k_{\lambda_{\nev}} \otimes (k_{\lambda_{\ev}} \otimes E_f) && W \; \cong \; k_{\mu_{\nev}} \otimes (k_{\mu_{\ev}} \otimes E_g)
\end{align*}

D'après la décomposition de la catégorie $\LL$-$\mathbf{mod}$ en blocs donnée au corollaire $\ref{bloc}$, si $(\lambda_{\nev},\,\chi_f^\lambda\,) \neq (\mu_{\nev},\,\chi_g^\mu\,)$, on a directement
$$ \Ext_\LL^1(V,W) \cong \Ext_\LL^1(k_\lambda \otimes E_f,k_\mu \otimes E_g) \cong \{0\} $$

Si par contre $(\lambda_{\nev},\,\chi_f^\lambda\,) = (\mu_{\nev},\,\chi_g^\mu\,)$, nos deux modules simples sont dans un même bloc d'extensions. Dans ce cas, par les isomorphismes (\ref{switchext}), on a
$$ \Ext_\LL^1(V,W) = \Ext_\LL^1(k_\lambda \otimes E_f,k_\mu \otimes E_g) \cong \Ext_\LL^1(k_{\lambda_{\ev}} \otimes E_f,k_{\mu_{\ev}} \otimes E_g) $$

Ainsi, il est alors suffisant de considérer les cas où on a affaire à deux objets simples de \linebreak $\LL_{\ev}$-$\mathbf{mod}$. En appliquant le théorème \ref{exteval}, on obtient que si le support des sections de $\F^\Gamma$ de $k_{\lambda_{\ev}} \otimes E_f$ et de $k_{\mu_{\ev}} \otimes E_g$ diffèrent de deux $\Gamma$-orbites de $\Max S$ ou plus, alors
$$ \Ext_\LL^1(V,W) \cong \{0\} $$
Si ce n'est pas le cas, comme on a supposé que $(\lambda,f) \neq (\mu,g) \in \LL^* \times \F^\Gamma$ et que $\lambda_{\nev} = \mu_{\nev}$, on peut appliquer la partie (2) du théorème \ref{exteval}. On obtient alors qu'il existe une seule $\Gamma$-orbite de $\Max S$ où $k_{\lambda_{\ev}} \otimes E_f$ et $k_{\mu_{\ev}} \otimes E_g$ ont des composantes non-isomorphes.

Il sera donc suffisant de décrire les classes d'équivalences d'extensions de $\LL$-modules entre deux modules simples d'évaluation non-isomorphes qui sont supportés sur une seule et même $\Gamma$-orbite de $\Max S$. La proposition suivante permet de clore la présente investigation. Notons que cette proposition plus terre-à-terre ne requiert pas l'hypothèse voulant que $\LL$ soit locale vis-à-vis des extensions.

\begin{Prop} \label{madeuxiemepropyo} Supposons que $k_{\lambda_{\ev}} \otimes E_f$ et $k_{\mu_{\ev}} \otimes E_g$ soient deux objets simples de \linebreak $\LL_{\ev}$-$\mathbf{mod}$ supportés sur la seule $\Gamma$-orbite de $M \in \Max S$ et qu'ils soient non-isomorphes. 

Sans perte de généralité, on peut supposer que ces présentations sont telles que les couples $(\lambda_{\ev},f)$ et $(\mu_{\ev},g)$  de $\LL^* \times \F^\Gamma$ vérifient les trois conditions de la proposition \ref{classmodules}. Notons qu'on a aussi supposé que $(\lambda_{\ev},f) \neq (\mu_{\ev},g)$. En particulier, ce sont des modules pour l'algèbre de Lie réductive de dimension finie $\g^M$. 

Posons $U = k_{\mu_{\ev} - \lambda_{\ev}} \otimes (E_f^* \otimes E_g)$ ainsi que
\begin{align*}
K_M &= \ker (\ev_M) \\
I &= \{s \in S \, | \, s(^\gamma M) = 0 \text{ pour tout } \gamma \in \Gamma\} \\
I^\Gamma &= I \cap S^\Gamma \\
N_M &= \{n \in K_M \, | \, I^\Gamma.\,n \subseteq K_M'\} 
\end{align*}
Alors, il y a un isomorphisme naturel
$$ \Ext_\LL^1(k_{\lambda_{\ev}} \otimes E_f,k_{\mu_{\ev}} \otimes E_g) \cong \Hom_{\g^M}(N_M/K_M',U) $$
\begin{pre} L'algèbre de Lie $\g^M$ est réductive (et de dimension finie), on donc peut écrire
$$ \g^M = \Z \oplus \SSS $$
où $\Z = Z(\g^M)$ et $\SSS = [\g^M,\g^M]$ est semi-simple.

Tel que mentionnné dans l'énoncé ci-haut, on adopte la notation
\begin{equation} \label{ecriturespdgyo}
U = k_{\mu_{\ev} - \lambda_{\ev}} \otimes (E_f^* \otimes E_g)
\end{equation}
Puisque $U \cong (k_{\lambda_{\ev}} \otimes E_f)^* \otimes (k_{\mu_{\ev}} \otimes E_g)$ et que $(\lambda_{\ev},f) \neq (\mu_{\ev},g) \in \LL^* \times \F^\Gamma$, le corollaire \ref{dimhom} donne que $U^{\g^M} = \{0\}$. Cette propriété se révèlera très importante un peu plus loin.

Notons $\rho_f$, $\tilde{\rho}_f$ et $\rho_g$ les actions de $\LL$ sur $E_f$, $E_f^*$ et $E_g$, respectivement. Puisque les conditions de la proposition \ref{classmodules} sont respectées pour les deux couples $(\lambda_{\ev},f)$ et $(\mu_{\ev},g)$, on sait que $\LL/\ker \rho_f$ et $\LL/\ker \rho_g$ sont des algèbres de Lie semi-simples. Il suit que $\LL/\ker \tilde{\rho}_f$ est aussi une algèbre de Lie semi-simple. 

Par ailleurs, $K_M \subseteq \ker \tilde{\rho}_f$ alors il y a un homomorphisme surjectif d'algèbres de Lie 
\begin{align*}
\varphi : \; \g^M \cong \;&\; \LL/K_M \twoheadrightarrow \LL / \ker \tilde{\rho}_f\\
&\ell + K_M  \mapsto \ell + \ker \tilde{\rho}_f
\end{align*}
En particulier, l'action d'un élément $x \in \g^M$ sur $E_f^*$ correspond exactement à l'action de $\varphi(x) \in \LL / \ker \tilde{\rho}_f$ sur $E_f^*$. Maintenant, comme $\Z = Z(\g^M)$ et que $\LL / \ker \tilde{\rho}_f$ est une algèbre de Lie semi-simple, on a 
$$ \varphi(\Z) \subseteq Z\big(\LL / \ker \tilde{\rho}_f\big) = \{0\} $$
C'est donc que $\Z$ agit comme 0 sur $E_f^*$. De la même manière, $\Z$ agit comme 0 sur $E_g$ et ainsi, $\Z$ agit comme 0 sur le produit tensoriel $E_f^* \otimes E_g$.

Ensuite, comme $\SSS = [\g^M,\g^M]$, on a que $\SSS$ agit comme 0 sur $k_{\mu_{\ev} - \lambda_{\ev}}$. On obtient finalement que l'écriture (\ref{ecriturespdgyo}) est telle que si $a \in k_{\mu_{\ev} - \lambda_{\ev}}$, $v \in E_f^* \otimes E_g$, alors pour $z + s \in \Z \oplus \SSS = \g^M$, on a
\begin{equation} \label{actionmaprop}
(z+s).(a \otimes v) = (z.a) \otimes v + a \otimes (s.v)
\end{equation}

Pour commencer véritablement la démonstration de la proposition, rappelons que \linebreak $\Ext_\LL^1(E,F) \cong \cohomo{\LL}{U}$ et qu'une des hypothèses est que $(\lambda_{\ev},f) \neq (\mu_{\ev},g) \in \LL^* \times \F^\Gamma$. 

Supposons d'abord que $\lambda_{\ev} \neq \mu_{\ev}$. Dans ce premier cas, $k_{\mu_{\ev} - \lambda_{\ev}} \ncong k_0$ et il existe donc un $z \in \Z$ tel que $\big(\mu_{\ev}-\lambda_{\ev}\big)(z) = 1$. Par l'équation (\ref{actionmaprop}), on voit que l'action de ce $z$ sur $U$ est l'identité. De plus, $U$ étant un produit tensoriel de deux représentations complètement réductibles, $U$ est lui-même une représentation de $\g^M$ qui est complètement réductible. La proposition \ref{suitespectrale} (1) donne alors que 
$$ \cohomo{\LL}{U} \cong \Hom_{\g^M}(K_M/K_M',U) $$
Supposons ensuite que $\lambda_{\ev} = \mu_{\ev}$. Dans ce second et dernier cas, $k_{\mu_{\ev} - \lambda_{\ev}} \cong k_0$, et alors $U \cong E_f^* \otimes E_g$.  L'équation (\ref{actionmaprop}) donne évidemment que $\Z.U = \{0\}$. La proposition \ref{suitespectrale} (2) justifie qu'on ait la suite exacte
$$ 0 \longrightarrow \cohomo{\g^M}{U} \stackrel{inf}{\longrightarrow} \cohomo{\LL}{U} \stackrel{res}{\longrightarrow} \Hom_{\g^M}(K_M/K_M',U) \longrightarrow \HHH^2(\g^M;U) $$
Puis, tel que noté dans les quelques lignes qui suivent l'équation \ref{ecriturespdgyo}, $U^{\g^M} = \{0\}$ dans ce contexte. Le théorème 10 de \cite{HSlie} donne alors que $\cohomo{\g^M}{U} \cong \{0\}$ et que $\HHH^2(\g^M;U) \cong \{0\}$. Sachant cela, la suite exacte ci-haut se réécrit
$$ 0 \longrightarrow \cohomo{\LL}{U} \stackrel{res}{\longrightarrow} \Hom_{\g^M}(K_M/K_M',U) \longrightarrow 0 $$
C'est donc que là aussi, on a un isomorphisme naturel 
$$ \cohomo{\LL}{U} \cong \Hom_{\g^M}(K_M/K_M',U) $$

Étudions la suite exacte courte de $\g^M$-modules suivante : 
$$ 0 \rightarrow N_M/K_M' \rightarrow K_M/K_M' \rightarrow K_M/N_M \rightarrow 0 $$
Le foncteur contravariant exact à gauche $\Hom_{\g^M}(-,U)$ donne alors lieu à une suite exacte de $\g^M$-modules (voir par exemple \cite{Rotman}, corollaire 6.62 et théorème 6.67) qui commence par
\begin{align*}
0 \rightarrow \Hom_{\g^M}(K_M/N_M,U&) \rightarrow \Hom_{\g^M}(K_M/K_M',U) \rightarrow \\
&\Hom_{\g^M}(N_M/K_M',U) \rightarrow \Ext_{\g^M}^1(K_M/N_M,U) 
\end{align*}
Par la proposition \ref{structureKab}, on sait que $K_M/N_M$ est un $\g^M$-module trivial. Conséquemment, si $\alpha \in \Hom_{\g^M}(K_M/N_M,U)$, alors $\im \alpha \subseteq U^{\g^M} = \{0\}$ ; ce qui implique que $\alpha = 0$ et donc, que 
\begin{equation} \label{premierpas}
\Hom_{\g^M}(K_M/N_M,U) \cong \{0\}
\end{equation}

D'un autre côté, le fait que $U$ soit complètement réductible combiné à $U^{\g^M} = \{0\}$ donne qu'en écrivant $U = \bigoplus_{k=0}^r V_k$ en somme directe finie de $\g^M$-modules simples de dimension finie, on trouve que chacun des $V_k$ est un $\g^M$-module simple non-trivial de dimension finie. Maintenant, puisque le foncteur $\Ext_{\g^M}^1(K_M/N_M,-)$ est additif en regard des sommes directes finies, la proposition \ref{mapremierepropyo} permet de conclure que
\begin{equation} \label{secondpas}
\Ext_{\g^M}^1(K_M/N_M,U) \cong \bigoplus_{k=0}^r \Ext_{\g^M}^1(K_M/N_M,V_k) \cong  \bigoplus_{k=0}^r \{0\} \cong \{0\}
\end{equation}

En utilisant les lignes (\ref{premierpas}) et (\ref{secondpas}), on voit que la dernière suite exacte peut en fait s'écrire
$$ 0 \rightarrow \Hom_{\g^M}(K_M/K_M',U) \rightarrow \Hom_{\g^M}(N_M/K_M',U) \rightarrow 0 $$
... on obtient ainsi des isomorphismes naturels
$$ \Ext_\LL^1(k_{\lambda_{\ev}} \otimes E_f,k_{\mu_{\ev}} \otimes E_g) \cong \Hom_{\g^M}(K_M/K_M',U) \cong \Hom_{\g^M}(N_M/K_M',U) $$
Le but de la proposition est donc atteint. \qed
\end{pre}
\end{Prop}

\begin{Rem} Dans la preuve ci-haut, il a d'abord été montré que si $E$ et $F$ étaient deux modules simples d'évaluations non-isomorphes supportés sur une seule orbite, alors 
$$ \Ext_\LL^1(E,F) \cong \Hom_{L/K_M}(K_M/K_M',E^* \otimes F) $$
Ce résultat intermédiaire apparaît comme un analogue à la situation des modules de \linebreak dimension 1. Voir, à ce propos, le corollaire \ref{extdim1}.  
\end{Rem}

\begin{Rem} Ce qu'il y a de bien au niveau théorique avec la dernière proposition, c'est que le $N_M/K_M'$ qui y apparait est un $\g^M$-module de dimension finie selon la proposition \ref{structureKab}. 

Notons aussi que ce module $N_M/K_M'$ ne dépend que de l'orbite de $M$ dans $\Max S$.
\end{Rem}

Soit $M \in \Max S$ et soit $N_M/K_M'$, le $\g^M$-module de dimension finie qui apparait dans la proposition précédente. Alors, il est complètement réductible et on peut l'écrire comme une somme de $\g^M$-modules irréductibles
\begin{equation} \label{constituantsnsurkprime}
N_M/K_M' = \bigoplus_{i=1}^n C_i
\end{equation}

\begin{Cor} \label{nsurkprimepourext} Supposons que $E$ et $F$ soient deux objets simples de $\LL_{\ev}$-$\mathbf{mod}$ qui soient non-isomorphes et qui soient supportés sur la seule $\Gamma$-orbite de $M \in \Max S$. 

Notons la décomposition de $E^* \otimes F$ en somme directe de $\g^M$-modules irréductibles
\begin{equation} \label{constituantstenseurevalsimples}
E^* \otimes F = \bigoplus_{j=1}^m D_j
\end{equation}

Alors
$$ \Ext_\LL^1(E,F) \cong \bigoplus_{i=1}^n \bigoplus_{j=1}^m \Hom_{\g^M}(C_i,D_j) $$
où les $C_i$ sont facteurs irréductibles de $N_M/K_M'$ tel que donné à la ligne (\ref{constituantsnsurkprime}).
\begin{pre} Ceci suit directement de la proposition précédente et du fait que le bifoncteur $\Hom_{\g^M}(-,-)$ est additif dans ses deux variables. Ceci est prouvé, par exemple, dans \cite{Rotman} à l'occasion des théorèmes 2.30 et 2.31. \qed  
\end{pre}
\end{Cor}

\begin{Cor} Supposons que $E$ et $F$ soient deux $\LL$-modules simples d'évaluation et de dimension finie qui soient non-isomorphes et qui soient supportés sur la seule $\Gamma$-orbite de $M \in \Max S$.

Suivant les notations des équations (\ref{constituantsnsurkprime}) et (\ref{constituantstenseurevalsimples}), on peut conclure que
$$ \Ext_\LL(E,F) = \{0\} \quad\; \Longleftrightarrow \;\quad \big\{[C_i]\big\}_{i=1}^n \cap \big\{[D_j]\big\}_{j=1}^m = \emptyset $$
\begin{pre} Ceci suit directement du corollaire précédent ainsi que de l'indispensable corollaire au lemme de Schur, le corollaire \ref{dimhom}. \qed
\end{pre}
\end{Cor}

\begin{Rem} Les trois derniers résultats montrent que pour connaître avec précision la structure des blocs d'extensions, il faudrait pouvoir connaître les facteurs irréductibles du module $N_M/K_M'$ associé à une orbite $\Gamma.M \subseteq \Max S$ arbitraire. 

En pratique, il semble ardu de seulement calculer des algèbres $\g^M$ lorsque celles-ci ne sont pas égales à $\g$ en entier (comme c'est le cas, par exemple, pour les formes tordues).
\end{Rem}

Il serait donc très utile à tous les niveaux de pouvoir réussir à déterminer les facteurs irréductibles du module $N_M/K_M'$ associée à une orbite $\Gamma.M$ arbitraire. Cependant, cela semble très difficile à faire car outre la proposition \ref{structureKab}, on ne sait pas grand chose sur ce $\g^M$-module de dimension finie.

Traitons maintenant le cas où les deux objets simples de $\LL_{\ev}$-$\mathbf{mod}$ supportés sur une seule et même $\Gamma$-orbite de $\Max S$ sont isomorphes. Autrement dit, si $V$ un objet simple de $\LL_{\ev}$-$\mathbf{mod}$ supporté sur une seule orbite, comment décrire $\Ext_\LL^1(V,V)$ ? 

Le résultat qui suit précise un petit peu cette question et reste indépendant de la propriété de localité vis-à-vis des extensions. 

\begin{Prop} \label{extiso} Soit $\LL$ une algèbre de courants tordue et soit $V$ un objet simple de \linebreak $\LL_{\ev}$-$\mathbf{mod}$ supporté sur la seule $\Gamma$-orbite de $M \in \Max S$. 

Posons $\g^M = \Z_M \oplus \SSS_M$ où $\Z_M = Z(\g^M)$ et $\SSS_M = [\g^M,\g^M]$. Posons $\mathfrak{Z}_M = \ev_M^{-1}(\Z_M)$ ; l'idéal de $\LL$ tel que $\ov{\ev}_M : \LL/\mathfrak{Z}_M \stackrel{\cong}{\longrightarrow} \SSS_M$. 

Alors
$$ \Ext_\LL^1(V,V) \cong \Hom_{\SSS_M}(\mathfrak{Z}_M/\mathfrak{Z}_M',V^* \otimes V) $$
\begin{pre} Par le corollaire \ref{sommelie}, on peut écrire $V$ sous la forme
$$ V \cong k_\lambda \otimes V_{\SSS_M} $$
où $\Z$ agit comme 0 sur $V_\SSS$. Par conséquent, $\Z$ agira nécessairement comme 0 sur 
$$ V^* \otimes V \cong k_{\lambda - \lambda} \otimes V_{\SSS_M}^* \otimes V_{\SSS_M} \cong V_{\SSS_M}^* \otimes V_{\SSS_M} $$
En particulier, $V^* \otimes V$ est un $\SSS_M$-module et $\mathfrak{Z}_M$ agit dessus comme 0. Aussi, puisque $\LL/\mathfrak{Z}_M$ est une algèbre de Lie semi-simple de dimension finie, la partie (1) de la proposition \ref{suitespectrale} s'applique et donne directement
$$ \cohomo{\LL}{V^* \otimes V} \cong \Hom_{\SSS_M}(\mathfrak{Z}_M/\mathfrak{Z}_M',V^* \otimes V) $$
Ceci est équivalent au résultat recherché par la proposition \ref{cohomoext}. \qed
\end{pre} 
\end{Prop}

À la lumière du théorème \ref{exteval}, de la proposition \ref{madeuxiemepropyo}, de la proposition \ref{extiso} et de la classification des blocs d'extensions d'une algèbre de courants tordue locale vis-à-vis des extensions, c'est possible d'écrire le résultat suivant :

\begin{Thm} \label{thmultime} Soit $\LL$ une algèbre de courants tordue locale vis-à-vis des extensions et fixons une décomposition 
$$ (\LL/\LL')^* = (\LL/\LL')_{\ev}^* \oplus (\LL/\LL')_{\nev}^* $$
Soient $V$ et $W$ deux objets simples de $\LL$-$\mathbf{mod}$. Sans perte de généralité, on peut écrire 
\begin{align*} 
V \cong k_\lambda \otimes E_f && W \cong k_\mu \otimes E_g
\end{align*}
où $(\lambda,f)$ et $(\mu,g)$ sont des couples de $\LL^* \times \F^\Gamma$ qui respectent les trois conditions de la proposition \ref{classmodules}. Écrivons aussi $\lambda = \lambda_{\ev} + \lambda_{\nev}$ et $\mu = \mu_{\ev} + \mu_{\nev}$ et
\begin{align*} 
V_{\ev} = k_{\lambda_{\ev}} \otimes E_f && W_{\ev} = k_{\mu_{\ev}} \otimes E_g
\end{align*}

Fixons $\{M_i\}_{i=1}^r \subseteq \Max S$ un ensemble complet de représentants des $r$ $\Gamma$-orbites distinctes de $\supp V_{\ev} \cup \supp W_{\ev}$. C'est donc que $V_{\ev} = \bigotimes_{i=1}^r V_i$ et  $W_{\ev} = \bigotimes_{i=1}^r W_i$ où les $V_i$ et $W_i$ sont des $\g^{M_i}$-modules simples. 

Soient $f_{V_{\ev}}$ et $f_{W_{\ev}}$ les sections $\Gamma$-invariantes à support fini de $\R$ correspondant aux classes d'isomorphismes de $V_{\ev}$ et $W_{\ev}$ respectivement. Alors

\begin{itemize}
\item[\em(i)] Si $(\lambda_{\nev},\chi_{V_{\ev}}) \neq (\mu_{\nev},\chi_{W_{\ev}})$, alors 
$$ \Ext_\LL^1(V,W) \cong \{0\} $$
\item[\em(ii)] Si $(\lambda_{\nev},\chi_{V_{\ev}}) = (\mu_{\nev},\chi_{W_{\ev}})$, alors :

\begin{itemize}
\item[\em(1)] si $f_{V_{\ev}}$ et $f_{W_{\ev}}$ diff\`erent sur plus d'une $\Gamma$-orbite de $\Max S$, alors 
$$ \Ext_\LL^1(V,W) \cong \{0\} $$
\item[\em(2)] si $f_{V_{\ev}}$ et $f_{W_{\ev}}$ ne diff\`erent que sur la $\Gamma$-orbite de $M_j \in \Max S$, alors 
$$ \Ext_\LL^1(V,W) \cong \Hom_{\g^{M_j}}(N_{M_j}/K_{M_j}',V_j^* \otimes W_j) $$
\item[\em(3)] Si si $f_{V_{\ev}} = f_{W_{\ev}}$, alors 
$$ \big((\LL/\LL')^*\big)^{\oplus\,r-1} \oplus \Ext_\LL^1(V,W) \cong \bigoplus_{i=1}^r \Hom_{\SSS_{M_i}}(\mathfrak{Z}_{M_i}/\mathfrak{Z}_{M_i}',V_i^* \otimes W_i) $$
\end{itemize}
\end{itemize}
\end{Thm}

\begin{Rem} \label{aubordduprecipice} Si quelqu'un arrivait à montrer, par exemple, que l'action de $Z(\g^M)$ sur le $\g^M$-module est triviale quelque soit l'idéal $M \in \Max S$, le théorème précédent serait beaucoup plus esthétique. 

En pratique, ceci ne risquerait toutefois pas de rendre les choses plus simples étant donnée la difficulté de calculer des ensembles d'homomorphismes tels ceux qui interviennent dans le théorème ci-haut.
\end{Rem}

\begin{Rem} Une autre conséquence directe du corollaire \ref{nsurkprimepourext} est la suivante. Soit $M \in \Max S$ et soit $N_M/K_M'$, le module associé à l'orbite $\Gamma.M$. Il importe de noter que c'est un $\g^M$-module. Ainsi, les classes d'isomorphismes de ses facteurs simples $\{[C_i]\}_{i=1}^n$ (en tant que $\g^M$-modules) sont tous dans le bloc d'extensions de $\LL$-$\mathbf{mod}$, celui du module trivial de dimension 1. 

En bref, pour tout $i \in \{1,...\,,n\}$, il est avéré que
$$ \llbracket C_i \rrbracket = \llbracket k_0 \rrbracket \quad \text{ dans } \LL\text{-}\mathbf{mod} $$

Notons que le résultat qui fait l'objet de la présente remarque peut aussi être interprété comme étant une conséquence du théorème \ref{thmultime}.
\end{Rem}             
\chapter{Applications}     
Dans ce chapitre, il est question de montrer quelques exemples plus concrets où la classification des extensions des algèbres de courant tordues permet soit de retrouver des résultats existants, soit d'exposer des résultats nouveaux.

Une des idées derrière l'étude récente des algèbres de courants tordues était que cette notion généralisait à la fois les formes tordues et les algèbres d'applications équivariantes.

\section{Algèbres d'applications équivariantes et formes tordues}
Cette section montre d'abord que les résultats de ce travail concordent avec les résultats de l'article \cite{NSext}, article qui est à la base du présent mémoire et dont plusieurs résultats sont repris. Ensuite, on s'attarde aux formes tordues qui sont le sujet de l'article \cite{LPtfa}.
\subsection{Algèbres d'applications équivariantes}
Certaines algèbres d'applications équivariantes concordent avec la notion d'algèbres de courants tordues de la section 2.1.1 tel qu'expliqué dans l'introduction de l'article  \cite{LauTCA}. Une algèbre d'applications équivariantes pour laquelle $\g$ est une algèbre de Lie simple est ce dont il est ici question.

Les résultats des derniers chapitres concordent avec les résultats obtenus par E.Neher et A.Savage dans leur article \cite{NSext} ; ce qui est normal puisque le présent mémoire se base sur cet article. Dans le contexte des algèbres d'applications équivariantes, les deux auteurs E.Neher et A.Savage ont pu exploiter le fait que l'action de $\Gamma$ sur $\g \otimes S$ était donnée directement par des actions $\Gamma \curvearrowright \g$ et $\Gamma \curvearrowright S$ pour préciser les résultats lorsque $\Gamma$ est abélien. Il n'a pas été possible d'en faire autant pour les algèbres de courants tordues.

Toute la section 4 de \cite{NSext} est dédiée à préciser les résultats sur les extensions des modules de $\LL$-$\mathbf{mod}$ lorsque $\LL$ est une algèbre d'applications équivariantes avec un groupe $\Gamma$ abélien. La forme précise de l'action leur prermet de faire des liens entre les $\Gamma$-gradations de $\g$, de $S$ et l'algèbre $\LL = (\g \otimes S)^\Gamma$. Ces méthodes sont cependant inapplicables dans la généralité des algèbres de courants tordues. J'ai tenté, en vain, d'obtenir des résultats similaires pour les algèbres de courants tordues, mais les résultats clés de la démarche de la section 4 de \cite{NSext} semblent être les équations (4.7) dans l'article. Celles-ci sont impossibles à recréer directement dans le cadre des algèbres de courants tordues. 

\begin{Rem} Je ne crois cependant pas qu'un résultat comme la proposition 4.9 de \cite{NSext} est totalement hors de portée pour certaines formes tordues. Le rapprochement manifeste qu'on peut faire entre la proposition 4.9 et la proposition \ref{trucmargaux} (de l'annexe B de ce mémoire) ne semble pas complètement banal. 
\end{Rem}

Un analogue à l'isomorphisme $\g^x \otimes (I/I^2)_0 \cong (K_x/D_x)_{\Xi_x}$ dans l'énoncé du théorème 4.5 de \cite{NSext} permettrait sans doute, dans le cas où $\Gamma$ est abélien et agit librement sur $\Max S$, de conclure un résultat semblable au corollaire 4.8 en suivant la démarche de la proposition 4.7. C'est que sous ces hypothèses, et dans le cas des algèbres d'applications équivariantes, on a, d'après le lemme (4.7), $K_x/D_x \cong \g^x \otimes (I/I^2)_0$. Ainsi, supposant qu'un tel isomorphisme soit encore valable dans le cas d'une algèbre de courants tordue, toute la preuve de la proposition 4.7 fonctionnerait telle quelle.

Quoi qu'il en soit, cette incapacité à adapter des méthodes pour obtenir des résultats comme ceux de la section 4 de \cite{NSext} est une des raisons pour lesquelles la section 6 de ce même article est beaucoup plus impressionnante que ce chapitre de mémoire.  

\subsection{Formes tordues}
Les algèbres de Lie qui sont des formes tordues (notion définie dans \cite{LPtfa}) cadrent parfaitement avec la notion d'algèbre de courants tordue de la section 2.1.1 de ce document. En fait, les hypothèses de base de l'article \cite{LauTCA} sur la classification des modules simples de dimension finie pour une algèbre de courants tordue sont seulement un peu plus faibles que les hypothèses de base dans \cite{LPtfa}.

Dans le cas des formes tordues, l'extension $S/S^\Gamma$ est même galoisienne. La définition d'une extension d'anneaux galoisienne est donnée en page six de \cite{CHR}, voir la définition 1.4. Ceci donne accès à de nombreux résultats plus précis pour les formes tordues ; notamment sur la structure de ces algèbres et de leurs modules. 

Il est vite établi dans \cite{LPtfa} que toutes les représentations irréductibles pour une forme tordue sont des représentations d'évaluation. Ceci vient du fait que les formes tordues sont parfaites (voir la preuve du lemme 2.7 dans \cite{LPtfa}). Ainsi, d'après la remarque \ref{parfaiteextlocal}, ces algèbres sont automatiquement locales vis-à-vis des extensions. Si $\LL$ est une forme tordue, les blocs d'extensions de $\LL$-$\mathbf{mod} = \LL_{\ev}$-$\mathbf{mod}$ sont donc classifiés par les caractères spectraux $\carspec^\Gamma$ comme le montrent les résultats de la section 2.3.2. Plus précisément, voici le résultat : 

\begin{Thm} Soit $\LL = (\g \otimes S)^\Gamma$ une forme tordue. Alors la catégorie $\LL$-$\mathbf{mod}$ des \linebreak $\LL$-modules simples de dimension finie se décompose en blocs comme suit :
$$ \text{$\LL$-$\mathbf{mod}$} = \bigoplus_{\chi \in \carspec^\Gamma} \Cat^\chi $$
En particulier, les blocs et blocs d'extensions de cette catégorie sont en bijection avec les caractères spectraux $\carspec^\Gamma$.
\end{Thm}

Une autre propriété notable des formes tordues est que toutes les applications d'évaluation $\ev_M : \LL \rightarrow \g$ (où $M \in \Max S$) sont surjectives (voir \cite{LPtfa}, proposition 3.3). C'est donc que tous les $\g^M$ qui apparaissent dans la théorie développée dans le chapitre 2 peuvent être remplacés par $\g$ dans les formules. Il faut cependant bien faire attention de toujours pouvoir distinguer un $\g = \g^{M_1}$ d'un $\g = \g^{M_2}$ si $M_1 \neq M_2 \in \Max S$ par exemple.

Dans le cas des formes tordues, certains calculs se rapprochent davantage. Dans l'énoncé de la proposition \ref{extiso} par exemple, on fixe un $M \in \Max S$ et il y est définit l'idéal \linebreak $\mathfrak{Z}_M = \ev_M^{-1}\big(Z(\g^M)\big)$ de $\LL$. Puisque $\LL$ est une forme tordue, $\g^M = \g$ est une algèbre de Lie simple, donc $Z(\g^M) = Z(\g) = \{0\}$. Ainsi $\mathfrak{Z}_M = \ker(\ev_M) = K_M$ et le résultat suivant devient accessible .

\begin{Prop} Soit $V$ un objet simple de $\LL$-$\mathbf{mod}$ supporté sur la seule $\Gamma$-orbite de l'idéal $M \in \Max S$. Alors, 
$$ \Ext_\LL^1(V,V) \cong \Hom_{\g}(K_M/K_M',V^* \otimes V) $$
\begin{pre} La partie (1) de la proposition \ref{suitespectrale} s'applique et donne directement le résultat. \qed
\end{pre} 
\end{Prop}

Sous les mêmes hypothèses, il y a entre $\Ext_\LL^1(V,V)$ et $\Hom_\g(N_M/K_M',V^* \otimes V)$ un lien un peu plus étroit. 

\begin{Cor} \label{montroisiemeresultatyo} Soit $V$ un objet simple de $\LL$-$\mathbf{mod}$ supporté sur la seule $\Gamma$-orbite de l'idéal $M \in \Max S$. Alors, il y a un isomorphisme d'espaces vectoriels
$$ \Ext_\LL^1(V,V) \cong (K_M/N_M)^* \oplus \Hom_\g(N_M/K_M',V^* \otimes V) $$
Notons que plusieurs tels isomorphismes existent et qu'ils dépendent d'un choix.
\begin{pre} Rappelons que le corollaire \ref{dimhom} donne que $(V^* \otimes V)^\g = k_0$. Puisque $V^* \otimes V$ est un $\g$-module de dimension finie, il est complètement réductible alors il est possible d'écrire
$$ V^* \otimes V = k_0 \oplus A $$
où $A$ est un $\g$-module de dimension fini tel que $A^\g = \{0\}$. 

C'est encore vrai que la suite exacte courte $\, 0 \rightarrow N_M/K_M' \rightarrow K_M/K_M' \rightarrow K_M/N_M \rightarrow 0$ induit une suite exacte de $\g$-modules commençant par :
\begin{align*}
0 \rightarrow \Hom_\g(K_M/N_M,V^* \otimes V&) \rightarrow \Hom_\g(K_M/K_M',V^* \otimes V) \rightarrow \\
&\Hom_\g(N_M/K_M',V^* \otimes V) \rightarrow \Ext_\g^1(K_M/N_M,V^* \otimes V) 
\end{align*}
Étant donné que $K_M/N_M$ est un $\g$-module trivial, la proposition \ref{mapremierepropyo} donne, ici aussi, que $\Ext_\g^1(K_M/N_M,V^* \otimes V) \cong \{0\}$ et puis, on a aussi
\begin{align*}
\Hom_\g(K_M/N_M,V^* \otimes V) &\cong \Hom_\g(K_M/N_M,k_0) \oplus \Hom_\g(K_M/N_M,A)\\
&\cong \Hom_\g(K_M/N_M,k_0)\\
&\cong \big((K_M/N_M)^*\big)^\g
\end{align*}
Encore une fois, puisque $K_M/N_M$ est un $\g$-module trivial, le module $(K_M/N_M)^*$ l'est aussi et alors 
$$ \big((K_M/N_M)^*\big)^\g = (K_M/N_M)^* $$ 
En général, $K_M/N_M$ est (sans doute) de dimension infinie alors on doit en rester là. Néanmoins, les faits établis jusqu'ici permettent d'écrire la suite exacte de $\g$-modules suivante
$$ 0 \rightarrow (K_M/N_M)^* \rightarrow \Hom_\g(K_M/K_M',V^* \otimes V) \rightarrow \Hom_\g(N_M/K_M',V^* \otimes V) \rightarrow 0 $$
En particulier, on a la suite exacte d'espaces vectoriels 
$$ 0 \rightarrow (K_M/N_M)^* \rightarrow \Ext_\LL^1(V,V) \rightarrow \Hom_\g(N_M/K_M',V^* \otimes V) \rightarrow 0 $$
Comme toutes les suites d'espaces vectoriels sont scindées, l'isomorphisme voulu existe, bien qu'il ne soit ni canonique, ni unique. \qed
\end{pre}
\end{Cor}

Pour conclure cette sous-section sur les formes tordues, voici comment c'est possible d'interpréter le théorème \ref{thmultime} pour le cadre théorique des formes tordues d'algèbres de Lie :

\begin{Thm} Soit $\LL = (\g \otimes S)^\Gamma$ une forme tordue et soient $V$ et $W$ deux objets simples de $\LL$-$\mathbf{mod}$. Sans perte de généralité, on peut écrire $V \cong E_f$ et $W \cong E_g$ où $f$ et $g$ sont des sections de $\F^\Gamma$.

Fixons $\{M_i\}_{i=1}^r \subseteq \Max S$ un ensemble complet de représentants des $r$ $\Gamma$-orbites distinctes de $\supp V \cup \supp W$. C'est donc que $V = \bigotimes_{i=1}^r V_i$ et  $W = \bigotimes_{i=1}^r W_i$ où les $V_i$ et $W_i$ sont des $\g = \g^{M_i}$-modules simples. Alors

\begin{itemize}
\item[\em(i)] Si $\chi_V \neq \chi_W$, alors 
$$ \Ext_\LL^1(V,W) \cong \{0\} $$
\item[\em(ii)] Si $\chi_V = \chi_W$, alors :

\begin{itemize}
\item[\em(1)] si $f$ et $g$ diff\`erent sur plus d'une $\Gamma$-orbite de $\Max S$, alors 
$$ \Ext_\LL^1(V,W) \cong \{0\} $$
\item[\em(2)] si $f$ et $g$ ne diff\`erent que sur la $\Gamma$-orbite de $M_j \in \Max S$, alors 
$$ \Ext_\LL^1(V,W) \cong \Hom_{\g}(N_{M_j}/K_{M_j},V_j^* \otimes W_j) $$
\item[\em(3)] si $f = g$, alors 
$$ \Ext_\LL^1(V,W) \cong \bigoplus_{i=1}^r \big((K_{M_i}/N_{M_i})^* \oplus \Hom_\g(N_{M_i}/K_{M_i}',V_i^* \otimes W_i)\big) $$
\end{itemize}
\end{itemize}
\end{Thm}

\section{Algèbres de Margaux}
\subsection{Les algèbres de Margaux et leur intérêt théorique}

Les formes tordues incluent l'importante classe d'algèbres de Lie des algèbres de multilacets lorsque $\g$ est simple. La notion générale de forme tordue comprend également des algèbres de Lie qui ne sont pas des algèbres de multilacets ; elles portent l'appellation d'algèbres de Margaux. Ceci est mentionné dans la section 4.2 de l'article \cite{LPtfa}. Cette observation vient de l'article \cite{GPmar} de P.Gille et A.Pianzola, un article de 2007. Les algèbres de Margaux restent assez énigmatiques, mais pourtant, leurs modules simples de dimension finie et les blocs d'extensions de leur catégorie de modules de dimension finie sont classifiés. L'intéret théorique de ces algèbres réside principalement dans le fait qu'elles sont des formes tordues, mais sans pour autant être des algèbres de multilacets.

En pratique, les algèbres de Margaux sont assez difficiles à appréhender, mais la construction de l'une d'entre elles est connue. Celle dont on connait une construction est probablement la moins complexe d'entre elles. Appelons cette algèbre $\MM$. En bref, l'algèbre de Margaux $\MM$ est une $\CC[t_1^{\pm 1},t_2^{\pm 1}] / \CC[t_1^{\pm 2},t_2^{\pm 2}]$ forme tordue de l'algèbre de Lie $\mathfrak{sl}_2\big(\CC[t_1^{\pm 2},t_2^{\pm 2}])$. Sa construction est brièvement décrite dans la section 4.2 de \cite{LPtfa} et aussi dans la section 5 de \cite{GPmar} (voir l'exemple 5.7). 

Dans ce qui suit, la construction connue de cette algèbre de Margaux $\MM$ est expliquée. Ne serait-ce que pour justifier cette lecture, il sera bon de garder en tête que la théorie exhibée dans la section 3.1.2 s'applique à cette algèbre. Il sera fait mention dans la sous-section subséquente, de comment il serait possible de préciser la description des blocs d'extension de la catégorie des $\MM$-$\mathbf{mod}$. Il convient de noter que \textbf{ces résultats sont nouveaux}.

Suivant les explications de \cite{LPtfa}, l'algèbre associative $M_2(\CC)$ et l'algèbre de Lie $\mathfrak{sl}_2(\CC)$ ont des groupe d'automorphismes isomorphes ; ceci fait en sorte qu'il y a une correspondance entre les $S/R$-formes tordues de $M(R)$, des algèbres associatives, et celles de $\mathfrak{sl}_2(R)$, des algèbres de Lie. Pour commencer, rappelons que toute algèbre associative $A$ est une algèbre de Lie avec le commutateur pour crochet. Notons $\lie(A)$ l'algèbre de Lie sur $A$ qui a pour crochet le commutateur. Ensuite, supposons maintenant que $A$ soit une $S/R$-forme de $M_2(R)$, alors
\begin{equation} \label{processussl2}
\big(\lie(A)\big)' \text{ est une } S/R \text{-forme tordue de } \mathfrak{sl}_2(R)
\end{equation}

C'est par cette procédure que sera construite l'algèbre de Margaux $\MM$.

\subsection{Construction de $\MM$}
D'abord, voici la notation de base qui sera utilisée dans le reste du chapitre :
\begin{itemize}
\item[\textbullet] $\MM$ est l'algèbre de Margaux dont il est question à la section 3.2.1.
\item[\textbullet] $k = \CC$.
\item[\textbullet] $\g = \mathfrak{sl}_2(\CC)$.
\item[\textbullet] $S = \CC[t_1^{\pm 1},t_2^{\pm 1}]$.
\item[\textbullet] $\Gamma = \ZZ_2 \times \ZZ_2$ et l'action de $(\ov{a},\ov{b}) \in \Gamma$ sur $S$ est donnée par extension $\CC$-linéaire de l'application de $S$ dans $S$ qui a pour effet
\begin{align*} 
t_1 &\mapsto (-1)^a \cdot t_1\\
t_2 &\mapsto (-1)^b \cdot t_2
\end{align*}
\end{itemize}

Cest hypothèses entraînent que $S^\Gamma = \CC[t_1^{\pm 2},t_2^{\pm 2}]$ et aussi que l'extension d'anneaux $S/S^\Gamma$ est galoisienne. Pour faire le lien avec la notation de (\ref{processussl2}) dans notre cas, on posera \linebreak $R = S^\Gamma = \CC[t_1^{\pm 2},t_2^{\pm 2}]$.

Il est bon de s'assurer de comprendre l'action de $\Gamma$ sur $\Max S$. Grâce au Nullstellensatz, (voir le théorème A.1.3 de \cite{GWsym} par exemple), nous savons que 
$$ \Max S = \{\,M_{a,b} = \; \langle \, t_1 - a,t_2 - b \, \rangle \unlhd \; S \; | \; a,b \in \CC^\times \} $$ 

Supposons que $M_{a,b} \in \Max S$ et que $(\ov{u},\ov{v}) \in \ZZ_2 \times \ZZ_2 = \Gamma$. On peut calculer que
\begin{align*}
{}^{(\ov{u},\ov{v})} M_{a,b} &= \; \langle \, t_1- (-1)^{u} \cdot a,t_2- (-1)^{v} \cdot b \, \rangle\\
&= M_{(-1)^{u} \cdot a,(-1)^{v} \cdot b}
\end{align*}
C'est donc que la $\Gamma$-orbite de $M_{a,b} \in \Max S$ est
\begin{equation}
\Gamma.M_{a,b} = \{M_{a,b},M_{-a,b},M_{a,-b},M_{-a,-b}\} \subseteq \Max S 
\end{equation}

Maintenant, si $M_{a,b} \in \Max S$, alors on a que
$$ M_{a,b} \cap S^\Gamma = \; \langle \, t_1^2 - a^2,t_2^2 - b^2 \, \rangle \; \in \Max S^\Gamma $$

\begin{Rem} En somme, un idéal maximal de $\Max S$ est la donnée d'un couple \linebreak $(a,b) \in (\CC^\times)^2$ et de plus, on a que
$$ \Gamma.M_{a,b} \neq \Gamma.M_{c,d} \quad \Longleftrightarrow \quad (a^2,b^2) \neq (c^2,d^2) \in (\CC^\times)^2 $$
\end{Rem}

Pour pouvoir donner la construction de l'algèbre de Margaux $\MM$, il faut d'abord introduire une algèbre associative particulière :

\begin{Def} \label{quatgen} Soit $\Q$, la $\CC[t_1^{\pm 2},t_2^{\pm 2}]$-algèbre générée par les symboles $T_1$ et $T_2$ qui soient sujet aux relations suivantes :
\begin{align*}
T_1^2 &= t_1^2\\
T_2^2 &= t_2^2\\
T_1T_2 &= -T_2 T_1 
\end{align*}
\end{Def}

\begin{Rem} \label{quaterniongen} $\Q$ est une algèbre de quaternions généralisée et c'est une algèbre d'Azumaya. Voir à ce propos l'exercice 8.10 (b) dans \cite{FDnoncom}.

Cet exercice montre aussi que si on pose $R = \CC[t_1^{\pm 2},t_2^{\pm 2}]$, alors $\Q$ est un $R$-module libre et l'ensemble $\{1,T_1,T_2,T_1T_2\}$ en forme une base.

Pour la définition d'une algèbre d'Azumaya, voir le tout début du chapitre 8 de la même référence. En particulier, $\Q$ est un $R$-module fidèle.
%
\end{Rem}

\begin{Rem} \label{cocyclequaternion} L'algèbre associative $\Q$ est en fait elle-même une $S/R$-forme tordue de $M_2(R)$ où $S = \CC[t_1^{\pm 1},t_2^{\pm 1}]$ et $R = \CC[t_1^{\pm 2},t_2^{\pm 2}]$ selon la procédure (\ref{processussl2}). 

La $S/R$-forme tordue $\big(\lie(\Q)\big)'$ est en fait isomorphe à une algèbre de multilacets. Voir $\LL_1$ dans la section 4.2 de \cite{LPtfa} ou dans l'exemple 5.7 de \cite{GPmar}. Ceci est lié au fait que $\Q$ est un $\Q$-module libre... Un peu plus de détails à propos de l'algèbre de Lie $\LL_1$ seront donnés dans l'annexe A dans laquelle le 1-cocycle associé à cette forme tordue est calculé.  
\end{Rem}

Dans le cas de l'algèbre $\MM$, on remplacera l'algèbre associative $\Q$ dans (\ref{processussl2}) par l'algèbre des endomorphismes d'un $\Q$-module projectif, mais non libre. 
 
Posons
\begin{equation} \label{modulemargaux}
O = \{(x,y) \in \Q_\Q \oplus \Q_\Q \; | \; (1+T_1)x = (1+T_2)y\}
\end{equation} 
où $\Q_\Q$ indique que $\Q$ est vu comme un $\Q$-module à droite. Le module $O$ est donc un $\Q$-module à droite, et ce par définition même.

\begin{Prop} \label{projnonlibre} $O$ est un $\Q$-module à droite projectif mais non-libre.

\begin{pre} D'abord, l'application
\begin{align*} 
F : \; \Q_\Q &\oplus \Q_\Q \quad \longrightarrow \quad \Q_\Q\\
(x&,y) \mapsto (1+T_1)x - (1+T_2)y
\end{align*}
est surjective. En effet, un calcul direct montre que $F\left(\frac{1-T_2}{2},\frac{T_1-1}{2}\right) = 1$. De plus, pour tout $q \in \Q$, la définition de $F$ donne que 
$$ F(xq,yq) = F(x,y)q $$ 
C'est donc un homomorphisme surjectif de $\Q$-modules à droite.

Par définition, on a que $O = \ker F$ alors la suite suivante de $\Q$-modules à droite
$$ 0 \rightarrow O \rightarrow \Q \oplus \Q \stackrel{F}{\longrightarrow}\Q \rightarrow 0 $$
est exacte. Elle est également scindée et pour le montrer, on peut considérer l'application 
\begin{align*}
\mathbf{sc} : \qquad \quad &\Q \quad \longrightarrow \quad \Q \oplus \Q\\
&\;q \mapsto \left(\frac{1-T_2}{2} q,\frac{T_1-1}{2} q\right) = \left(\frac{1-T_2}{2},\frac{T_1-1}{2}\right)q
\end{align*}
Il est facile de vérifier que $\mathbf{sc}$ est un morphisme de $\Q$-modules à droite. C'est aussi une section de $F$ puisque si $q \in \Q$, alors 
\begin{align*}
\big(F \circ \mathbf{sc}\big)(q) &= F\left(\frac{1-T_2}{2} q,\frac{T_1-1}{2} q\right)\\
&= \left(F\left(\frac{1-T_2}{2},\frac{T_1-1}{2}\right)\right) \cdot q\\
&= 1 \cdot q\\
&= q
\end{align*}

En particulier, le $\Q$-module à droite $O$ est projectif.

Le fait que $O$ soit un $Q$-module à droite non-libre est établi par la proposition 3.20 de \cite{GPmar}. Leur preuve s'inspire de la preuve de la proposition 1 de \cite{OSproj}, un article de M.Ojanguren et R.Sridharan paru en 1971. \qed
\end{pre}
\end{Prop}

\begin{Rem} \label{modfidele} Puisque $\Q$ est une algèbre d'Azumaya (voir la remarque \ref{quaterniongen}), $\Q$ est fidèle en tant que $R$-module. C'est donc aussi le cas pour le $R$-module $\Q \oplus \Q$. De cela, on peut aussi vérifier montrer que $O$ est un $R$-module à droite fidèle. 
\end{Rem}

\begin{Rem} Combinant les conclusions de la proposition \ref{projnonlibre} et de la remarque \ref{modfidele}, il est établi que $O$ est un $\Q$-module de type fini, projectif et fidèle. La proposition 8.3 de \cite{FDnoncom} donne, ensuite, directement que $\End_\Q(O)$ est une $R$-algèbre d'Azumaya.

Ceci signifie, puisque $O \hookrightarrow \Q \oplus \Q$, que $\End_\Q(O)$ est une forme tordue pour $M_2(R)$. En fait, on peut même dire que c'est une $S/R$-forme tordue de $M_2(R)$.
\end{Rem}

\begin{Def} Suivant (\ref{processussl2}), l'algèbre de Margaux $\MM$ est définie comme
\begin{equation} \label{premiermargaux}
\MM = \Big(\lie\big(\End_\Q(O)\big)\Big)'
\end{equation}
\end{Def}

\begin{Rem} Le fait que $\MM$ ne soit pas une algèbre de multilacets est directement relié au fait que le module $O$ construit ci-haut est un $\Q$-module projectif, mais non-libre. Pour plus de détails à ce sujet, voyez les détails entourant l'exemple 5.7 de \cite{GPmar}.
\end{Rem}

\subsection{Description de $\MM$}

$\MM$ est une $S/R$-forme de $\mathfrak{sl}_2(R) \cong \mathfrak{sl}_2(\CC) \otimes R$. Cette description manque cependant un peu de précision. Ce qui suit vise à concrétiser, dans la mesure du possible, la description de $\MM$ par l'équation (\ref{premiermargaux}).

Dans la preuve de la proposition \ref{projnonlibre}, il est établi que la suite exacte 
$$ 0 \rightarrow O \rightarrow \Q \oplus \Q \stackrel{F}{\longrightarrow}\Q \rightarrow 0 $$ 
est scindée et on donne une section $\mathbf{sc} : \Q \rightarrow \Q \oplus \Q$. Ainsi, si on identifie $O = \ker F$ et $\Q = \im (\mathbf{sc})$, il y a un isomorphisme de $\Q$-module à droite 
\begin{equation} \label{decompositionisomar}
O \oplus \Q \; = \; \ker F \oplus \im (\mathbf{sc}) \; \cong \; \Q \oplus \Q
\end{equation}
On peut considérer le morphisme $(\mathbf{sc} \circ F) : \Q \oplus \Q \rightarrow \Q \oplus \Q$. Comme $F$ est surjective, on a $\im (\mathbf{sc}) = \im (\mathbf{sc} \circ F)$. Soit $(a,b) \in \Q \oplus \Q$, alors on peut écrire
\begin{align*}
&\big(\mathbf{sc} \circ F\big)(a,b) \quad = \quad \left(\frac{1-T_2}{2},\frac{T_1-1}{2}\right)F(a,b)\\
&\big(\mathbf{sc} \circ F\big): \; \Vv{a}{b} \mapsto \Vv{\frac{1-T_2}{2}}{\frac{T_1-1}{2}}\big((1+T_1)a - (1+T_2)b\big)\\
&\qquad \qquad \;\,\, \Vv{a}{b} \mapsto \Vv{\frac{1-T_2}{2}\big((1+T_1)a - (1+T_2)b\big)}{\frac{T_1-1}{2}\big((1+T_1)a - (1+T_2)b\big)}\\
&\qquad \qquad \;\,\, \Vv{a}{b} \mapsto \frac{1}{2}\matt{(1-T_2)(1+T_1)}{-(1-T_2)(1+T_2)}{(T_1-1)(1+T_1)}{(1-T_1)(1+T_2)} \Vv{a}{b}\\
&\qquad \qquad \;\,\, \Vv{a}{b} \mapsto \frac{1}{2}\matt{(1-T_2)(1+T_1)}{t_2^2-1}{t_1^2-1}{(1-T_1)(1+T_2)} \Vv{a}{b}
\end{align*}

\begin{Def} La matrice qui apparaît à la ligne précédente est assez importante dans la suite. Celle-ci aura donc un symbole ; 
$$ \Omega = \frac{1}{2}\matt{(1-T_2)(1+T_1)}{t_2^2-1}{t_1^2-1}{(1-T_1)(1+T_2)} $$
\end{Def}

\begin{Lem} La matrice $\Omega$ est \textbf{idempotente}, c'est à-dire qu'elle vérifie 
$$ \Omega^2 = \Omega \in M_2(\Q) $$
En particulier, la matrice $(\Id_2 - \Omega) \in M_2(\Q)$ est aussi idempotente.
\begin{pre} On pourrait faire les calculs directement à la main. Cela se fait aussi avec des logiciels de calculs spécialisés tel le logiciel gratuit « Sage ». 

Une autre façon de montrer que $\Omega^2 = \Id_2$ est de considérer la composition des applications $(\mathbf{sc} \circ F) \circ (\mathbf{sc} \circ F)$. Puisque $\mathbf{sc}$ est une section de $F$, on a $F \circ \mathbf{sc} = \Id : \Q \rightarrow \Q$, on peut écrire
$$ (\mathbf{sc} \circ F) \circ (\mathbf{sc} \circ F) = \mathbf{sc} \circ (F \circ \mathbf{sc}) \circ F = \mathbf{sc} \circ F \; : \quad \Q \oplus \Q \rightarrow \Q \oplus \Q $$ 
Puisque $\Omega$ est la matrice de $\mathbf{sc} \circ F$, on conclut que $\Omega^2 = \Id_2 \in M_2(\Q)$ ; il s'agit donc d'une matrice idempotente. Pour terminer cette preuve, on peut écrire
$$ (\Id_2 - \Omega)^2 = \Id_2^2 - 2\Omega + \Omega^2 = \Id_2 - 2\Omega + \Omega = \Id_2 - \Omega $$
Ainsi, la matrice $(\Id_2 - \Omega)$ est elle aussi idempotente.
\qed
\end{pre}
\end{Lem}

Soit donc $(a,b) \in \Q \oplus \Q$. On peut évidemment écrire
\begin{equation} \label{decompositionmatricemar}
\Vv{a}{b} = (\Id_2 - \Omega) \Vv{a}{b} + \; \Omega \Vv{a}{b}
\end{equation}

\begin{Lem} Le premier terme de la décomposition (\ref{decompositionmatricemar}) est dans $O = \ker F$ et le deuxième terme de cette décomposition est dans $\im (\mathbf{sc}) \cong \Q_\Q$.

En particulier, l'isomorphisme (\ref{decompositionisomar}) est rendu explicite par la décomposition (\ref{decompositionmatricemar}).
\begin{pre} D'après ce qui a été expliqué plus haut, et d'après la définition même de $\Omega$, le deuxième terme de (\ref{decompositionmatricemar}) est dans $\im (\mathbf{sc}) = \im (\mathbf{sc} \circ F)$.
Pour ce qui est du premier terme, il suffit de calculer qu'il est envoyé sur $0 \in \Q$ par $F$. C'est bien le cas puisque
\begin{align*}
\big(F \circ (\Id_2 - \Omega)\big) &= F \circ Id_2 - F \circ (\mathbf{sc} \circ F) \\
&= F - (F \circ \mathbf{sc}) \circ F\\
&= F - F\\
&= 0
\end{align*}
Le résultat est donc démontré. Ensuite, puique la décomposition (\ref{decompositionmatricemar}) est unique pour chaque élément $(a,b) \in \Q \oplus \Q$, c'est bel et bien l'isomorphisme de $\Q$-modules à droite (\ref{decompositionisomar}) qui est rendu explicite par l'équation (\ref{decompositionmatricemar}). \qed
\end{pre}
\end{Lem}

\begin{Cor} Le module à droite $O$ est le sous-module de $\Q_\Q \oplus \Q_\Q$ obtenu comme l'image de la multiplication à gauche des éléments de $\Q \oplus \Q$ (placés en colonnes) par la matrice $(\Id_2 - \Omega) \in M_2(\Q)$. 
\end{Cor}

\begin{Rem} \label{pi} Le corollaire précédent décrit $O$ avec une meilleure précision. Le rôle que doit jouer la matrice $(\Id_2 - \Omega) \in M_2(Q)$ dans une description plus explicite de $\MM$ est donc important. 

Dans ce qui suit, la notation $\pi = (\Id_2 - \Omega) \in M_2(Q)$ prévaudra. Un calcul donne
$$ \pi \; = \; \frac{1}{2}\matt{(1-T_1)(1+T_2)}{1-t_2^2}{1-t_1^2}{(1-T_2)(1+T_1)} \; \in M_2(Q) $$
\end{Rem}

\begin{Rem} \label{piomega} Reprenons en gardant l'idée en tête que les matrices $\pi$ et $\Omega$ sont des applications ($R$-linéaires) de $\Q \oplus \Q$ dans soi-même. En fait, puisqu'on voit ici $\Q \oplus \Q$ comme un $\Q$-module à droite, les multipications à gauche par $\pi$ ou $\Omega$ sont des homomorphismes de $\Q$-modules à droite. Autrement dit, on peut écrire 
$$ \pi,\, \Omega \, \in \End_\Q(\Q_\Q \oplus \Q_\Q) $$ 

En particulier, ces matrices idempotentes représentent des « projections » de $\Q \oplus \Q$ sur les sous-modules $O$ et $\Q = \im (\mathbf{sc})$, respectivement. En effet, c'est facile de vérifier que 
\begin{align*}
\pi |_{\im \pi} = \Id : \im \pi \rightarrow \im \pi && \Omega |_{\im \Omega} = \Id : \im \Omega \rightarrow \im \Omega
\end{align*}
\end{Rem}

D'après l'exercice 26 du chapitre 0 de \cite{FDnoncom}, on comprend qu'il y a un isomorphisme d'anneaux $\End_\Q(\Q_\Q) \cong \Q$. Il faut savoir que dans notre cas, $\Q = \Q_\Q$ est vu comme un module à droite. En fait, c'est surtout que tout endomorphisme $f \in \End_\Q(\Q_\Q)$ est entièrement déterminé par $f(1) \in \Q$ et si $f, g \in \End_\Q(\Q_\Q)$, alors 
$$ \big(f \circ g\big) (1) = f\big(g(1)\big) = f\big(1 \cdot g(1)\big) = f(1)g(1) $$ 
La multiplication dans $\Q$ correspond donc à la composition dans $\End_\Q(\Q)$. On a donc un isomorphisme naturel d'anneaux
\begin{align*}
\End_\Q(\Q_\Q) \; &\cong \quad \Q\\
f \quad &\leftrightarrow \; f(1) 
\end{align*}

Ensuite, puisque le foncteur $\Hom_\Q(-,-)$ est additif dans ses deux composantes par rapport aux sommes directes finies, on a un isomorphisme linéaire naturel
\begin{equation} \label{dude}
\End_\Q(\Q_\Q \oplus \Q_\Q) \cong \End_\Q(\Q_\Q)^{\oplus \, 4} \cong \Q^{\oplus \, 4}
\end{equation}

\begin{Prop} Il y a un isomorphisme naturel d'anneau
\begin{equation} \label{grossecorrespondanceend}
\End_\Q(\Q_\Q \oplus \Q_\Q) \cong M_2(\Q)
\end{equation}
où la composition des endomorphismes correspondent à la multiplication des matrices.

\begin{pre} L'isomorphisme linéaire naturel (\ref{dude}) donne que $\End_\Q(\Q_\Q \oplus \Q_\Q) \cong M_2(\Q)$ comme des $R$-modules. Il reste à montrer que la composition des endomorphismes correspond à la multiplication des matrices.

L'exercice 10 du chapitre 0 de \cite{FDnoncom} donne que chaque élément de $f \in \End_\Q(\Q_\Q \oplus \Q_\Q)$ s'identifie à une matrice $M_f \in M_2(\Q)$ dont les entrées sont issues de l'anneau $\Q \cong \End_\Q(\Q_\Q) $. Cet exercice donne aussi que sous cette identification, on a l'application de l'endomorphisme $f \in \End_\Q(\Q_\Q \oplus \Q_\Q)$ à $(a,b) \in \Q \oplus \Q$ correspond à la multiplication du vecteur colonne $(a,b)$ par la matrice $M_f$. En bref, c'est que
$$ a,b \in \Q \qquad \Longrightarrow \qquad f(a,b) \; = \; M_f \Vv{a}{b} \; \in \Q \oplus \Q $$
\qed
\end{pre}
\end{Prop}

\begin{Rem} Cette proposition permet de mieux appréhender les conclusions de la remarque \ref{piomega}, à savoir que les matrices idempotentes $\pi$ et $\Omega$ pouvaient être vues comme des éléments de $\End_\Q(\Q_\Q \oplus \Q_\Q)$.
\end{Rem}

Étant donné que $O \oplus \Q_\Q \cong \ker F \oplus \im (\mathbf{sc}) \cong \Q_\Q \oplus \Q_\Q$, on a un autre isomorphisme linéaire naturel
\begin{equation}
\End_\Q(\Q_\Q \oplus \Q_\Q) \cong \End_\Q(O) \oplus \Hom_\Q(O,\Q_\Q) \oplus \Hom_\Q(\Q_\Q,O) \oplus \End_\Q(\Q_\Q)
\end{equation}
Selon l'isomorphisme d'anneau (\ref{grossecorrespondanceend}), la remarque précédente et la remarque \ref{piomega}, c'est facile de voir que comme anneaux, on a un isomorphisme naturel
\begin{equation}
\End_\Q(O) \, \cong \; \pi \, M_2(\Q) \, \pi \; = \, \left\{ \pi \, A \, \pi \; | \; A \in M_2(\Q)\right\}
\end{equation}
où la composition des endomorphismes correspondent à la multiplication des matrices.

\begin{Rem} D'après l'exercice 27 du chapitre 0 de \cite{FDnoncom}, l'ensemble $\pi M_2(\Q) \pi$ est bel et bien un anneau unitaire. 
\end{Rem}

\begin{Cor} \label{descriptionultimemar} L'algèbre de Margaux $\MM$ définie par l'équation (\ref{premiermargaux}) s'écrit de façon un peu plus explicite comme
$$ \MM = \big(\lie(\pi M_2(\Q) \pi)\big)' $$
où $\,\pi \in M_2(\Q)$ est définie dans la remarque (\ref{pi}).
\end{Cor}


\begin{Rem} Le véritable but des efforts mis pour l'obtention du précédent corollaire était d'avoir un point de départ pour tenter de construire, un isomorphisme explicite 
\begin{equation} \label{isocheri}
\MM \otimes_R S \; \cong \; \mathfrak{sl}_2(S)
\end{equation}
Expliciter un tel isomorphisme permettrait de calculer précisément le 1-cocycle de \linebreak $\HHH^1\big(\Gamma,\Aut_{S-\lie}(\mathfrak{sl}_2(S)\big)$ qui correspond à la classe d'isomorphismes de $\MM$ en tant que forme tordue de $\mathfrak{sl}_2(R)$.

Il s'est trouvé que la complexité des calculs dans un anneau tel que $M_2(\Q)$ rendaient les choses assez difficiles. Mes essais pour expliciter des éléments arbitraires de $\MM$ avec le logiciel de calcul « Sage » n'ont abouti à rien, mais il faut dire que des problèmes de simplifications de calculs ont freiné tout espoir d'arriver à des expressions de taille raisonnable. À moins de pouvoir deviner un bon isomorphisme ou d'avoir une toute autre approche, il semble que l'obtention du cocycle correspondant à l'algèbre de Margaux $\MM$ n'est pas pour demain.
\end{Rem}

Pour parvenir à calculer un isomorphisme qui réalise (\ref{isocheri}), donc à calculer le 1-cocycle associé à $\MM$, il faudrait d'abord trouver une stratégie pertinente. J'espérais pouvoir faire calculer un élément arbitraire de $\MM$ moyennant sa description donnée au corollaire \ref{descriptionultimemar} et y reconnaître une correspondance avec une base de $\mathfrak{sl}_2(S)$. Cette non-stratégie aura au moins fonctionnée pour un exemple plus simple de forme tordue qui fait l'objet de l'annexe A.

\subsection{Blocs d'extensions de $\MM$-$\mathbf{mod}$, conjectures}

Puisque $\MM$ est une forme tordue avec $\g = \mathfrak{sl}_2(\CC)$, toutes les particularités de la section 3.1.2 s'appliquent à la description des blocs d'extension de la catégorie $\MM$-$\mathbf{mod}$ des $\MM$-modules de dimension finie.
 
En particulier, les blocs et blocs d'extensions de $\MM$-$\mathbf{mod}$ sont en bijection avec les caractères spectraux $\carspec^\Gamma$ de $\MM$. Rappelons que, dans ce contexte, ceux-ci sont des applications 
$$ \Max \left(\CC[t_1^{\pm 1},t_2^{\pm 1}]\right) = \{M_{a,b} \; | \; a,b \in \CC^\times \} \;\;\; \longrightarrow \bigsqcup_{M \in \Max S} \Bloc(\MM \, | \, \Gamma.M) $$ 
qui sont $\Gamma$-invariante sous une action de $\Gamma$ et dont le support est fini.

\begin{Rem} Rappelons que $\Bloc(\MM \, | \, \Gamma.M)$ désigne l'ensemble des blocs de la catégorie des $\MM$-modules de dimension finie dont tous les facteurs sont des modules d'évaluation avec support égal à $\Gamma.M$.
\end{Rem}

\begin{Rem} Les caractères spectraux sont déterminés par leur valeur sur un représentant de chaque $\Gamma$-orbite de $\Max S$. Pour déterminer ceux-ci, il doit donc suffire de spécifier l'image des idéaux 
$$ \big\{M_{a,b} \; | \; \Imagine a > 0 \text{ ou } a \in \RR_{> 0} \; \text{ et puis } \; \Imagine b > 0 \text{ ou } b \in \RR_{> 0} \big\} \subseteq \Max S $$
afin de décrire tous les caractères spectraux de l'algèbre de Margaux $\MM$.
\end{Rem}

Il faut maintenant pouvoir décrire les ensembles $\Bloc(\MM \, | \, \Gamma.M)$. Soit $P = \ZZ$, le groupe abélien des poids maximaux intégraux de $\mathfrak{sl}_2(\CC)$. Fixons $M \in \Max S$ et soit $P_+^M$ une copie de l'ensemble $P_+ = \NN$ des poids dominant intégraux de $\mathfrak{sl}_2(\CC)$. À la toute fin de l'article \cite{LPtfa}, on explique que les $\MM$-modules simples de dimension finie sont en bijection avec l'ensemble des applications 
$$ f : \; \{M_{a,b} \; | \; a,b \in \CC^\times \} \;\;\; \longrightarrow \bigsqcup_{M \in \Max S} P_+^M = P_+ $$ 
qui ont un support fini et qui sont telles que $f(M_{\pm a,\pm b})$ est constant dans $P_+$.

Puis, la remarque 5.12 de \cite{NSext} fait remarquer qu'il suffirait de trouver la relation d'équivalence à considérer sur l'ensemble $P_+$ qui correspond à partitionner les classes d'isomorphismes de $\MM$-modules simples d'évaluation supportées sur une même orbite selon ce qu'ils sont dans un même bloc d'extensions de $\MM$. La remarque décrit aussi un cadre d'hypothèses sous lesquelles cette relation d'équivalence sur $P_+$ est celle donnée par les co-ensembles de ses éléments dans le groupe additif $P/Q$, où $Q = 2\ZZ$ est le sous-groupe des racines de $\mathfrak{sl}_2(\CC)$.

\begin{Conj} La relation d'équivalence sur $P_+$ qui le partitionne selon ce que les classes d'isomorphismes des modules simples correspondantes sont dans un même bloc d'extension de $\MM$ est celle donnée par les co-ensembles de $P/Q$.
\end{Conj}

\begin{Conj} Les blocs de l'algèbre $\MM$ sont en bijection avec les sous-ensembles de cardinalité finie de l'ensemble suivant
$$ \big\{(a,b) \in (\CC^\times)^2 \; | \; \Imagine a > 0 \text{ ou } a \in \RR_{> 0} \; \text{ et puis } \; \Imagine b > 0 \text{ ou } b \in \RR_{> 0} \big\} \times \{0,1\} $$
\begin{pre} À partir de la conjecture précédente, il suffit de remarquer que dans le cas de $\mathfrak{sl}_2(\CC)$, nous avons 
$$ P/Q \; \cong \; \ZZ_2 \qquad \text{en tant que groupes} $$
Enfin, le premier ensemble n'est qu'un ensemble complet et non-redondant de représentants pour les $\Gamma$-orbites de $\Max \left(\CC[t_1^{\pm 1},t_2^{\pm 1}]\right)$. \flushright \textsc{Fin de la preuve à partir de la conjecture précédente}
\end{pre}
\end{Conj}

Voici maintenant quelques éléments d'indices qui portent à croire à une telle chose :

\begin{itemize}
\item[\textbullet] La proposition 4.9 de \cite{NSext} est, je crois, est encore valable dans ce contexte. Pour l'algèbre de Margaux $\MM$, il est facile de vérifier que l'action de $\Gamma \curvearrowright \Max S$ est libre. 

Malheureusement, les preuves de la section 4 de l'article \cite{NSext} dépendent de la spécificité de l'action de $\Gamma$ sur une algèbre d'applications équivariantes. Ceci permet aux auteurs d'obtenir les isomorphismes (4.7) qui ne sont plus valables dans le contexte des formes tordues. Cependant, le fait que l'action de $\Gamma \curvearrowright \Max \CC[t_1^{\pm 1},t_2^{\pm 1}]$ est libre dans le cas de $\MM$, laisse croire que le résultat tienne.
\item[\textbullet] En référence directe à la notation du lemme 4.5 de \cite{NSext},  un isomorphisme
$$ \big(\g \otimes (M/M^2)\big)^\Gamma \cong K_M/D_M \cong \g^M \otimes \big((M/M^2)^\Gamma\big) $$
est envisageable dans le cas de l'algèbre de Margaux $\MM$. Si l'existence d'un tel isomorphisme était avérée, on aurait directement accès au résultat de la proposition 4.7 et au corollaire 4.8 (toujours de l'article \cite{NSext}).

Un autre indice qui semble appuyer le résultat du corollaire 4.8 est la construction donnée à la proposition 3.4 de \cite{CMspec}. Un début de ce que je crois pouvoir être une preuve valable de cette construction est présentée dans l'annexe B.
\item[\textbullet] Voici quelques faits concernant les modules simples de dimension finie de $\mathfrak{sl}_2(\CC)$ :
\begin{itemize}
\item[$\rightarrow$] Ils sont déterminés à isomorphismes près par leur dimension qui correspond à leur poids maximal additionné de 1. 
\item[$\rightarrow$] Si $V(d)$ est un $\mathfrak{sl}_2(\CC)$-module simple de dimension $d \in \NN \bs \{0\}$, alors $\big(V(d)\big)^* \cong V(d)$. 
\item[$\rightarrow$] Si $d, e \in \NN \bs \{0\}$ et $d \geq e$, alors $V(d) \otimes V(e) \cong \bigoplus_{i=0}^e V(d-e+2i)$. 

Il s'agit de la formule de Clebsch-Gordan ; on en fait mention par exemple dans \cite{Humphreys} à l'exercice 7 de la page 126.
\end{itemize}
\item[\textbullet] Une généralisation de la proposition 1.2 de \cite{CMspec} tient. Il s'agit du corollaire A.4 de \cite{NSext} qu'on peut utiliser en fixant $U = \mathfrak{sl}_2(\CC)$.
\item[\textbullet] Supposant que la conjecture B.2 de ce mémoire tienne. La formule de Clebsch-Gordan permet alors de montrer que tous les modules d'évaluation de $\MM$ supportés sur une même orbite qui ont des dimensions de même parité sont dans un même bloc d'extension de $\MM$.

Finalement, puisque tous les $\LL$-modules d'évaluation supportés sur $\Gamma.M$ ne sont probablement pas tous dans un même bloc d'extensions, il doit y avoir deux de ces blocs d'extensions. 
\end{itemize}

Il reste donc un minimum de travail à accomplir pour avoir la descritption des blocs d'extensions de la catégorie $\MM$-$\mathbf{mod}$.

\chapter*{Conclusion}         
\phantomsection\addcontentsline{toc}{chapter}{Conclusion} 

Ce travail montre que l'étude des blocs d'extensions des catégories de modules de dimension finie pour les algèbres de courants tordues peut être faite selon la même approche que l'analogue pour le cas plus particulier des algèbres d'applications équivariantes. En effet, les blocs d'extensions des algèbres de courants tordues peuvent aussi être donnés en termes des blocs dans les catégories de modules d'évaluations appropriées et d'une pièce d'information concernant les modules de dimension 1. Les résultats de classification de la fin du chapitre 2 sont toutefois beaucoup plus esthétiques lorsque spécifiés au cas des formes tordues comme décrit dans la section 3.1. Ceci tient du fait que ce type d'algèbres de Lie de dimension infinie n'a pas de modules non-trivial de dimension 1. Presque bizzarement, les modules de dimension 1 représentent toujours la principale source de difficulté dans cette étude.

Un des défis importants restera celui de comprendre si la propriété de localité vis-à-vis des extensions, exigée des algèbres dès les travaux de E.Neher et A.Savage, peut être vérifiée pour toutes les algèbres de courants tordues ou non. En outre, cette description des blocs d'extensions dépend directement de la forme de la classification des objets simples dans la catégorie. Il y a place à amélioration quand à la mise en phase de la classification des objets simples donnée dans \cite{LauTCA} et la description des blocs donnée dans ce mémoire. Pour l'instant, la description des blocs d'extensions de ce mémoire reste en phase avec une description analogue à celle des classes d'isomorphismes d'objets simples de \cite{NSSreps} dans le cas des algèbres d'application équivariantes ; description moins attrayante puisqu'elle implique un choix. 

Malgré cela, un point théorique est tout de même amélioré quant à la détermination des extensions entre des modules simples d'évaluation arbitraires,. Soit $\LL = (\g \otimes S)^\Gamma$ une algèbre de courants tordue ; la proposition \ref{madeuxiemepropyo} révèle qu'étudier les modules de la forme $N_M/K_M'$ pour chaque orbite d'idéal maximal $\Gamma.M \subseteq \Max S$ est suffisant. La remarque \ref{aubordduprecipice} fait mention d'une vague conjecture, mais qui illustre ce propos. Pour une orbite $\Gamma.M \subseteq \Max S$ arbitraire, le $\g^M$-module de dimension finie bien spécifique $N_M/K_M'$ reste cependant très difficile à appréhender et ce, même dans des cas concrets et les plus simples possibles.

La catégorie $\LL$-$\mathbf{mod}$ des $\LL$-modules de dimension finie est abélienne alors les notions de groupe de Grothendieck ont un sens pour cette catégorie. Si $\LL$ est une algèbre de courants tordue, il serait bien intéressant de pouvoir relier d'une façon ou d'une autre le groupe de Grothendieck de $\LL$-$\mathbf{mod}$ à ceux des catégories $\g^M$-$\mathbf{mod}$ où $M \in \Max S$. Ce groupe est engendré par les classes d'isomorphismes des objets simples et est soumis aux relations $[E] = [V_1]+[V_2]$ pour toute suite exacte courte scindée $0 \rightarrow V_2 \rightarrow E \rightarrow V_1 \rightarrow 0$ dans la catégorie. Les modules de l'algèbre $\LL$ restent donc encore bien incompris, mais ce mémoire fournit un premier pas vers une meilleure compréhension des représentations de dimension finie de cette algèbre de Lie de dimension infinie.


\appendix                       
\chapter{Calcul de cocycle pour une forme tordue de $\mathfrak{sl}_2(\CC) \otimes \CC[t_1^{\pm 2},t_2^{\pm 2}]$}     

La théorie sur les formes tordues explique qu'il y a une bijection entre les classes d'isomorphismes des formes tordues d'une algèbre préalablement fixée et un ensemble de cohomologie galoisienne. Voir à ce sujet \cite{Bgalois} aux pages 118 et 119. Dans cette annexe, se trouve le calcul de la classe de cohomologie correspondante à  la classe d'isomorphismes une forme tordue d'une algèbre de Lie dont il a été question dans la section 3.2. 

Il s'agit du cas de l'algèbre obtenue en prenant $A = \Q$ dans (\ref{processussl2}). J'appellerai cette algèbre $\LL_1$ comme l'ont fait les auteurs P.Gille, A.Pianzola et M.Lau dans leurs articles \cite{GPmar} et \cite{LPtfa}.

La notation utilisée dans la section 3.2 est reprise et complémentée quelque peu dans cette annexe. En voici les éléments importants :
\begin{itemize}
\item[\textbullet] $k = \CC$
\item[\textbullet] $\zeta_4 \in \CC$ est un choix de racine quatri\`eme primitive de l'unit\'e.
\item[\textbullet] $\g = \mathfrak{sl}_2(\CC)$
\item[\textbullet] $S = \CC[t_1^{\pm 1},t_2^{\pm 1}]$
\item[\textbullet] $\Gamma = \ZZ_2 \times \ZZ_2$ et l'action de $(\ov{a},\ov{b}) \in \Gamma$ sur $S$ est donnée par extension $\CC$-linéaire de l'application de $S$ dans $S$ qui a pour effet
\begin{align*} 
t_1 &\mapsto (-1)^a \cdot t_1\\
t_2 &\mapsto (-1)^b \cdot t_2
\end{align*}
\item[\textbullet] $R = S^\Gamma = \CC[t_1^{\pm 2},t_2^{\pm 2}]$
\item[\textbullet] $\Q$ est la $R$-algèbre de quaternions généralisés de la définition \ref{quatgen}.
\item[\textbullet] $\LL_1 = \lie(\Q)'$ est une $R$-algèbre de Lie. 
\end{itemize}

\begin{annexeRem} $\LL_1$ est en fait une $S/R$-forme tordue de $\mathfrak{sl}_2(R) \cong \mathfrak{sl}_2(\CC) \otimes \CC[t_1^{\pm 2},t_2^{\pm 2}]$.

Ce qui suit sont les calculs qui mènent à l'obtention d'un 1-cocyle dont la classe de cohomologie dans $\HHH^1\Big(\Gamma,\Aut_{S-\lie}\big(\mathfrak{sl}_2(S)\big)\Big)$ correspond à la classe d'isomorphismes de $\LL_1$ en tant que \linebreak $R$-algèbre.
\end{annexeRem}

\begin{annexeRem} Tout le contenu de la remarque \ref{quaterniongen} est pertinent pour ce qui suit. En bref, $\Q$ est une $R$-alg\`ebre libre de rang 4, générée par $T_1$ et $T_2$ avec les relations $T_1^2 = t_1^2 \in R$, $T_2^2 = t_2^2 \in R$ et $T_1T_2 = -T_2T_1 \in \Q$. J'écrirai 
$$ \Q = R.1+R.T_1+R.T_2+R.T_1T_2 $$
\end{annexeRem}

\begin{annexeRem} Suivant ce qui est prescrit dans \cite{Bgalois}, le permier objectif pour mener à bien le calcul est de trouver un isomorphisme de $S$-algèbre de Lie entre les algèbres 
$$ \LL_1 \otimes_R S = \lie(\Q)'  \otimes_R S \qquad \qquad et \qquad \qquad \mathfrak{sl}_2(S) \cong \mathfrak{sl}_2(R) \otimes_R S $$ 

La théorie sur les formes tordues assure que si $\varphi : \LL_1 \otimes S \rightarrow \mathfrak{sl}_2(S)$ était un isomorphisme de $S$-algèbres de Lie, alors le cocycle défini par 
\begin{align} 
\alpha^{[\LL_1]} : \quad &\Gamma \;\;\;\;\; \longrightarrow \;\;\;\;\; \Aut_{S-\lie}\big(\mathfrak{sl}_2(S)\big) \notag \\ 
&\gamma \;\; \longmapsto \;\; \alpha_\gamma^{[\LL_1]} \, = \, \varphi \circ (\gamma \cdot \varphi^{-1}) \tag{A.4}
\end{align}
\addtocounter{laid}{1}
sera un cocycle dont la classe de cohomologie correspond à la classe d'isomorphismes de la $R$-algèbre de Lie $\LL_1$ en tant que $S/R$-forme tordue de $\mathfrak{sl}_2(R)$.
\end{annexeRem}

Ainsi, il faut d'abord trouver un isomorphisme $\varphi$. Pour en trouver un, il faut avoir une idée de générateurs de la $R$-algèbre de Lie $\LL_1$, puis comprendre les crochets de Lie entre ces générateurs, ou autrement dit, comprendre plus explicitement $\LL_1$.

\begin{annexeLem} $\LL_1 = R.T_1 + R .T_2 + R.T_1T_2$. 

En particulier, $\LL_1 \otimes_R S = S.T_1 + S.T_2 + S.T_1T_2$
\begin{pre} Un élément arbitraire de $\LL_1 = \lie(Q)'$ en est un de la forme
$$ [a_0 + a_1T_1 + a_2T_2 + a_3T_1T_2 \, , \, b_0 + b_1T_1 + b_2T_2 + b_3T_1T_2] $$
où les $a_i$ et $b_i$ sont des éléments de $R$. Cet élément arbitraire est égal à...
\begin{align*}
& \;\; a_0b_0 + a_0b_1T_1 + a_0b_2T_2 + a_0b_3T_1T_2 + a_1b_0T_1 + a_1b_1T_1^2 + a_1b_2T_1T_2 + a_1b_3T_1^2T_2\\
&+ a_2b_0T_2 - a_2b_1T_1T_2 + a_2b_2T_2^2 -a_2b_3T_1T_2^2 + a_3b_0T_1T_2 - a_3b_1T_1^2T_2 +a_3b_2T_1T_2^2\\ &- a_3b_3T_1^2T_2^2 - \big(\, a_0b_0 + a_1b_0T_1 + a_2b_0T_2 + a_3b_0T_1T_2 +a_0b_1T_1 + a_1b_1T_1^2 + a_2b_1T_1T_2\\
& + a_3b_1T_1^2T_2 + a_0b_2T_2 - a_1b_2T_1T_2 + a_2b_2T_2^2 -a_3b_2T_1T_2^2 + a_0b_3T_1T_2 - a_1b_3T_1^2T_2\\
&+ a_2b_3T_1T_2^2 - a_3b_3T_1^2T_2^2 \, \big)\\
= \; & \;\; 2a_1b_2T_1T_2 -2a_2b_1T_1T_2 + 2a_1b_3T_1^2T_2 - 2a_2b_3T_1T_2^2 - 2a_3b_1T_1T_2^2 + 2a_3b_2T_1T_2^2\\
= \; & \;\; 2(a_1b_2-a_2b_1)T_1T_2 + 2(a_1b_3 - a_3b_1)T_1^2T_2 + 2(a_3b_2 - a_2b_3)T_1T_2^2\\
= \; & \;\; 2t_2^2(a_3b_2 - a_2b_3)T_1 + 2t_1^2(a_1b_3 - a_3b_1)T_2 + 2(a_1b_2-a_2b_1)T_1T_2
\end{align*}
Il faut donc spécifier 6 paramètres dans $R$ pour déterminer un élément de $\LL_1$. Maintenant, on peut directement vérifier que...
\begin{itemize}
\item[\textbullet] $a_1 = \frac{1}{\sqrt{2}}, \; a_2 = 0, \; a_3 = 0, \; b_1 = 0, \; b_2 = \frac{1}{\sqrt{2}}, \; b_3 = 0 \;$ donnent que $\; T_1T_2 \in \LL_1$.
\item[\textbullet] $a_1 = \frac{1}{\sqrt{2}}t_1^{-2}, \; a_2 = 0, \; a_3 = 0, \; b_1 = 0, \; b_2 = 0, \; b_3 = \frac{1}{\sqrt{2}} \;$ donnent que $\; T_2 \in \LL_1$.
\item[\textbullet] $a_1 = 0, \; a_2 = 0, \; a_3 = \frac{1}{\sqrt{2}}t_2^{-2}, \; b_1 = 0, \; b_2 = \frac{1}{\sqrt{2}}, \; b_3 = 0 \;$ donnent que $\; T_1 \in \LL_1$.
\end{itemize}

Ceci entraîne directement que les éléments arbitraires de la $R$-algèbre de Lie $\LL_1$ sont de la forme $c_1T_1 + c_2T_2 + c_3T_1T_2$ où $c_1,c_2,c_3 \in R = \CC[t_1^{\pm 2},t_2^{\pm 2}]$.
\qed
\end{pre}
\end{annexeLem}

Il faut maintenant calculer les crochets de Lie entre les paires issues de ces trois générateurs de $\LL_1$. Voici les calculs : 
\begin{align}
[T_1,T_2] &= T_1T_2 - T_2T_1 = 2T_1T_2 \tag{A.6}\\
[T_1T_2,T_2] &= T_1T_2^2-T_2T_1T_2 = -2T_1T_2^2 = 2t_2^2T_1 \tag{A.7}\\
[T_1T_2,T_1]&= T_1T_2T_1 - T_1^2T_2 = -2T_1^2T_2 = -2t_1^2T_2 \tag{A.8}
\end{align}
\addtocounter{laid}{3}

On voit dans les relations ci-hautes que le crochet de $T_1T_2$ avec un des $T_i$ a pour effet de le changer en une l'autre générateur, à une constante près. C'est aussi manifeste que le crochet de des deux $T_i$ entre eux donnent un multiple scalaire de $T_1T_2$.

Des observations similaires peuvent être faites avec $\mathfrak{sl}_2(S) = S.h + S.e + S.f$ puisque :
\begin{align*}
[e - f,e + f] &= h - (-h) = 2h\\
[h,e - f] &= 2e-(-2f) = 2(e + f)\\
[h,e + f] &= 2e - 2f = 2(e - f)
\end{align*}

Puisque $\{h, e - f, e + f\}$ est aussi une base de $\mathfrak{sl}_2(S)$ comme $S$-algèbre de Lie, on peut penser à faire une association 
\begin{align*}
T_1T_2 \quad &\longleftrightarrow \quad \omega h\\
T_1 \quad &\longleftrightarrow \quad \lambda (e - f)\\
T_2 \quad &\longleftrightarrow \quad \mu (e + f)
\end{align*}
où $\omega, \lambda, \mu \in S$.

Pour réussir à trouver l'isomorphisme voulu, il suffira donc de trouver un triplet $(\omega,\lambda,\mu)$ dans $S \times S \times S$ qui rend compte des mêmes constantes de structures que dans les équations (A.6), (A.7) et (A.8). Il faut donc solutionner...
\begin{align}
2 \omega h &= [\lambda(e - f), \mu(e + f)] = \lambda \mu [e - f, e + f] = 2 \lambda \mu h \tag{A.9}\\
2t_2^2 \lambda (e - f) &= [\omega h, \mu (e+f)] = \omega \mu [h, e + f]  = 2 \omega \mu (e - f) \tag{A.10}\\
-2t_1^2 \mu (e+f) &= [\omega h, \mu (e - f)] = \omega \lambda [h, e - f]  = 2 \omega \lambda (e + f) \tag{A.11}
\end{align}
\addtocounter{laid}{3}

De (A.9), on obtient $\omega = \lambda \mu$. En substituant dans ce qu'on obtient des deux autres équations, on trouve
\begin{align*}
2t_2^2 \lambda &= 2 \lambda \mu^2 \qquad \qquad \qquad \qquad -2t_1^2\mu = \lambda^2 \mu\\
t_2^2 &= \mu^2 \qquad \qquad \qquad \qquad \qquad \,-t_1^2 = \lambda^2
\end{align*}
Si donc, on pose $\lambda = \zeta_4t_1 \in S$ et $\mu = t_2 \in S$, ces équations seront respectées. Ceci implique que $\omega = \zeta_4t_1t_2 \in S$. 

Par construction, c'est donc que l'association 
\begin{align*}
T_1T_2 \quad &\longleftrightarrow \quad \zeta_4 h\\
T_1 \quad &\longleftrightarrow \quad \zeta_4 t_1 (e - f)\\
T_2 \quad &\longleftrightarrow \quad t_2 (e + f)
\end{align*}
respecte le crochet de Lie des deux côtés. De plus, puisque c'est une bijection $S$-linéaire entre des $S$-bases des algèbres $\LL_1 \otimes_R S$ et $\mathfrak{sl}_2(S)$, elle donne lieu à un isomorphisme de $S$-algèbres de Lie. L'application 
\begin{align*}
\varphi :  &\qquad \LL_1 \otimes_R S \quad \longrightarrow \quad \mathfrak{sl}_2(S) \\
& \qquad \quad \,\,\,\, T_1 \,\,\,\,\,\, \longmapsto \quad \zeta_4t_1\,(e-f) \\
& \qquad \quad \,\,\,\, T_2 \,\,\,\,\,\, \longmapsto \quad t_2\,(e+f) \\
& \qquad \quad \,\, T_1T_2 \,\, \longmapsto \quad \zeta_4t_1t_2\,h \\
\end{align*}
est donc un $S$-isomorphisme d'algèbres de Lie qui fera l'affaire pour le calcul du cocycle.

\begin{annexeRem} Tout est maintenant en place pour le calcul direct du cocycle $\alpha^{[\LL_1]}$ puisqu'un isomorphisme $\varphi$ a été obtenu. Voir l'équation (A.4) au besoin. 
\end{annexeRem}

\begin{annexeRem} Le groupe $\Gamma = \ZZ_2 \times \ZZ_2$ agit sur les deux $S$-algèbres de Lie $\mathfrak{sl}_2(S)$ et $\LL_1 \otimes_R S$ via son action sur $S$ ; c'est-à-dire via son action sur les coefficients. 
\end{annexeRem}

\begin{annexeLem} 
Fixons $(\bar{a},\bar{b}) \in \Gamma = \ZZ_2 \times \ZZ_2$. 

Toujours selon (A.4), l'automorphisme de $S$-algèbre de Lie $\alpha_{(\bar{a},\bar{b})}^{[\LL_1]}$ de $\mathfrak{sl}_2(S)$ exprimé en termes d'une matrice dans la base ordonnée $\big(\, \zeta_4t_1t_2\,h \, ,\, \zeta_4t_1\,(e - f) \, ,\, t_2\,(e+f) \, \big)$ est 
$$ \alpha_{(\bar{a},\bar{b})}^{[\LL_1]} \; = \;
\begin{pmatrix}
(-1)^{a+b} & 0 & 0 \\ 
0 & (-1)^{a} & 0 \\
0 & 0 & (-1)^{b}
\end{pmatrix} \; :  \quad \mathfrak{sl}_2(S) \rightarrow \mathfrak{sl}_2(S) $$
\begin{pre} Il s'agit de voir comment $\alpha_{(\bar{a},\bar{b})}^{[\LL_1]} \in \Aut_{S-\lie}\big(\mathfrak{sl}_2(S)\big)$ transforme les éléments de la base $\big(\, \zeta_4t_1t_2\,h \, ,\, \zeta_4t_1\,(e - f) \, ,\, t_2\,(e+f) \, \big)$.

Rappelons d'abord que 
\begin{align*} 
\alpha_{(\bar{a},\bar{b})}^{[\LL_1]} &= \varphi \circ \big((\bar{a},\bar{b}) \cdot \varphi^{-1}\big) \\
&=\varphi \circ (\bar{a},\bar{b}). \circ \varphi^{-1} \circ (\bar{a},\bar{b}).
\end{align*} 

Cet automorphisme opère donc comme suit sur la base :
\begin{align*} 
&\quad \Vvv{\zeta_4t_1t_2\,h}{\zeta_4t_1\,(e-f)}{t_2\,(e+f)} \quad 
\stackrel{(\bar{a},\bar{b}).}{\longmapsto} \quad \Vvv{(-1)^{a+b} \, \zeta_4t_1t_2\,h}{(-1)^a \, \zeta_4t_1\,(e-f)}{(-1)^b \, t_2\,(e+f)} \quad 
\stackrel{\varphi^{-1}}{\longmapsto} \quad \Vvv{(-1)^{a+b} \, T_1T_2}{(-1)^a \, T_1}{(-1)^b \, T_2} \\
&\stackrel{(\bar{a},\bar{b}).}{\longmapsto} \quad \Vvv{(-1)^{a+b} \, T_1T_2}{(-1)^a \, T_1}{(-1)^b \, T_2} \quad 
\stackrel{\varphi}{\longmapsto} \quad \Vvv{(-1)^{a+b} \, \zeta_4t_1t_2\,h}{(-1)^a \, \zeta_4t_1\,(e-f)}{(-1)^b \, t_2\,(e+f)}
\end{align*}
\qed
\end{pre}
\end{annexeLem}
 
\begin{annexeCor} L'automorphisme de $S$-algèbre de Lie $\alpha_{(\bar{a},\bar{b})}^{[\LL_1]}$ de $\mathfrak{sl}_2(S)$ exprimé en termes d'une matrice dans la base ordonnée $(\, h \, ,\, e \, ,\, f \,)$ est
$$ \alpha_{(\bar{a},\bar{b})}^{[\LL_1]} \; = \; \begin{pmatrix}
(-1)^{a+b} & 0 & 0 \\ 
0 & \frac{(-1)^{a}+(-1)^{b}}{2} & (-\zeta_4)t_1^{-1}t_2 \frac{(-1)^{a}-(-1)^{b}}{2} \\
0 & \zeta_4t_1t_2^{-1} \frac{(-1)^{a}-(-1)^{b}}{2} & \frac{(-1)^{a}+(-1)^{b}}{2}
\end{pmatrix} \; : \; \mathfrak{sl}_2(S) \rightarrow \mathfrak{sl}_2(S) $$
\begin{pre} Il suffit de faire le changement de base approprié. La matrice $M \in GL\big(\mathfrak{sl}_2(S)\big)$ qui représente le changement de base $(\, h \, ,\, e \, ,\, f \,) \, \mapsto \, \big(\, \zeta_4t_1t_2\,h \, ,\, \zeta_4t_1\,(e - f) \, ,\, t_2\,(e+f) \, \big)$ et son inverse sont
\begin{align*}
M = \begin{pmatrix}
\zeta_4t_1t_2 & 0 & 0 \\ 
0 & \zeta_4t_1 & t_2 \\
0 & -\zeta_4t_1 & t_2
\end{pmatrix} 
&& 
M^{-1} = \begin{pmatrix}
(-\zeta_4)t_1^{-1}t_2^{-1} & 0 & 0 \\ 
0 & \frac{1}{2}(-\zeta_4)t_1^{-1} & \frac{1}{2}\zeta_4 t_1^{-1} \\
0 & \frac{1}{2}t_2^{-1} & \frac{1}{2}t_2^{-1}
\end{pmatrix}
\end{align*}
Il ne reste qu'à calculer le produit matriciel suivant :
\begin{align*} & \qquad \qquad \qquad \qquad \qquad
M^{-1} \; \begin{pmatrix}
(-1)^{a+b} & 0 & 0 \\ 
0 & (-1)^{a} & 0 \\
0 & 0 & (-1)^{b}
\end{pmatrix} \; M \\
&= \begin{pmatrix}
(-\zeta_4)t_1^{-1}t_2^{-1} & 0 & 0 \\ 
0 & \frac{1}{2}(-\zeta_4)t_1^{-1} & \frac{1}{2}\zeta_4 t_1^{-1} \\
0 & \frac{1}{2}t_2^{-1} & \frac{1}{2}t_2^{-1}
\end{pmatrix} \cdot 
\begin{pmatrix}
(-1)^{a+b} & 0 & 0 \\ 
0 & (-1)^{a} & 0 \\
0 & 0 & (-1)^{b}
\end{pmatrix} \cdot 
\begin{pmatrix}
\zeta_4t_1t_2 & 0 & 0 \\ 
0 & \zeta_4t_1 & t_2 \\
0 & -\zeta_4t_1 & t_2
\end{pmatrix} \\
&= \begin{pmatrix}
(-\zeta_4)t_1^{-1}t_2^{-1} & 0 & 0 \\ 
0 & \frac{1}{2}(-\zeta_4)t_1^{-1} & \frac{1}{2}\zeta_4 t_1^{-1} \\
0 & \frac{1}{2}t_2^{-1} & \frac{1}{2}t_2^{-1}
\end{pmatrix} \cdot 
\begin{pmatrix}
(-1)^{a+b} \, \zeta_4t_1t_2 & 0 & 0 \\ 
0 & (-1)^{a} \, \zeta_4t_1 &  (-1)^{a} \, t_2 \\
0 & -(-1)^{b} \, \zeta_4t_1 & (-1)^{b} \, t_2
\end{pmatrix} \\
&= \begin{pmatrix}
(-1)^{a+b} & 0 & 0 \\ 
0 & \frac{(-1)^{a}+(-1)^{b}}{2} & (-\zeta_4)t_1^{-1}t_2 \frac{(-1)^{a}-(-1)^{b}}{2} \\
0 & \zeta_4t_1t_2^{-1} \frac{(-1)^{a}-(-1)^{b}}{2} & \frac{(-1)^{a}+(-1)^{b}}{2}
\end{pmatrix}
\end{align*}
\qed
\end{pre}
\end{annexeCor}

Le cocycle $\alpha^{[\LL_1]}$ a été calculé avec succès ! Fixant la base $(\, h \, ,\, e \, ,\, f \,)$ de $\mathfrak{sl}_2(S)$, ce 1-cocycle peut se voir comme l'application qui envoie les éléments de $\Gamma$ à la forme matricielle de l'automorphisme correspondant dans la base $(\, h \, ,\, e \, ,\, f \,)$ ;

\begin{itemize}
\item[\textbullet] $(\bar{0},\bar{0})$ et $(\bar{1},\bar{1}) \quad$ sur la matrice$\quad \begin{pmatrix} 1 & 0 & 0 \\ 0 & 1 & 0 \\ 0 & 0 & 1 \end{pmatrix}$.
\item[\textbullet] $(\bar{1},\bar{0})$ et $(\bar{0},\bar{1}) \quad$ sur la matrice $\quad \begin{pmatrix} -1 & 0 & 0 \\ 0 & 0 & (-\zeta_4)t_1^{-1}t_2 \\ 0 & \zeta_4t_1t_2^{-1} & 0\end{pmatrix}$.
\end{itemize}

\chapter{Une idée pour les formes tordues}

Voici une petite conjecture et ce que je pense en être le début d'une démonstration valable. Le résultat serait une condition suffisante pour que deux modules simples d'évaluation supportés sur une seule orbite aient des extensions non-triviales dans le cas où l'algèbre de courants tordue $\LL$ vérifie la condition $\g^M = \g$ pour tout $M \in \Max S$. C'est le cas pour les formes tordues (voir \cite{LPtfa}, proposition 3.3).

Avant de montrer le résultat, il faut rappeler une notion simple de géométrie algébrique.

\begin{annexeDef} Soit $M \in \Max S$. Alors une \textbf{dérivation de $S$ au point $M$} est une application $k$-linéaire 
$$ \partial : \; S \longrightarrow k $$
qui satisfait à la propriété que pour tout $s_1,s_2 \in S$, on ait 
$$ \partial(s_1s_2) = \partial(s_1) \, s_2(M) + s_1(M) \, \partial(s_2) $$

L'ensemble des dérivations de $S$ au point $M$ est un sous-espace vectoriel de $\Hom_k(S,k)$.
\end{annexeDef}

L'ensemble des dérivations de $S$ au point $M \in \Max S$ est naturellement isomorphe comme espace vectoriel à $(M/M^2)^*$, l'espace tangent associé à $M$. La démonstration est la même que celle donnée dans \cite{GWsym} où on remplace $f \in \CC[t]$ par $s \in S$ et $a \in \CC$ par $M \in \Max S$.

La condition suffisante mentionnée un peu plus haut provient de la proposition qui suit. Elle est inspirée de la proposition 3.4 de l'article \cite{CMspec}.

\begin{annexeConj} \label{trucmargaux} Supposons que l'algèbre de courants tordue $\LL$ satisfasse la propriété que $\g^M = \g$ pour chaque $M \in \Max S$. Soient $V_1$ et $V_2$ des $\g^M$-modules simples de dimension finie. 

Posons $V = V_1^* \otimes V_2$. Alors, il y a une application linéaire injective
\begin{align*}
\beta : \; \Hom_{\g}(\g,V) \otimes (M/M^2)^* \; \longrightarrow \; \cohomo{\LL}{V}
\end{align*}

\em \textbf{Idée de preuve :} Identifions $(M/M^2)^*$ aux dérivations de $S$ au point $M$. Pour ce qui suit, un élément de $(M/M^2)^*$ sera donc une dérivation de $S$ au point $M$. 

L'application $\beta$ envoie un tenseur simple $f \otimes \partial \in \Hom_{\g}(\g,V) \otimes (M/M^2)^*$ sur la restriction à $\LL$ de l'application
\begin{align*}
\tilde{\beta}_{f \otimes \partial} : \quad \g &\otimes S \quad \longrightarrow \quad V\\
x &\otimes s \; \longmapsto \; \partial(s)f(x)
\end{align*}
Notons que cette définition fait de l'application $\beta$ une application $k$-linéaire bien définie. Il faut ensuite montrer que les fonctions linéaires qui composent l'image de $\beta$ sont bien des dérivations de $\LL$ dans $V$. 

Chaque $\tilde{\beta}_{f \otimes \partial}$ est $k$-linéaire par construction. Montrons ensuite que chaque $\tilde{\beta}_{f \otimes \partial}$ est une dérivation de $\g \otimes S$ dans $V$. Soient donc $x \otimes s$ et $y \otimes r \in \g \otimes S$. On peut alors écrire
\begin{align*}
\Big(\tilde{\beta}_{f \otimes \partial}\Big)\big([x \otimes s, y \otimes r]\big) &= \Big(\tilde{\beta}_{f \otimes \partial}\Big)\big([x,y] \otimes sr\big)\\
&= \partial(sr)f\big([x,y]\big)\\
&= \big(\partial(s) \, r(M) + s(M) \, \partial(r)\big)\,f\big([x,y]\big)\\
&= \partial(s) \, r(M)f\big([x,y]\big) + \partial(r) \, s(M)f\big([x,y]\big)\\
&= - \, \partial(s) \, r(M)f\big(y.x\big) + \partial(r) \, s(M)f\big(x.y\big)\\
&= - \, \partial(s) \, r(M)\,y.f\big(x\big) + \partial(r) \, s(M)\,x.f\big(y\big)\\
&= - \, \partial(s) \, (y \otimes r).f(x) + \partial(r) \, (x \otimes s).f(y)\\
&= - (y \otimes r).\big(\partial(s)f(x)\big) + (x \otimes s).\big(\partial(r)f(y)\big)\\
&= (x \otimes s).\big(\partial(r)f(y)\big) - (y \otimes r).\big(\partial(s)f(x)\big)\\
&= (x \otimes s).\Big(\big(\tilde{\beta}_{f \otimes \partial}\big)(y \otimes r)\Big) - (y \otimes r).\Big(\big(\tilde{\beta}_{f \otimes \partial}\big)(x \otimes s)\Big)
\end{align*}
Ainsi, $\tilde{\beta}_{f \otimes \partial} \in \Der(\g \otimes S,V)$ et en conséquence, sa restriction à $\LL \subseteq \g \otimes S$ appartient à l'ensemble de dérivations $\Der(\LL,V)$. En composant avec le passage au quotient $\Der(\LL,V) \rightarrow \cohomo{\LL}{V}$, on obtient bien une application entre les bons ensembles.

Pour terminer, il faudrait montrer que l'application $\beta$ est injective, mais \textbf{je n'ai pas réussi à le faire}. Cependant, je peux montrer que $\tilde{\beta} : \Hom_\g(\g,V) \otimes (M/M^2)^* \rightarrow \cohomo{\g \otimes S}{V}$ est injective. 

Fixons une base $\{\ov{m}_j\}_{j \in J}$ de $(M/M^2)$ et fixons $\{\partial_j\}_{j \in J}$ la base duale associée  à cette dernière. Notons que les $\partial_j$ sont interprétés ici comme des dérivations de l'algèbre $S$ au point $M$. Fixons également une base $\{f_i\}_{i \in I}$ de $\Hom_\g(\g,V)$.

Soit $\sum_{i,j} c_{ij} \, f_i \otimes \partial_j \in \ker \tilde{\beta}$ (la somme étant finie). Alors son image est une dérivation intérieure de $\g \otimes S$ dans $V$. Fixons donc $v \in V$ tel que 
$$ \sum_{i,j} c_{i,j} \, \tilde{\beta}_{f_i \otimes \partial_j} = d^v \in \IDer(\g \otimes S,V)$$ 
Alors, pour tout $x \in \g$, on doit avoir
\begin{align*}
\sum_{i,j} c_{ij} \, \Big(\tilde{\beta}_{f_i \otimes \partial_j}\Big)\big(x \otimes (\textstyle{\sum_k}\, m_k)\big) &= \big(x \otimes (\textstyle{\sum_k}\, m_k)\big).v \\
\sum_{i,j} c_{ij} \, \partial_j\big(\textstyle{\sum_k}\, m_k\big) f_i(x) &= \big(\textstyle{\sum_k}\, m_k\big)(M) \; x.v \\
\sum_{i,j} c_{ij} \, f_i(x) &= 0 \\
\sum_i \Big(\sum_j c_{ij}\Big) \, f_i(x) &= 0
\end{align*}
Puisque les $\{f_i\}_{i \in I}$ forment une base de $\Hom_\g(\g,V)$, on obtient que $\sum_j c_{i,j} = 0$ pour chacun des indices $i$. Alors on conclut que
$$ \sum_{i,j} c_{ij} \, f_i \otimes \partial_j = \sum_i \Big(\sum_j c_{ij}\Big) \, f_i \otimes \partial_j = 0 \in \Hom_\g(\g,V) \otimes (M/M^2)^*$$
Ainsi, $\tilde{\beta}$ est injective. 

Le problème avec cette méthode pour prouver l'injectivité de $\beta$ est que les éléments de la forme $x \otimes (\sum_k m_k)$, où $x \in \g$ est arbitraire, ne sont pas forcément dans $\LL$ alors on ne peut pas utiliser ces éléments pour aboutir à la conclusion analogue à celle ci-haut dans le cas de $\g \otimes S$. \flushright \textsc{Fin de l'idée de preuve.}
\end{annexeConj}

\begin{annexeRem} La preuve de l'injectivité de $\tilde{\beta}$ ci-haut généralise la construction de V.Chari et A.Moura donnée à la proposition 3.4 de leur article \cite{CMspec}. La différence est qu'ici, l'algèbre de Lie en question est $\g \otimes S$ au lieu de $\g \otimes k[t^{\pm 1}]$.
\end{annexeRem}

\begin{annexeConj} Soit $M \in \Max S$ et soient $V_1$ et $V_2$ des $\g^M$-modules simples de dimension finie. Posons $V = V_1^* \otimes V_2$. Alors
$$ \left.\begin{array}{c} \Hom_{\g^M}(\g^M,V) \neq \{0\} \\ (M/M^2)^* \neq \{0\}\end{array}\right] \quad \Longrightarrow \quad \Ext^1_\LL(V_1,V_2) \neq \{0\} $$
\end{annexeConj}

\begin{annexeRem} Ces deux conjectures donneraient un sens à la possible « interprétation géométrique » à laquelle E.Neher et A.Savage faisaient référence suite à leur lemme 4.5 dans \cite{NSext}. Voir la remarque 4.6 de ce même article.
\end{annexeRem}                

~\nocite{EWliealg}

\bibliography{Biblio}                 

\end{document}